\documentclass[3p,authoryear,11pt,fleqn,numbers,square,sort&compress]{elsarticle}
\usepackage{amssymb,amsbsy,amsmath,amsfonts,amssymb,amscd}
\usepackage[english]{babel}
\usepackage{color}
\usepackage{float}
\usepackage{amssymb,natbib}
\usepackage{graphicx}
\usepackage[latin1]{inputenc}
\usepackage{stmaryrd}
\usepackage{mathrsfs}
\usepackage{booktabs}
\usepackage{multirow,booktabs}
\usepackage{array}
\usepackage{subfig}
\usepackage{geometry}
\usepackage[ruled]{algorithm2e}
\usepackage[hidelinks,colorlinks=true,citecolor=blue,linkcolor=red,urlcolor=green,pdfstartview=]{hyperref}
\usepackage{pifont}
\usepackage{framed}
\usepackage{nomencl}
\usepackage{setspace}
\usepackage{etoolbox}

\makenomenclature
\renewcommand\nomgroup[1]{\item[\bfseries
\ifstrequal{#1}{P}{Peridynamics}{\ifstrequal{#1}{N}{Finite element}{\ifstrequal{#1}{O}{Other Symbols}{}}}]}
\biboptions{comma,round}
\newtheorem{e-proposition}[theorem]{Proposition}

\newtheorem{e-definition}[theorem]{Definition\rm}

\def\og{\leavevmode\raise.3ex\hbox{$\scriptscriptstyle\langle\!\langle$~}}
\def\fg{\leavevmode\raise.3ex\hbox{~$\!\scriptscriptstyle\,\rangle\!\rangle$}}

\journal{Computer Methods in Applied Mechanics and Engineering}
\input{tcilatex}

\graphicspath{{fig/}}

\begin{document}

\begin{frontmatter}



\title{Hybrid FEM and Peridynamic simulation of hydraulic fracture propagation \\ in saturated porous media}

\author[label1,label2,label3]{Tao Ni\corref{cor1}}
\ead{hhunitao@hhu.edu.cn}

\author[label5]{Francesco Pesavento}
\author[label3,label4]{Mirco Zaccariotto}
\author[label3,label4]{Ugo Galvanetto\corref{cor1}}
\ead{ugo.galvanetto@unipd.it}
\author[label1,label2]{Qi-Zhi Zhu}
\author[label5,label6]{Bernhard A. Schrefler}

\cortext[cor1]{Corresponding author}
\address[label1]{College of Civil and Transportation Engineering, Hohai University, 210098, Nanjing, China}
\address[label2]{Key Laboratory of Ministry	of Education for Geomechanics and Embankment Engineering, Hohai University, Nanjing, China}
\address[label3]{Industrial Engineering Department, University of Padova, via Venezia 1, Padova, 35131, Italy}
\address[label4]{Center of Studies and Activities for Space (CISAS)-G. Colombo, via Venezia 15, Padova, 35131, Italy}
\address[label5]{Department of Civil, Environmental and Architectural Engineering, University of Padova, via Marzolo 9, Padova, 35131, Italy}
\address[label6]{Institute for Advanced Study, Technische Universität München, Lichtenbergstrasse 2a, D-85748 Garching b. München, Germany}

\begin{abstract}
This paper presents a hybrid modelling approach for simulating hydraulic fracture propagation in saturated porous media: ordinary state-based peridynamics is used to describe the behavior of the solid phase, including the deformation and crack propagation, while FEM is used to describe the fluid flow and to evaluate the pore pressure. Classical Biot poroelasticity theory is adopted. The proposed approach is first verified by comparing its results with the exact solutions of two examples. Subsequently, a series of pressure- and fluid-driven crack propagation examples are solved and presented. The phenomenon of fluid pressure oscillation is observed in the fluid-driven crack propagation examples, which is consistent with previous experimental and numerical evidences. All the presented examples demonstrate the capability of the proposed approach in solving problems of hydraulic fracture propagation in saturated porous media.
\end{abstract}

\begin{keyword}
	Peridynamics \sep Hydraulic fracture propagation \sep Saturated porous media  \sep Finite element method 
\end{keyword}

\end{frontmatter}

\section{Introduction}

Hydraulic fracture (HF) is an effective technology to improve the recovery
rate of natural resources (oil, gas etc.) from reservoirs with
low-permeability \citep{zhou2019phase}. The cracks produced by hydraulic
fracturing connect the resources in the reservoirs with the wellbores, which
allows to extract a vast amount of resources inaccessible in the past and
provides great economic benefits \citep{zhou2018phase}. However, cracks may
also lead to leakage of the natural resources and to groundwater
contamination \citep{vidic2013impact,ren2017equivalent}. Furthermore, the
hydraulic fracturing mining process may also result in air and surface
pollution around the mining area \citep{mikelic2013phase}. Therefore,
accurate prediction of hydraulic fracture consequences in unconventional
hydrocarbon reservoirs has a great relevance in guiding the exploitation of
natural resources \citep{turner2013non,ouchi2015fully,ouchi2017effect}, as
well as in balancing the vast economic benefits with the environmental risks.

The numerical study of the HF propagation in deformable porous media has
become a major research topic in mechanical, energy and environmental
engineering during the past decades %
\citep{mikelic2013phase,lee2016pressure,zhou2018phase}. HF is a typical
hydro-mechanical coupled problem involving deformation and fracture
evolution of the solid phase as well as fluid flow in the fractured area. To
solve this type of coupled system, a variety of complex numerical models
have been developed \citep{lecampion2018numerical}, which can be simply
classified into three categories: continuous, discrete and hybrid
approaches. The typical continuous approaches mainly refer to the
computational techniques based on the finite element method. These
approaches are potentially appropriate for irregular body shapes with
heterogeneous material properties and non-linear behavior. In %
\citep{simoni2003cohesive,schrefler2006adaptive,secchi2007mesh}, the
adaptive re-meshing technique is used to follow the crack surfaces during
the hydraulic fracturing process. Re-meshing is usually computationally
expensive and it is difficult to develop robust schemes considering mass and
momenta conservation \citep{secchi2007mesh}. Cohesive zone models with
interface elements are developed for the HF cases of known fracture paths %
\citep{chen2009cohesive,sarris2011modeling,carrier2012numerical,sarris2013numerical,yao2015pore}%
. A fine mesh is required to accurately simulate both the toughness and
viscosity of the HF propagation. The simulation results are seriously
affected by the details of the cohesive constitutive law in the cases with a
large cohesive zone \citep{chen2009cohesive,needleman2014some,yao2015pore}.
Cohesive zone models can also be used in more general cases with unknown
crack paths by either inserting cohesive interface elements between all
finite elements of the mesh or inserting them adaptively between the
yielding finite elements \citep{simoni2003cohesive,lecampion2018numerical}.
By performing the enrichment locally in the elements intersected by a
discontinuity, the extended finite element method (XFEM) has the ability to
simulate the crack propagation \citep{belytschko2009review}, which can be
used to solve HF problems %
\citep{rethore2007two,lecampion2009extended,gordeliy2013coupling,vahab2018x}%
. Combined with the cohesive zone models, XFEM can produce the solutions of
plane-strain HF propagation cases \citep{faivre20162d,carrier2012numerical}
and is widely used to analyse static or dynamic HF problems %
\citep{mohammadnejad2013extended,wang2015numerical,paul20183d,cao2018porous}%
. Nevertheless, the XFEM-based approaches are still facing the challenges of
dealing with multi-cracks and 3D non-planar crack problems %
\citep{karihaloo2003modelling,gupta2014simulation}. Another continuous
method for HF problems is the phase field method (PFM), in which a
continuous scalar variable called the \textit{field order parameter} or 
\textit{crack field} is introduced to describe a smooth transition between
the intact and the fully broken regions \citep{ambati2015review}. The PFM is
usually solved in the FE framework and in the discretization the crack is
located by smearing it over few elements \citep{miehe2010thermodynamically}.
The extension of the PFM to poroelastic media is available in %
\citep{miehe2015minimization,mikelic2015phase,wilson2016phase,lee2016pressure}%
. But while the smeared nature of the PFM allows it to model complex
fracture behavior, it also limits its capability to accurately represent the
displacement discontinuities across the crack surface . Other continuous
approaches for simulating HF problems are those based on the boundary
element method %
\citep{siebrits2002efficient,yamamoto2004multiple,gordeliy2011fixed,rungamornrat2005numerical,xu2013interaction}
and on meshless techniques \citep{samimi2016fully,douillet2016mesh}.

The discrete approaches for HF are usually based on the discrete element
method (DEM), including the particle models %
\citep{damjanac2016application,damjanac2013three} and the lattice models %
\citep{grassl20152d,damjanac2016application}. In the DEM, the particles with
different sizes and shapes are assembled to represent the solid medium, the
deformation is described through the displacement of the particles while the
internal forces are expressed by the contacts between the particles. The
loss of contact indicates the formation of cracks. By introducing the
expression of fluid flow through the transport elements or finite volume
scheme, the DEM is applied to simulate the hydraulic fracturing behavior %
\citep{damjanac2016application,damjanac2013three,grassl20152d}.

In order to combine the advantages of different approaches, some hybrid
approaches have been proposed, for example, the combined finite-discrete
element method in \citep{zhao2014numerical,yan2016combined}, the hybrid
discrete-fracture/matrix and multi-continuum model in \citep{jiang2016hybrid}
and a hybrid discrete-continuum method in \citep{zhang2017fully}. In
addition, other attempts have been implemented, such as the numerical
manifold method \citep{yang2018hydraulic,wu2014extension}, the discretized
virtual internal bond method \citep{peng2018hydraulic} and the distinct
lattice spring model \citep{milanese2016avalanches,cao2017interaction}.
Although a large number of numerical approaches have been proposed for HF
problems, they still have some shortcomings and are facing challenges in
dealing with the problems of multi-cracks, crack bifurcation or leak-off
particularly for 3D conditions.

The theory of peridynamics (PD) was firstly introduced by Silling in 2000 %
\citep{silling2000reformulation}. It is a non-local continuum theory based
on integro-differential equations. It is particularly suited to simulate
crack propagation in structures since it allows cracks to grow naturally
without resorting to external crack growth criteria. The first version of
the PD theory, called Bond-Based PD (BB-PD), had a strong limitation because
the Poisson's ratio could only assume a fixed value %
\citep{silling2005meshfree,galvanetto2016effective,ni2019static,wang2019improved,wang2019size}%
. The present paper makes use of the Ordinary State-Based PD theory (OSB-PD) %
\citep{silling2007peridynamic}. However, some other modified versions of the
theory were also developed for solid materials with any Poisson's ratio %
\citep{zhu2017peridynamic,wang20183,chowdhury2015micropolar,zhang2018axisymmetric,zhang2019modified,diana2019bond}%
. The application of the PD models to HF problems can be traced back to %
\citep{turner2013non}, in which the problems of consolidation and leak-off
in HF are solved with an extended SB-PD model. The hydro-mechanical coupling
in \citep{turner2013non} is introduced according to Biot theory with a given
pore pressure field, without considering the fluid flow in the porous media.
\textcolor{blue}{Subsequently, the PD theory for single-phase fluid flow in porous media was developed in %
\citep{katiyar2014peridynamic,jabakhanji2015peridynamic}, in which non-local fluid transport in domains with discontinuities is considered in the context of peridynamics. Based on the proposed PD porous flow model and the extended SB-PD model, a fully
coupled PD model for HF problems was presented in \citep{ouchi2015fully}. In the fully coupled model, it assumes that a fracture space is created in any two adjacent material points when they separate beyond the critical stretch and the damage values for both the points reach a preset value. Such material points are referred as \textquotedblleft dual points\textquotedblright , not only acting as porous material but also representing the fracture space for fluid to flow in between them. The permeability at \textquotedblleft dual points\textquotedblright \ is computed by using the cubic law, where the fracture width used in the cubic law is defined according to the numerically evaluated reservoir porosity. Later, the fully coupled PD model is applied to the analysis of several typical HF cases in \citep{ouchi2017effect,ouchi2017peridynamics}. Several other works on the application of PD in hydro-mechanical coupled problems can be found in \citep{nadimi20163d,oterkus2017fully,zhang2019coupling}. In \citep{nadimi20163d}, a state-based viscoelastic peridynamic approach is developed and applied in the simulation of the fracturing process under hydraulic loading without mentioning how to compute the fluid flow field. In \citep{oterkus2017fully}, a fully coupled poroelastic formulation is presented in the context of bond-based PD theory and applied to fluid-filled fractures, and the non-local governing equation of fluid flow is also given, but the fluid flowing in cracks seems not be considered. The implicit implementation of the fully coupled PD model from \citep{ouchi2015fully} is made in \citep{zhang2019coupling} and used for solving consolidation problem and dynamic analysis of saturated porous media without considering crack propagation. In addition to above work, a numerical framework of coupling peridynamics with a edge-based finite element model for simulating HF process in Karoo is introduced in brief in \citep{turner2013coupled}. However, no detailed implementation is mentioned, and so far, no subsequent work has been published.}

PD-based numerical approaches are usually more computationally expensive
than those making use of local mechanics and FEM %
\citep{galvanetto2016effective,zaccariotto2018coupling}. 
\textcolor{blue}{For that reason and inspired by \citep{milanese2016avalanches}, in this
	paper, a hybrid modelling approach making use of FEM and PD is proposed for
	simulating the hydraulic fracturing in saturated porous media.} PD is used
to describe the deformation of the solid skeleton and to capture the
advancement of cracks, while FEM equations are used to describe the fluid
mechanics. The hydro-mechanical coupling is based on Biot theory. In the FEM
equations, the coupling terms are considered as local expression, whereas in
PD equations, the pore pressures are applied non-locally. A staggered scheme
is adopted to solve the coupled system. In each solution sequence, the
pressure field is solved first by using an implicit time integration scheme
from \citep{zienkiewicz1977finite,smith2013programming}, and then the
adaptive dynamic relaxation method %
\citep{Underwood1983dynamic,kilic2010adaptive,Ni2019Coupling} is used to
compute the displacement field. The solution algorithms for solving the
transient analysis of consolidation and fluid flow are also presented. The
relevant software is developed in MATLAB. With the developed software, a
series of numerical examples are performed to validate the proposed approach
and to demonstrate its main features. The main contributions of the paper
are:

\begin{itemize}
\item {\textcolor{blue}{The FEM and peridynamic models are combined to
solve the hydro-mechanical coupled problems;}}

\item {The coupling is described in detail and the coupling matrix of a PD
bond is presented;}

\item {An algorithm for computing the aperture of the cracks is presented,
which is much simpler than the normal displacement jump approach used in %
\citep{lee2016pressure};}

\item {The proposed approach is applied successfully in the simulation of HF
propagation and bifurcation. In the fluid-driven cases, the pressure
oscillation phenomenon discussed in \citep{cao2018porous} is observed, which
indicates the capability of the proposed method to describe the dynamic
behaviors in HF problems.}
\end{itemize}

The contents of this paper are organized as follows. In $sect.$ 2, the
OSB-PD theory is first briefly summarized and extended for poroelastic media
by using the effective stress principle, then the governing equations for
the fluid flow in fractured porous media are introduced. $Sect.$ 3 presents
the numerical implementation in detail. In $sect.$ 4, the accuracy of the
numerical simulation is verified by using two examples, then a series of
pressure- and fluid-driven HF cases are presented for further describing the
performance of the hybrid approach. Finally, $sect.$ 5 concludes the paper.

\section{Theoretical basis}

\subsection{The extended ordinary state-based peridynamic model for
saturated porous media}

\subsubsection{Basic concepts}

\begin{figure}[h]
\begin{center}
\includegraphics[scale=0.5]{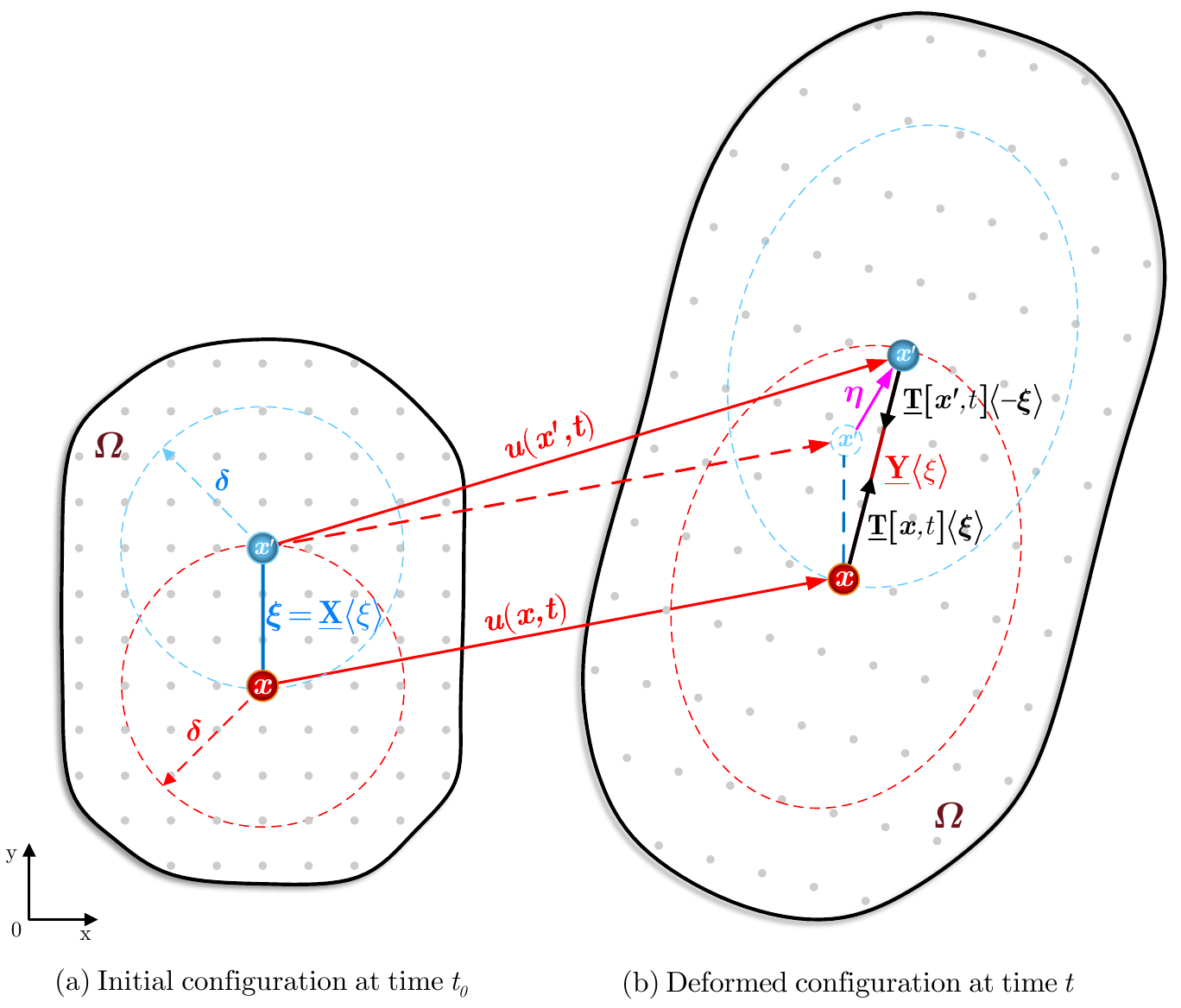}
\end{center}
\caption{The positions of two material points in the (a) initial and (b)
deformed configurations.}
\label{fig1}
\end{figure}

As shown in Fig. \ref{fig1}, a domain $\Omega $ is modelled by OSB-PD, in
which each point interacts with all the other points around it within a
prescribed horizon radius $\delta $ %
\citep{silling2007peridynamic,le2014two,sarego2016linearized,bobaru2016handbook,van2018objectivity,zhang2019new}%
. Let us consider two points $\boldsymbol{x}$ and $\boldsymbol{x}^{\prime }$
in the initial configuration of domain $\Omega $ (see Fig. \ref{fig1}a). The
bond connecting these two points, that represents their relative position,
can be expressed as:%
\begin{equation}
\boldsymbol{\xi =x}^{\boldsymbol{\prime }}-\boldsymbol{x}  \label{1.1}
\end{equation}

If $\boldsymbol{u}$ and $\boldsymbol{u}^{\boldsymbol{\prime }}$ are the
displacements of $\boldsymbol{x}$ and $\boldsymbol{x}^{\prime }$\
respectively, then the relative displacement vector can be given by:%
\begin{equation}
\boldsymbol{\eta }=\boldsymbol{u}^{\boldsymbol{\prime }}-\boldsymbol{u}
\label{1.2}
\end{equation}

The reference position vector state and deformation vector state are defined
as $\underline{\boldsymbol{X}}\left\langle \boldsymbol{\xi }\right\rangle $
and $\underline{\boldsymbol{Y}}\left\langle \boldsymbol{\xi }\right\rangle $%
, respectively, and are given by the following expressions:%
\begin{equation}
\begin{array}{ccc}
\underline{\boldsymbol{X}}\left\langle \boldsymbol{\xi }\right\rangle =%
\boldsymbol{\xi } & , & \underline{\boldsymbol{Y}}\left\langle \boldsymbol{%
\xi }\right\rangle =\boldsymbol{\xi +\eta }%
\end{array}
\label{1.3}
\end{equation}

The reference position scalar state and deformation scalar state are defined
as:%
\begin{equation}
\begin{array}{ccc}
\underline{x}=\left\Vert \underline{\boldsymbol{X}}\right\Vert & , & 
\underline{y}=\left\Vert \underline{\boldsymbol{Y}}\right\Vert%
\end{array}
\label{1.4}
\end{equation}%
where $\left\Vert \underline{\boldsymbol{X}}\right\Vert $ and $\left\Vert 
\underline{\boldsymbol{Y}}\right\Vert $\ are the norms of $\underline{%
\boldsymbol{X}}$ and $\underline{\boldsymbol{Y}}$, representing the lengths
of the bond in its initial and deformed states, respectively. For later use,
a unit state in the direction of $\underline{\boldsymbol{Y}}$ is defined as:%
\begin{equation}
\underline{\boldsymbol{M}}\left\langle \boldsymbol{\xi }\right\rangle =\frac{%
\underline{\boldsymbol{Y}}}{\left\Vert \underline{\boldsymbol{Y}}\right\Vert 
}  \label{1.5}
\end{equation}

In addition, the extension scalar state is defined to describe the
longitudinal deformation of the bond:%
\begin{equation}
\underline{e}=\underline{y}-\underline{x}  \label{1.6}
\end{equation}

In the deformed configuration, shown in Fig. \ref{fig1}b, the force density
vector states\ at time $t$ are defined as \textbf{\b{T}}$\left[ \boldsymbol{x%
},t\right] \left\langle \boldsymbol{\xi }\right\rangle $\ and \textbf{\b{T}}$%
\left[ \boldsymbol{x}^{\prime },t\right] \left\langle \boldsymbol{-\xi }%
\right\rangle $\ along the deformed bond, and their values can be different.
Thus, the OSB-PD equation of motion is given as \citep{van2018objectivity}: 
\begin{equation}
\rho \boldsymbol{\ddot{u}}\left( \boldsymbol{x},t\right) =\int\nolimits_{%
\mathcal{H}_{x}}\left\{ \text{\textbf{\b{T}}}\left[ \boldsymbol{x},t\right]
\left\langle \boldsymbol{\xi }\right\rangle -\text{ \textbf{\b{T}}}\left[ 
\boldsymbol{x}^{\prime },t\right] \left\langle \boldsymbol{-\xi }%
\right\rangle \right\} dV_{x^{\prime }}+\boldsymbol{b}\left( \boldsymbol{x}%
,t\right)  \label{1.7}
\end{equation}%
where $dV_{x^{\prime }}$ is the infinitesimal volume associated to point $%
\boldsymbol{x}^{\prime }$. $\boldsymbol{b}$ is the force density of the
external body force. $\mathcal{H}_{x}$ is the neighborhood associated with
the material point $\boldsymbol{x}$, which is usually a circle in 2D and a
sphere in 3D, it is mathematically defined as: $\mathcal{H}_{x}\left( \delta
\right) =\left\{ x\in \Omega\text{:}\left\Vert \boldsymbol{\xi }\right\Vert
\leqslant \delta \right\} $.

\subsubsection{Linear isotropic solid material}

According to %
\citep{silling2007peridynamic,bobaru2016handbook,van2018objectivity}, the 3D
elastic strain energy density at point $\boldsymbol{x}$ in a linear OSB-PD
isotropic material can be expressed as: 
\begin{equation}
W\left( \theta ,\underline{e}^{d}\right) =\frac{\kappa\theta ^{2}}{2}+\frac{%
15\mu }{2m}\int\nolimits_{\mathcal{H}_{x}}\underline{\mathit{w}}\text{ }%
\underline{e}^{d}\text{ }\underline{e}^{d}\text{d}V_{x^{\prime }}
\label{1.8}
\end{equation}%
where $k$ and $\mu $ are bulk modulus and shear modulus. $\theta $\ is the
volume dilatation value, $\underline{e}^{d}$ is the deviatoric extension
state, and they are defined as: 
\begin{equation}
\theta =\frac{3}{m}\int\nolimits_{\mathcal{H}_{x}}\left( \underline{\mathit{w%
}}\text{ }\underline{x}\text{ }\underline{e}\right) \text{d}V_{x^{\prime }}
\label{1.9}
\end{equation}%
\begin{equation}
\underline{e}^{d}=\underline{e}-\frac{\theta \underline{x}}{3}  \label{1.10}
\end{equation}%
in which $m$ is called the weighted volume, given as: 
\begin{equation}
m=\int\nolimits_{\mathcal{H}_{x}}\underline{\mathit{w}}\left\Vert 
\boldsymbol{\xi }\right\Vert ^{2}\text{d}V_{x^{\prime }}  \label{1.11}
\end{equation}

In the above formulae, $\underline{\mathit{w}}$ is an influence function,
which is usually adopted as $\underline{\mathit{w}}=\exp \left( -\frac{%
\left\Vert \boldsymbol{\xi }\right\Vert ^{2}}{\delta ^{2}}\right) $. Based
on the above notions, the force density vector state \textbf{\b{T}} can be
defined as:%
\begin{equation}
\text{\textbf{\b{T}}}=\underline{t}\cdot \underline{\boldsymbol{M}}
\label{1.12}
\end{equation}%
where $\underline{t}$ is called the force density scalar state. In %
\citep{van2018objectivity}, $\underline{t}$ is defined in 3D context as:%
\begin{equation}
\underline{t}=3\kappa\theta \frac{\underline{\mathit{w}}\text{ }\underline{x}%
}{m}+15\mu \underline{e}^{d}\frac{\underline{\mathit{w}}}{m}\text{ }
\label{1.13}
\end{equation}

The expression of PD force density scalar state in plane strain conditions
is also introduced \citep{Ni2019Coupling}: 
\begin{equation}
\underline{t}=2(\kappa-\frac{\mu }{3})\theta \frac{\underline{\mathit{w}}%
\text{ }\underline{x}}{m}+8\mu \underline{e}^{d}\frac{\underline{\mathit{w}}%
}{m}\text{ }  \label{1.13_1}
\end{equation}
\textcolor{blue}{where $\theta$ is caculated by:
\begin{equation}
\theta =\frac{2}{m}\int\nolimits_{\mathcal{H}_{x}}\left( \underline{\mathit{w}}\text{ }\underline{x}\text{ }\underline{e}\right) \text{d}V_{x^{\prime }}
\label{1.13_2}
\end{equation}}

\subsubsection{Failure criterion}

To describe the material failure and crack propagation, a failure criterion
is required in the PD-based numerical models. 
\textcolor{blue}{Three different types of failure criteria are commonly used in SB-PD models. The first is called \textquotedblleft critical bond stretch\textquotedblright \ criterion \citep{silling2005meshfree}. It was derived in the context of the BB-PD theory based on mode I fracture. The second is the bond-level energy based failure criterion derived as in the critical bond stretch criterion but involving deviatoric deformation \citep{foster2011energy}. The bond-level energy based criterion is thoroughly introduced and successfully applied to the quantitative fracture analysis of isotropic and orthotropic materials in \citep{zhang2018state,zhang2019state}. The last criterion is an alternative bond failure criterion in terms of strain invariants introduced in \citep{warren2009non}, which can be used as a bridge to directly apply the damage model of the classical theory of mechanics in peridynamics. As discussed in \citep{dipasquale2017discussion}, all the three criteria can be suitable in the analysis of mode I fracture problems and reproduce the experimentally observed behavior. In the considered HF problems, the fracture mainly occurs with the opening crack propagation, hence for simplicity and following the recommendation in \citep{littlewood2015roadmap}, the \textquotedblleft critical bond
stretch\textquotedblright \ criterion is adopted for the SB-PD model in this
paper.} 

The stretch value of bond $\boldsymbol{\xi }$ is defined as %
\citep{Ni2019Coupling}:%
\begin{equation}
s\left\langle \boldsymbol{\xi }\right\rangle =\frac{\underline{e}%
\left\langle \boldsymbol{\xi }\right\rangle }{\underline{x}\left\langle 
\boldsymbol{\xi }\right\rangle }  \label{1.13.1}
\end{equation}%
and the critical stretch value can be expressed as %
\citep{zaccariotto2018coupling,ni2018peridynamic,zhang2020experimental}: 
\begin{equation}
s_{c}=\left\{ 
\begin{array}{lll}
\sqrt{\frac{5G_{c}}{6E\delta }} & , & 3D \\ 
\sqrt{\frac{5G_{c}}{12E\delta }} & , & \text{plane strain}%
\end{array}%
\right.  \label{1.13.2}
\end{equation}%
where $G_{c}$ is the critical energy release rate for mode I fracture. Eq. (%
\ref{1.13.2}) was described in %
\citep{zaccariotto2018coupling,ni2018peridynamic} for the case of BB-PD. Its
use in SB-PD is not fully justified, but it is usually accepted %
\citep{dipasquale2017discussion,littlewood2015roadmap}.

In addition, a scalar variable $\mathit{\varrho }$ is introduced in the
OSB-PD models to indicate the connection status of the bonds %
\citep{zaccariotto2018coupling,ni2018peridynamic}: 
\begin{equation}
\underline{\mathit{\varrho }}\left\langle \boldsymbol{\xi }\right\rangle 
\boldsymbol{=}\left\{ 
\begin{array}{ccc}
1 & , & \text{if }s\left\langle \boldsymbol{\xi }\right\rangle <s_{c}\text{
, for all }0<\overline{t}<t \\ 
0 & , & \text{otherwise}%
\end{array}%
\right.  \label{1.13.3}
\end{equation}%
thus the damage value $\varphi _{x}$ at point $\boldsymbol{x}$\ in the
system can be defined as: 
\begin{equation}
\varphi _{x}=1-\frac{\int\nolimits_{\mathcal{H}_{x}}\underline{\mathit{w}}%
\left\langle \boldsymbol{\xi }\right\rangle \text{ }\underline{\mathit{\
\varrho }}\left\langle \boldsymbol{\xi }\right\rangle \text{d}V_{x^{\prime }}%
}{\int\nolimits_{\mathcal{H}_{x}}\underline{\mathit{w}}\left\langle 
\boldsymbol{\xi }\right\rangle \text{d}V_{x^{\prime }}}  \label{1.13.4}
\end{equation}

\textcolor{blue}{\begin{figure}[h]
	\begin{center}
		\includegraphics[scale=0.5]{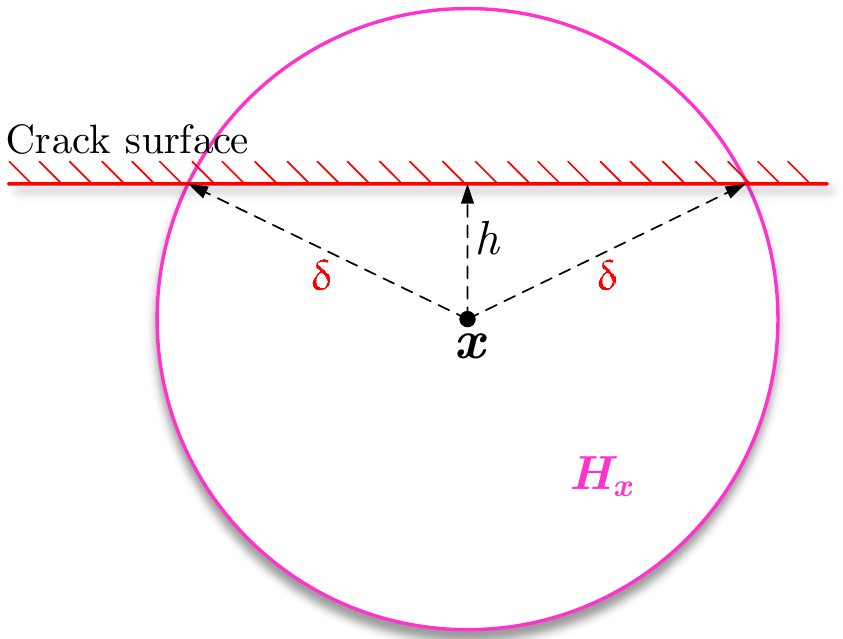}
	\end{center}
	\caption{Neighborhood of point $\boldsymbol{x}$ crossed by a crack.}
	\label{fig2}
\end{figure}
Let us consider a material point $\boldsymbol{x}$ existing in a cracked solid modeled by peridynamics, as shown in Fig. \ref{fig2}, the neighborhood of point $\boldsymbol{x}$ is crossed by a crack. The area on the upper side of the crack can be given as:
\begin{equation}
S_{up}=\frac{1}{2 }\arccos \left(\left( \frac{h}{\delta }
	\right) ^{2}-1\right)\delta^{2} -h\sqrt{\delta ^{2}-h^{2}} \label{1.2_1}
\end{equation}
The damage value $\varphi _{x}$ at the point $\boldsymbol{x}$ can be represented geometrically with the ratio of the $S_{up}$ to the total area of the neighborhood, which is given as:
\begin{equation}
\varphi _{x}=\frac{1}{2\pi }\left[ \arccos \left( 2\left( \frac{h}{\delta }
\right) ^{2}-1\right) -2\frac{h}{\delta }\sqrt{1-\left( \frac{h}{\delta }
	\right) ^{2}}\right]  \label{1.2_2}
\end{equation}
\begin{figure}[h]
	\begin{center}
		\includegraphics[scale=0.8]{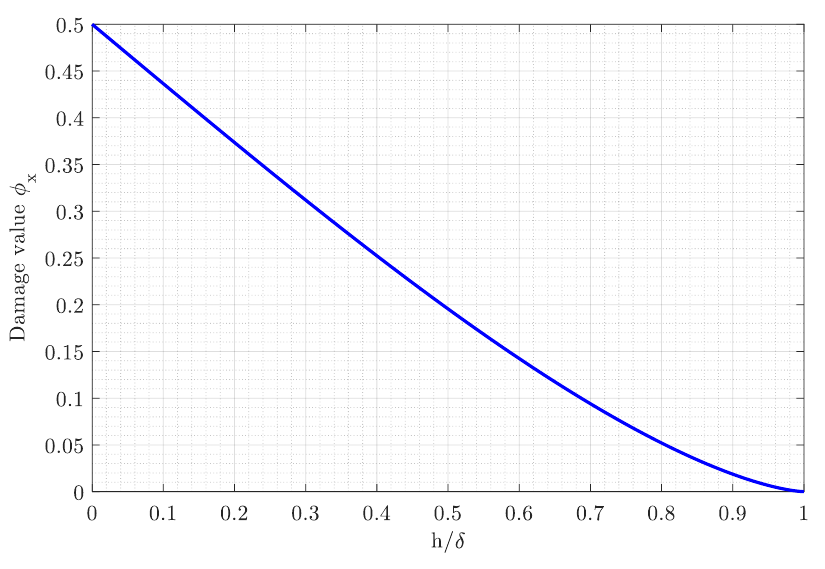}
	\end{center}
	\caption{Evolution of damage value $\varphi _{x}$ as the position changes.}
	\label{fig2_2}
\end{figure}
The evolution of damage value $\varphi _{x}$ as the position changes is expressed in Fig. \ref{fig2_2}. Thus, the value of $\varphi _{x}$ theoretically ranges from 0 to 0.5 at the continuous level. }

\subsubsection{Effective stress principle in ordinary state-based
peridynamics}

\bigskip In the theory of single-phase fluid saturated porous media %
\citep{lewis1998finite,zienkiewicz1999computational}, the effective stress
principle is expressed as:%
\begin{equation}
\boldsymbol{\sigma }^{tot}=\boldsymbol{\sigma }^{eff}-\alpha \boldsymbol{I}p
\label{1.14}
\end{equation}%
where $\boldsymbol{\sigma }^{tot}$ and $\boldsymbol{\sigma }^{eff}$ are the
total and the effective stress tensors, $\alpha $\ is the Biot coefficient, $%
p$\ is the pore pressure, and $\boldsymbol{I}$ is the unit tensor. In Eq.(%
\ref{1.14}), the negative sign is introduced as it is a general convention
to take tensile components of stress as positive %
\citep{zienkiewicz1999computational}. 

In Eq.(\ref{1.13}), the force density scalar state has been divided into a
dilatational term, $3\kappa\theta \frac{\underline{\mathit{w}}\text{ }%
\underline{x}}{m}$, and a deviatoric term, $15\mu \underline{e}^{d}\frac{%
\underline{\mathit{w}}}{m}$. In classical continuum theory, the volume
dilatation is defined by $\theta = \varepsilon _{xx}+\varepsilon
_{yy}+\varepsilon _{zz}$. By neglecting the coefficient part
related to the discretization, the dilatational term in the force density
scalar state should correspond to volume stress, $\left( \sigma _{xx}+\sigma
_{yy}+\sigma _{zz}\right) /3$, in classical continuum solid theory. Thus, by
using the effective stress principle expressed in Eq.(\ref{1.14}), the pore
pressure term can be directly added to the dilatational term and the
expression of the force density scalar state for the 3D problem of saturated
porous media can be obtained as %
\citep{turner2013non,york2018advanced,zhang2019coupling}:%
\begin{equation}
\begin{array}{c}
\underline{t}^{tot}=3\left( \kappa\theta -\alpha p\right) \frac{\underline{%
\mathit{w}}\text{ }\underline{x}}{m}+15\mu \underline{e}^{d}\frac{\underline{%
\mathit{w}}}{m} \\ 
\Downarrow \\ 
\underline{t}^{tot}=\underline{t}^{eff}-3\alpha p\frac{\underline{\mathit{w}}%
\text{ }\underline{x}}{m}%
\end{array}
\label{1.15}
\end{equation}%
in which, $\underline{t}^{eff}$ corresponds to $\underline{t}$ in Eq. (\ref%
{1.13}).

Thus, the coupled peridynamic equation of motion will be given as:%
\begin{equation}
\begin{array}{l}
\rho \boldsymbol{\ddot{u}}\left( \boldsymbol{x},t\right) =\int\nolimits_{%
\mathcal{H}_{x}}\left\{ \text{\textbf{\b{T}}}\left[ \boldsymbol{x},p,t\right]
\left\langle \boldsymbol{\xi }\right\rangle -\text{ \textbf{\b{T}}}\left[ 
\boldsymbol{x}^{\prime },p^{\prime },t\right] \left\langle \boldsymbol{-\xi }%
\right\rangle \right\} dV_{x^{\prime }}+\boldsymbol{b}\left( \boldsymbol{x}%
,t\right) \\ 
=\int\nolimits_{\mathcal{H}_{x}}\left\{ \text{\textbf{\b{T}}}\left[ 
\boldsymbol{x},t\right] \left\langle \boldsymbol{\xi }\right\rangle -\text{ 
\textbf{\b{T}}}\left[ \boldsymbol{x}^{\prime },t\right] \left\langle 
\boldsymbol{-\xi }\right\rangle \right\} dV_{x^{\prime }}-3\alpha
\int\nolimits_{\mathcal{H}_{x}}\left[ p\frac{\underline{\mathit{w}}\text{ }%
\underline{x}}{m\left( \boldsymbol{x}\right) }\underline{\boldsymbol{M}}%
\left\langle \boldsymbol{\xi }\right\rangle -\text{ }p^{\prime }\frac{%
\underline{\mathit{w}}\text{ }\underline{x}}{m\left( \boldsymbol{x}^{\prime
}\right) }\underline{\boldsymbol{M}}\left\langle -\boldsymbol{\xi }%
\right\rangle \right] dV_{x^{\prime }}+\boldsymbol{b}\left( \boldsymbol{x}%
,t\right)%
\end{array}
\label{1.16}
\end{equation}%
where $p$ and $p^{\prime }$\ are the values of pore pressure at nodes $%
\boldsymbol{x}$ and $\boldsymbol{x}^{\prime }$.

In addition, the expression of the force density scalar state for the 2D
problem will be obtained similarly:%
\begin{equation}
\underline{t}^{tot}=\underline{t}^{eff}-2\alpha p\frac{\underline{\mathit{w}}%
\text{ }\underline{x}}{m}  \label{1.17}
\end{equation}%
in which $\underline{t}^{eff}$ is the force density scalar state in plane
strain cases expressed in Eq. (\ref{1.13_1}).

\subsection{Governing equations for flow in fractured porous media}

\begin{figure}[h!]
\centering  
\subfloat[]{\includegraphics[scale=0.4]{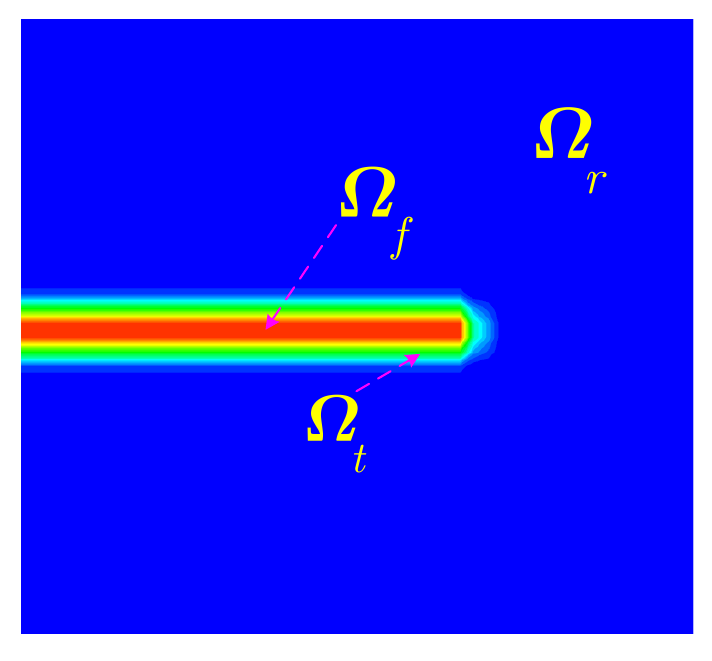}	%
\label{fig2_1:sub1}}\hspace{0.5in} \subfloat[]{%
\includegraphics[scale=0.4]{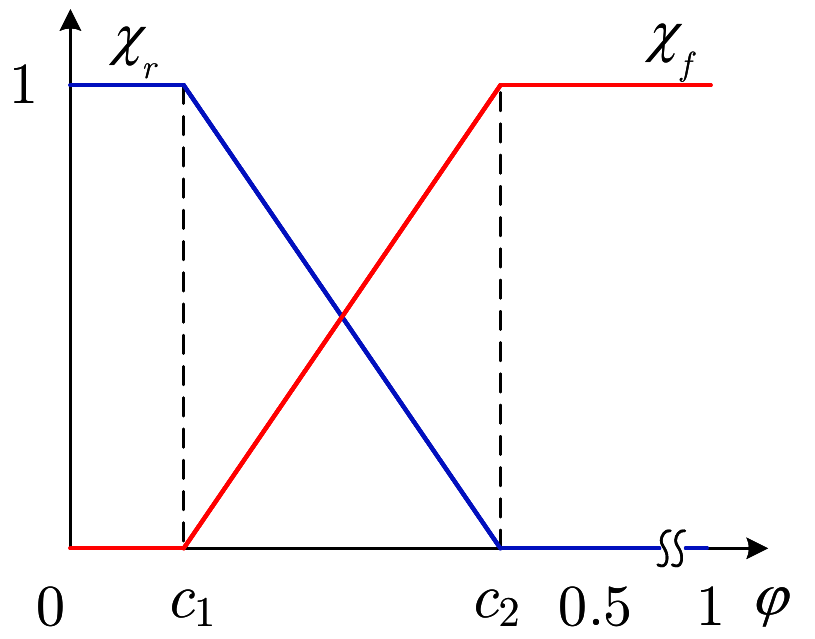}	%
\label{fig2_1:sub2}}
\caption{(a) The definition of the $\Omega _{r}$, $\Omega _{t}$ and $\Omega
_{f}$; (b) linear indicator functions $\protect\chi _{r}$ and $\protect\chi %
_{f}$.}
\label{fig2_1}
\end{figure}

To formulate the governing equation for flow in the domain $\Omega $, the
whole domain $\Omega $ is divided into three parts: $\Omega _{r}$, $\Omega
_{f}$\ and $\Omega _{t}$ (see Fig. \ref{fig2_1:sub1}), representing the
unbroken domain (reservoir domain), the fracture domain and the transition
domain between $\Omega _{r}$ and $\Omega _{f}$. Following references %
\citep{lee2016pressure,zhou2018phase,zhou2019phase}, we use the peridynamic
damage field ($\varphi $ in Eq. (\ref{1.13.4})) as an indicator. As shown in
Fig. \ref{fig2_1:sub2}, two threshold values ($c_{1}$ and $c_{2}$) are set
to identify the three flow domains: the reservoir domain is defined as $\varphi
\leqslant c_{1}$, the fracture domain as $\varphi \geqslant c_{2}$ and
the transition domain as $c_{1}\leq \varphi \leq c_{2}$.

We use Darcy's law to describe the flow field in the saturated porous media,
then the mass balance in the reservoir domain is given as:%
\begin{equation}
\rho _{r}s_{r}\frac{\partial p}{\partial t}+\rho _{r}\alpha _{r}\frac{%
\partial \varepsilon _{v}}{\partial t}+\rho _{r}\nabla \cdot \left[ \frac{%
k_{r}}{\mu_{w}}\left( -\nabla p+\rho _{r}g\right) \right] =q_{r}
\label{2.2.1}
\end{equation}%
where $\alpha _{r}$,\ $s_{r}$, $k_{r}$, $q_{r}$ and $\rho _{r}$\ are the
Biot coefficient, storage coefficient, permeability, source term and the
density of the media in the reservoir domain, respectively; $\mu_{w}$ is the
viscosity coefficient of the fluid in the reservoir domain. $g$\ is gravity
and $\varepsilon _{v}$\ is volumetric strain. The storage coefficient is
given as: 
\begin{equation}
s_{r}=\frac{\left( \alpha _{r}-n_{r}\right) \left( 1-\alpha _{r}\right) }{%
K_{r}}+\frac{n_{r}}{K_{w}}  \label{2.2.2}
\end{equation}%
where $K_{r}$ and $K_{w}$ are the bulk moduli of solid and fluid in the
reservoir domain, $n_{r}$\ is the porosity.

\textcolor{blue}{In order to use cubic law to evaluate the permeability in
the fracture domain, the fracture domain needs to be assumed fully filled by
fluid and the porosity $n_{f}=1$ is adopted. Supposing that the fluid in the
fracture is incompressible, the volumetric strain term in Eq. (\ref{2.2.1})
vanishes \citep{van2019monolithic,lee2016pressure}:} 
\begin{equation}
\rho _{f}s_{f}\frac{\partial p}{\partial t}+\rho _{f}\nabla \cdot \left[ 
\frac{k_{f}}{\mu_{w}}\left( -\nabla p+\rho _{f}g\right) \right] =q_{f}
\label{2.2.3}
\end{equation}%
where $s_{f}$, $k_{f}$, $q_{f}$ and $\rho _{f}$\ are the storage
coefficient, permeability, source term and the density of the fluid in the
fracture domain.

Following the descriptions in %
\citep{lee2016pressure,zhou2018phase,zhou2019phase}, two linear indicator
functions are defined to connect the transition domain with the reservoir
and fracture domains: $\chi _{r}$ and $\chi _{f}$, which satisfy:
\begin{equation}
\chi _{r}\left( \cdot ,\varphi \right) :=\chi _{r}\left( \boldsymbol{x}%
,t,\varphi \right) =1\text{ in }\Omega _{r}\left( t\right) \text{ and }\chi
_{r}\left( \cdot ,\varphi \right) =0\text{ in }\Omega _{f}\left( t\right)
\label{2.2.4}
\end{equation}%
\begin{equation}
\chi _{f}\left( \cdot ,\varphi \right) :=\chi _{f}\left( \boldsymbol{x}%
,t,\varphi \right) =1\text{ in }\Omega _{f}\left( t\right) \text{ and }\chi
_{f}\left( \cdot ,\varphi \right) =0\text{ in }\Omega _{r}\left( t\right)
\label{2.2.5}
\end{equation}

In the transition domain, the two linear functions are defined as (see Fig. \ref{fig2_1:sub2}):%
\begin{equation}
\chi _{r}\left( \cdot ,\varphi \right) =\frac{c_{2}-\varphi }{c_{2}-c_{1}},%
\text{ and }\chi _{f}\left( \cdot ,\varphi \right) =\frac{\varphi -c_{1}}{%
c_{2}-c_{1}}  \label{2.2.6}
\end{equation}

The detailed illustration of the linear indicator functions $\chi _{r}$ and $%
\chi _{f}$ can be found in %
\citep{lee2016pressure,zhou2018phase,zhou2019phase}. Fluid and solid
properties in the transition domain can be obtained by interpolating those
properties of reservoir and fracture domains by using the indicator
functions: $\rho _{T}=\rho _{r}\chi _{r}+\rho _{f}\chi _{f}$, $\alpha
_{T}=\alpha _{r}\chi _{r}+\alpha _{f}\chi _{f}$, $n_{T}=n_{r}\chi
_{r}+n_{f}\chi _{f}$, $k_{T}=k_{r}\chi _{r}+k_{f}\chi _{f}$, $%
s_{T}=s_{r}\chi _{r}+s_{f}\chi _{f}$. \textcolor{blue}{It should be noted in
particular that the Biot coefficient in the fracture domain is taken as
$\alpha_{f}=1$.} Therefore, the governing equation for the flow in the
transition domain can be given as:%
\begin{equation}
\rho _{T}s_{T}\frac{\partial p}{\partial t}+\rho _{T}\alpha _{T}\frac{%
\partial \varepsilon _{v}}{\partial t}+\rho _{T}\nabla \cdot \left[ \frac{%
k_{T}}{\mu_{w}}\left( -\nabla p+\rho _{T}g\right) \right] =q_{T}
\label{2.2.7}
\end{equation}

Moreover, the cubic law is used to evaluate the permeability of the fracture
domain \citep{zimmerman1994hydraulic,cao2017interaction}:%
\begin{equation}
k_{f}=\frac{1}{12}a^{2}  \label{2.2.8}
\end{equation}%
where $a$ is the aperture of the crack, which can be obtained by using the
method described in $sect.$ $3.3$.

\section{Discretization and numerical implementation}

\subsection{FEM discretization of the governing equation for the fluid flow}

Using the Galerkin finite element method, the discretized governing
equations of the fluid flow in classical continuum theory assume the
following form \citep{lewis1998finite}:%
\begin{equation}
\boldsymbol{S\dot{p}}+\boldsymbol{Q}^{T}\boldsymbol{\dot{u}}+\boldsymbol{Hp=q%
}^{w}  \label{3.1.1}
\end{equation}%
where: $\boldsymbol{Q}$ is the coupling matrix, $\boldsymbol{H}$ is the
permeability matrix, $\boldsymbol{S}$ is the compressibility matrix, all of
them can be obtained by assembling the corresponding element matrices. When
using the shape functions $\boldsymbol{N}_{u}$ (for displacement) and $%
\boldsymbol{N}_{p}$ (for pressure), the matrices in Eq.(\ref{3.1.1}) can be
given as \citep{milanese2016avalanches,peruzzo2019dynamics}:%
\begin{equation}
\boldsymbol{Q=}\int\nolimits_{\Omega }\left( \boldsymbol{LN}_{u}\right)
^{T}\alpha \boldsymbol{mN}_{p}\text{d}\Omega  \label{3.1.2}
\end{equation}%
\begin{equation}
\boldsymbol{H=}\int\nolimits_{\Omega }\left( \nabla \boldsymbol{N}%
_{p}\right) ^{T}\frac{k}{\mu_{w}}\left( \nabla \boldsymbol{N}_{p}\right) 
\text{d}\Omega  \label{3.1.3}
\end{equation}%
\begin{equation}
\boldsymbol{S=}\int\nolimits_{\Omega }\boldsymbol{N}_{p}^{T}s\boldsymbol{N}%
_{p}\text{d}\Omega  \label{3.1.4}
\end{equation}%
in which $\boldsymbol{L}$ is the differential operator defined as: 
\begin{equation}
\boldsymbol{L=}\left[ 
\begin{array}{ccc}
\frac{\partial }{\partial x} & 0 & 0 \\ 
0 & \frac{\partial }{\partial y} & 0 \\ 
0 & 0 & \frac{\partial }{\partial z} \\ 
\frac{\partial }{\partial y} & \frac{\partial }{\partial x} & 0 \\ 
0 & \frac{\partial }{\partial z} & \frac{\partial }{\partial y} \\ 
\frac{\partial }{\partial z} & 0 & \frac{\partial }{\partial x}%
\end{array}%
\right]  \label{3.1.5}
\end{equation}%
$\boldsymbol{m}$\ is a vector used instead of the unit tensor $\boldsymbol{I}
$ in Eq.(\ref{1.14}) given as: 
\begin{equation}
\boldsymbol{m}=\left[ 1,1,1,0,0,0\right] ^{T}  \label{3.1.9}
\end{equation}

Note that the coupling term $\left( \rho \alpha \frac{\partial \varepsilon
_{v}}{\partial t}\right) $ caused by volumetric change will be removed from
the governing equation for the flow in the fracture domain.

\subsection{Discretization of the peridynamic equations}

After discretization, the peridynamic equation of motion of the current node 
$\boldsymbol{x}_{i}$ is written as:%
\begin{equation}
\begin{array}{l}
\rho \boldsymbol{\ddot{u}}_{i}^{t}=\sum_{j=1}^{N_{H_{i}}}\left\{ \text{%
\textbf{\b{T}}}\left[ \boldsymbol{x}_{i},t\right] \left\langle \boldsymbol{%
\xi }_{ij}\right\rangle -\text{ \textbf{\b{T}}}\left[ \boldsymbol{x}_{j},t%
\right] \left\langle \boldsymbol{-\xi }_{ij}\right\rangle \right\} \cdot
V_{j} \\ 
-3\alpha \sum_{j=1}^{N_{H_{i}}}\left[ p_{i}\frac{\underline{\mathit{w}}%
\left\langle \boldsymbol{\xi }_{ij}\right\rangle \text{ }\underline{x}%
\left\langle \boldsymbol{\xi }_{ij}\right\rangle }{m\left( \boldsymbol{x}%
_{i}\right) }\underline{\boldsymbol{M}}\left\langle \boldsymbol{\xi }%
_{ij}\right\rangle -\text{ }p_{j}\frac{\underline{\mathit{w}}\left\langle 
\boldsymbol{\xi }_{ij}\right\rangle \text{ }\underline{x}\left\langle 
\boldsymbol{\xi }_{ij}\right\rangle }{m\left( \boldsymbol{x}_{j}\right) }%
\underline{\boldsymbol{M}}\left\langle -\boldsymbol{\xi }_{ij}\right\rangle %
\right] \cdot V_{j}+\boldsymbol{b}_{i}^{t}%
\end{array}
\label{3.1.10}
\end{equation}%
where $N_{H_{i}}$ is the number of family nodes of $\boldsymbol{x}_{i}$, $%
\boldsymbol{x}_{j}$\ is $\boldsymbol{x}_{i}$'s family node, $V_{j}$\ is the
volume of node $\boldsymbol{x}_{j}$.

Under the assumption of small deformation, the global form of Eq.(\ref%
{3.1.10}) can be written as:%
\begin{equation}
\boldsymbol{M}^{PD}\boldsymbol{\ddot{u}}+\boldsymbol{K}^{PD}\boldsymbol{u}-%
\boldsymbol{Q}^{PD}p=f^{PD}  \label{3.1.11}
\end{equation}%
in which $\boldsymbol{M}^{PD}$, $\boldsymbol{K}^{PD}$ and $\boldsymbol{Q}%
^{PD}$\ are the mass, stiffness and \textquotedblleft
coupling\textquotedblright\ matrices of the PD domain, respectively. Note
that $\boldsymbol{M}^{PD}$\ is usually taken as lumped mass matrix. The
method of obtaining $\boldsymbol{K}^{PD}$ of OSB-PD equations can be found
in \citep{sarego2016linearized}. Assume that $\underline{\boldsymbol{M}}%
\left\langle \boldsymbol{\xi }_{ij}\right\rangle =\left[ M_{x},M_{y},M_{z}%
\right] $, then the \textquotedblleft coupling\textquotedblright\ matrix for 
$\boldsymbol{\xi }_{ij}$ can be given as:%
\begin{equation}
\boldsymbol{Q_{ij}^{PD}=}3\alpha \underline{\mathit{w}}\left\langle 
\boldsymbol{\xi }_{ij}\right\rangle \underline{x}\left\langle \boldsymbol{%
\xi }_{ij}\right\rangle V_{i}V_{j}\left[ 
\begin{array}{cc}
\frac{M_{x}}{m\left( \boldsymbol{x}_{i}\right) } & \frac{M_{x}}{m\left( 
\boldsymbol{x}_{j}\right) } \\ 
\frac{M_{y}}{m\left( \boldsymbol{x}_{i}\right) } & \frac{M_{y}}{m\left( 
\boldsymbol{x}_{j}\right) } \\ 
\frac{M_{z}}{m\left( \boldsymbol{x}_{i}\right) } & \frac{M_{z}}{m\left( 
\boldsymbol{x}_{j}\right) } \\ 
-\frac{M_{x}}{m\left( \boldsymbol{x}_{i}\right) } & -\frac{M_{x}}{m\left( 
\boldsymbol{x}_{j}\right) } \\ 
-\frac{M_{y}}{m\left( \boldsymbol{x}_{i}\right) } & -\frac{M_{y}}{m\left( 
\boldsymbol{x}_{j}\right) } \\ 
-\frac{M_{z}}{m\left( \boldsymbol{x}_{i}\right) } & -\frac{M_{z}}{m\left( 
\boldsymbol{x}_{j}\right) }%
\end{array}%
\right]  \label{3.1.12}
\end{equation}%
where the influence function $\underline{\mathit{w}}\left\langle \boldsymbol{%
\xi }_{ij}\right\rangle $ used for the coupling bonds should be specified as 
$\underline{\mathit{w}}=1$.

\subsection{Coupling of the solid and fluid}

\begin{figure}[h!]
\centering  
\subfloat[Effect of fluid on
solid.]{				\includegraphics[scale=0.55]{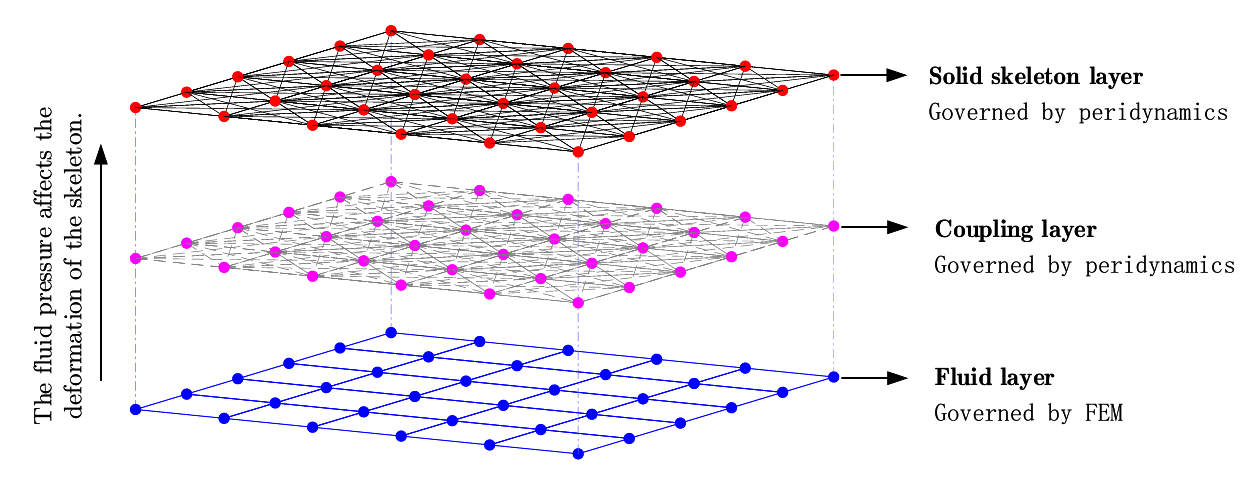}				\label{fig3_2:sub1}}%
\\
\subfloat[Effect of solid on
fluid.]{				\includegraphics[scale=0.55]{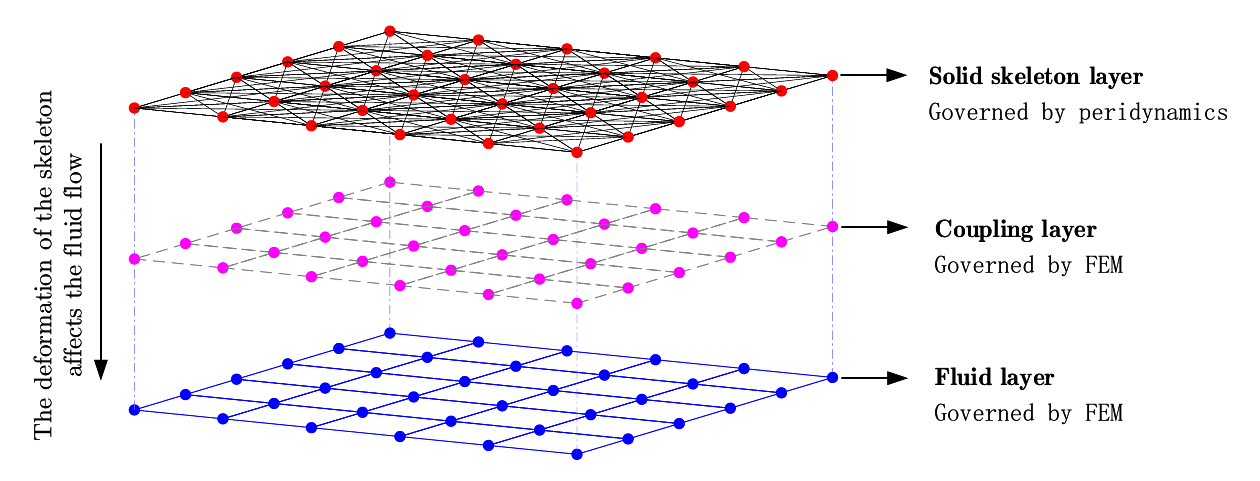}				\label{fig3_2:sub2}}
\caption{Schematic diagram of bidirectional influence in the coupling
process.}
\label{fig3_2}
\end{figure}

This section describes the coupling process in planar conditions. In a plane
discretization, see Fig. \ref{fig3_2}, the solid portion is discretized by
PD nodes and the fluid portion is discretized by 4-node FE element. PD nodes
and FE nodes share the same node coordinates. 
\textcolor{blue}{A bidirectional interaction between solid and liquid
	discretization is realized to link the local and non-local equations (Eqs.(\ref{3.1.1}) and (\ref{3.1.11})).}
\textcolor{blue}{The same shape functions are used for generating the coupling
matrix and the other fluid matrices}, which are defined in terms of the
normalized natural domain (i.e. $-1\leqslant \xi \leqslant 1$ and $%
-1\leqslant \eta \leqslant 1$) and four Gauss points are used in the
numerical integration.%
\begin{equation}
\boldsymbol{N}_{u}=\boldsymbol{N}_{p}\Rightarrow \left\{ 
\begin{array}{c}
N_{1}=\frac{1}{4}\left( 1-\xi \right) \left( 1-\eta \right) \\ 
N_{2}=\frac{1}{4}\left( 1-\xi \right) \left( 1+\eta \right) \\ 
N_{3}=\frac{1}{4}\left( 1+\xi \right) \left( 1+\eta \right) \\ 
N_{4}=\frac{1}{4}\left( 1+\xi \right) \left( 1-\eta \right)%
\end{array}%
\right.  \label{3.13}
\end{equation}

\begin{figure}[h!]
\begin{center}
\includegraphics[scale=0.375]{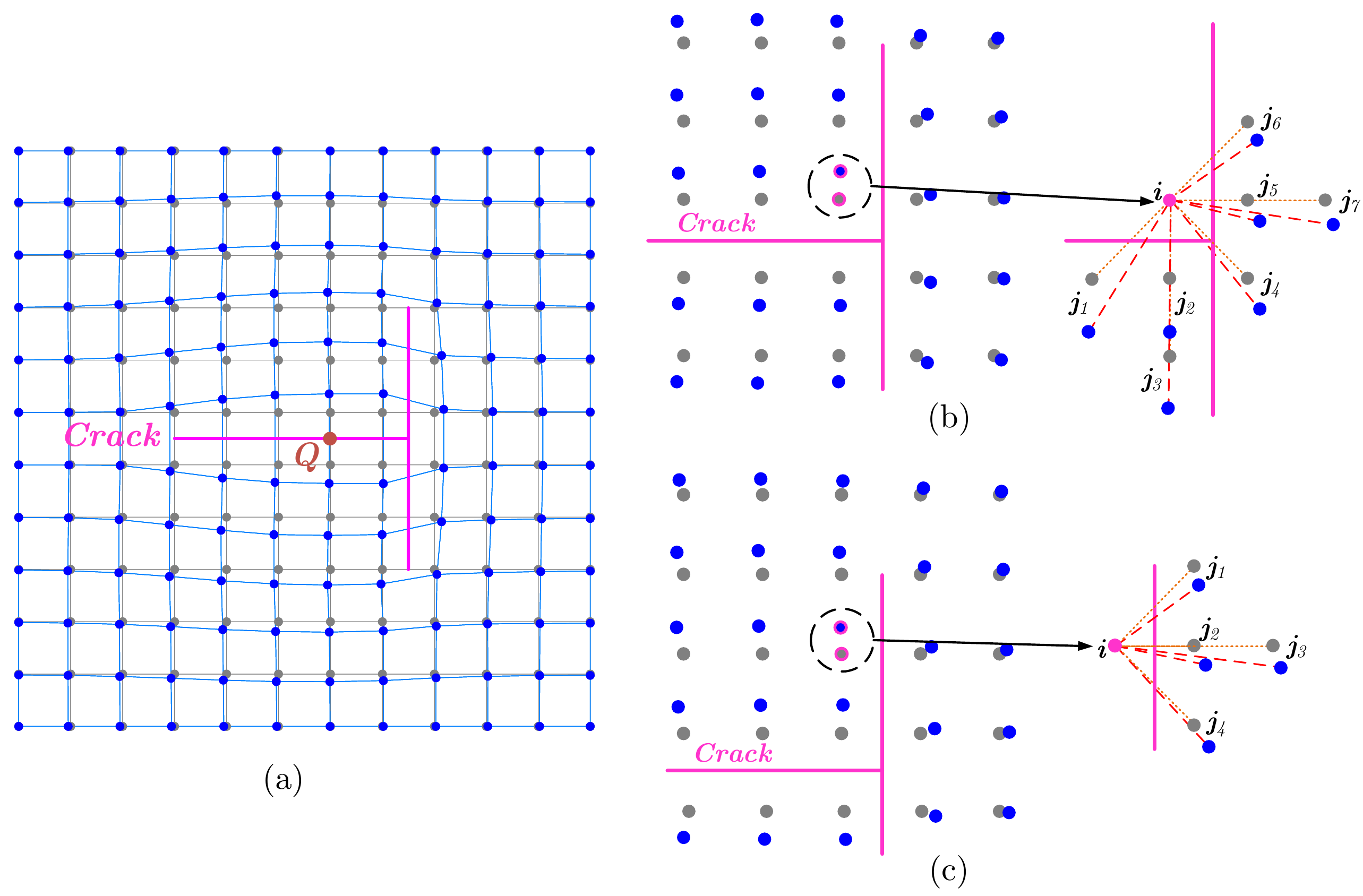}
\end{center}
\caption{\textcolor{blue}{Schematic diagram for the calculation of the
aperture value at PD node.}}
\label{fig3_1}
\end{figure}

\begin{figure}[h!]
\begin{center}
\includegraphics[scale=0.5]{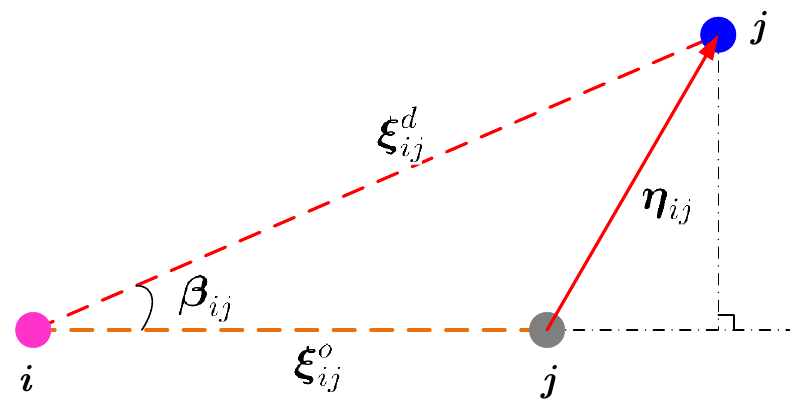} 
\end{center}
\caption{\textcolor{blue}{The decomposition of the relative displacement
vector along and perpendicular to the direction of the original bond.}}
\label{fig3_3}
\end{figure}

\textcolor{blue}{As shown in Fig. \ref{fig3_1}(a), let us consider a plate notched by a T-shaped crack with a fluid injection at the point $Q$ of the horizontal crack. Under the action of the pore pressure, the solid skeleton deforms gradually. In Fig. \ref{fig3_1}(a), the initial mesh is grey and the deformed mesh is blue, at the meantime, the initial position of the nodes are represented by grey particles and the blue particles are the deformed position. To use Eq.(\ref{2.2.8}) to evaluate the permeability
coefficient in the crack domain, the aperture values at the PD nodes near
crack surfaces need to be calculated. Here we choose two representative nodes as reference to explain the procedure for calculating the aperture values. As shown in Figs. \ref{fig3_1}(b) and \ref{fig3_1}(c), the magenta node
is the current PD node marked as NODE $\boldsymbol{i}$, while the nodes in its neighborhood that need to be used to calculate the aperture value at NODE $\boldsymbol{i}$ are marked as NODES $\boldsymbol{j}$ ($m_{ratio}=2$). The right part of Figs. \ref{fig3_1}(b) and \ref{fig3_1}(c) are showing the relative position changes of NODES $\boldsymbol{j}$ to NODE $\boldsymbol{i}$. The relative positions of NODE $\boldsymbol{i}$ and of one of the NODES $\boldsymbol{j}$ can be shown as in Fig. \ref{fig3_3}. The relative displacement vector $\eta_{ij}$ as in Fig. \ref{fig3_3} can be decomposed into two components: along and perpendicular to the direction of  the original bond supposing that they represent the opening displacement and the dislocation of the crack, respectively. Hence the aperture value related to the bond $\boldsymbol{\xi}_{ij}$ can be expressed geometrically as:
\begin{equation}
a_{ij}=\left\Vert \boldsymbol{\xi }_{ij}\right\Vert ^{d}\cos \beta_{ij}-\left\Vert \boldsymbol{\xi }_{ij}\right\Vert ^{o} \label{3.14}
\end{equation}
The decomposition of the relative displacement vector can avoid updating the permeability coefficients at closed cracks with dislocations and can also reduce the error in updating the permeability coefficients at opening cracks with non negligible dislocations (see the bonds $\boldsymbol{\xi }_{ij_{2}}$, $\boldsymbol{\xi }_{ij_{3}}$ in Fig. \ref{fig3_1}(c) and $\boldsymbol{\xi }_{ij_{5}}$, $\boldsymbol{\xi }_{ij_{7}}$ in Fig. \ref{fig3_1}(b)). 
The computational procedure to evaluate the
aperture value of NODE $\boldsymbol{i}$ is presented in Algorithm \ref{al_1}. The algorithm considers all broken bonds connected with NODE $\boldsymbol{i}$ and, for each of them, evaluates
the aperture in the initial direction of the bond. Then, the aperture of the
node is the average of bond apertures. Subsequently, the permeability and storage coefficients related to NODE $\boldsymbol{i}$ are updated correspondingly. Although this method may be improved, it is found in the subsequent applications to provide reasonable results.}


\begin{algorithm}[tph!]
	\caption{Routine for computing aperture value at the NODE $\boldsymbol{i}$.}
	\label{al_1}
	\LinesNumbered 
	\KwIn{$\left\Vert \boldsymbol{\xi}_{ij}\right\Vert ^{o}$, $\left\Vert \boldsymbol{\xi}_{ij}\right\Vert ^{d}$, $\beta _{ij}$}
	\KwOut{$a_i$}
	$\left\Vert \boldsymbol{\xi }_{ij}\right\Vert ^{o}$ is the original length of bonds $\boldsymbol{\xi }_{ij}$\;
	$\left\Vert \boldsymbol{\xi }_{ij}\right\Vert ^{d}$ is the deformed length of bonds $\boldsymbol{\xi }_{ij}$\;
	$\beta _{ij}$ is the angle between the initial bond and the deformed bond\;
	$a_i$ is the aperture value of NODE $\boldsymbol{i}$\; 
	$n$ is a cumulative indicator \;
	Initialisation of $n$: $n=0$, and $a_i$: $a_i=0$\;
	\For{all bonds}{
		\If{$\frac{\left\Vert \boldsymbol{\xi }_{ij}\right\Vert ^{d}-\left\Vert \boldsymbol{%
					\xi }_{ij}\right\Vert ^{o}}{\left\Vert \boldsymbol{\xi }_{ij}\right\Vert ^{o}}%
			\geqslant s_{c}$ and $(\left\Vert \boldsymbol{\xi }_{ij}\right\Vert ^{d}\cos \beta_{ij}) -\left\Vert \boldsymbol{%
				\xi }_{ij}\right\Vert ^{o}\geqslant 0$} {
			$a_i=a_i+(\left\Vert \boldsymbol{\xi }_{ij}\right\Vert ^{d}\cos \beta_{ij})-\left\Vert \boldsymbol{\xi }_{ij}\right\Vert ^{o}$\;
			$n=n+1$\;
		}
		
	}
	$a_i=a_i/n$\;
\end{algorithm}

\newpage

\subsection{Solution algorithms}

To ensure the correctness and accuracy of the implementation of the PD-based
hydraulic fracture model, transient consolidation and fluid flow problems
are studied first. Thus, besides the algorithm for the solution of hydraulic
fracturing problems, algorithms for transient consolidation analysis and
transient fluid flow problems are also presented.

\subsubsection{Algorithm for transient consolidation analysis}

In transient consolidation problems, the inertia term governing the motion
of the solid can be neglected \citep{lewis1998finite}, and Eqs.(\ref{3.1.1}
and \ref{3.1.12}) are combined and rewritten as: 
\begin{equation}
\left[ 
\begin{array}{cc}
\boldsymbol{0} & \boldsymbol{0} \\ 
\boldsymbol{Q}^{T} & \boldsymbol{S}%
\end{array}%
\right] \left[ 
\begin{array}{c}
\boldsymbol{\dot{u}} \\ 
\boldsymbol{\dot{p}}%
\end{array}%
\right] +\left[ 
\begin{array}{cc}
\boldsymbol{K}^{PD} & \boldsymbol{-Q}^{PD} \\ 
\boldsymbol{0} & \boldsymbol{H}%
\end{array}%
\right] \left[ 
\begin{array}{c}
\boldsymbol{u} \\ 
\boldsymbol{p}%
\end{array}%
\right] =\left[ 
\begin{array}{c}
\boldsymbol{f} \\ 
\boldsymbol{q}^{w}%
\end{array}%
\right]  \label{3.2.1}
\end{equation}

To solve the coupled system, the finite differences in time are used and the
following expression is reached \citep{lewis1998finite}:%
\begin{equation}
\begin{array}{r}
\left[ 
\begin{array}{cc}
\vartheta \boldsymbol{K}^{PD} & \boldsymbol{-}\vartheta \boldsymbol{Q}^{PD}
\\ 
\boldsymbol{Q}^{T} & \boldsymbol{S}+\vartheta \Delta t\boldsymbol{H}%
\end{array}%
\right] ^{n+\vartheta }\left[ 
\begin{array}{c}
\boldsymbol{u} \\ 
\boldsymbol{p}%
\end{array}%
\right] ^{n+1}=\left[ 
\begin{array}{cc}
\left( \vartheta -1\right) \boldsymbol{K}^{PD} & \left( 1\boldsymbol{-}%
\vartheta \right) \boldsymbol{Q}^{PD} \\ 
\boldsymbol{Q}^{T} & \boldsymbol{S}-\left( 1-\vartheta \right) \Delta t%
\boldsymbol{H}%
\end{array}%
\right] ^{n+\vartheta }\left[ 
\begin{array}{c}
\boldsymbol{u} \\ 
\boldsymbol{p}%
\end{array}%
\right] ^{n} \\ 
+\left[ 
\begin{array}{c}
\boldsymbol{f} \\ 
\Delta t\boldsymbol{q}^{w}%
\end{array}%
\right] ^{n+\vartheta }%
\end{array}
\label{3.2.2}
\end{equation}%
in which $0\leqslant \vartheta \leqslant 1$\ is the parameter for time
integration.

\subsubsection{Algorithm for transient fluid flow problem}

Neglecting the coupling term caused by volumetric change in Eq.(\ref{3.1.1}%
), the typical governing equation for fluid flow in saturated porous media
can be expressed as:%
\begin{equation}
\boldsymbol{S\dot{p}}+\boldsymbol{Hp=q}^{w}  \label{3.3.1}
\end{equation}

To solve this equation, the following implicit time integration iteration is
often used \citep{zienkiewicz2000finite}:%
\begin{equation}
\boldsymbol{p}^{n+1}=\left[ \boldsymbol{S}+\vartheta \Delta t\boldsymbol{H}%
\right] ^{-1}\left\{ \left[ \boldsymbol{S}-\left( 1-\vartheta \right) \Delta
t\boldsymbol{H}\right] \boldsymbol{p}^{n}-\Delta t\boldsymbol{q}^{w}\right\}
\label{3.3.2}
\end{equation}

In order to obtain a stable solution, $\vartheta $\ is usually taken as $%
0.5\leqslant \vartheta \leqslant 1$\ .

\subsubsection{Algorithm for dynamic solution of hydraulic fracturing
problems}

Two classes of algorithms, \textquotedblleft monolithic\textquotedblright\ and \textquotedblleft staggered\textquotedblright\, are usually
considered to solve the hydro-mechanical coupled system \citep{Walter2013}.
According to %
\citep{lewis1998finite,zienkiewicz1999computational,zienkiewicz2000finite,peruzzo2019dynamics}%
, Eqs.(\ref{3.1.1}) and (\ref{3.1.11}) are combined and rewritten as: 
\begin{equation}
\left[ 
\begin{array}{cc}
\boldsymbol{M}^{PD} & \boldsymbol{0} \\ 
\boldsymbol{0} & \boldsymbol{0}%
\end{array}%
\right] \left[ 
\begin{array}{c}
\boldsymbol{\ddot{u}} \\ 
\boldsymbol{\ddot{p}}%
\end{array}%
\right] +\left[ 
\begin{array}{cc}
\boldsymbol{0} & \boldsymbol{0} \\ 
\boldsymbol{Q}^{T} & \boldsymbol{S}%
\end{array}%
\right] \left[ 
\begin{array}{c}
\boldsymbol{\dot{u}} \\ 
\boldsymbol{\dot{p}}%
\end{array}%
\right] +\left[ 
\begin{array}{cc}
\boldsymbol{K}^{PD} & \boldsymbol{-Q}^{PD} \\ 
\boldsymbol{0} & \boldsymbol{H}%
\end{array}%
\right] \left[ 
\begin{array}{c}
\boldsymbol{u} \\ 
\boldsymbol{p}%
\end{array}%
\right] =\left[ 
\begin{array}{c}
\boldsymbol{f} \\ 
\boldsymbol{q}^{w}%
\end{array}%
\right]  \label{3.4.1}
\end{equation}%
which in concise form can be written as:%
\begin{equation}
\boldsymbol{M\ddot{a}}+\boldsymbol{C\dot{a}}+\boldsymbol{Ka}=\boldsymbol{f}
\label{3.4.2}
\end{equation}

The time integration of such system described in Eq.(\ref{3.4.2}) can be
carried out by using the \textquotedblleft Truncated Taylor series
collocation algorithm\textquotedblright\ presented in \textit{sect.} 18.3.3
of \citep{zienkiewicz2000finite}. A quadratic expansion is used (p=2). The
solution of the unknown variable $\boldsymbol{a=}\left[ \boldsymbol{u},%
\boldsymbol{p}\right] ^{T}$ at time step $n+1$ is given by:%
\begin{equation}
\boldsymbol{a}^{n+1}=-\boldsymbol{A}^{-1}(\boldsymbol{f}^{n+1}+\boldsymbol{C%
\dot{a}}^{n+1}+\boldsymbol{M\ddot{a}}^{n+1})  \label{3.4.3}
\end{equation}%
where

\begin{equation}
\begin{array}{l}
\boldsymbol{A}=\frac{2}{\beta _{2}\Delta t^{2}}\boldsymbol{M+}\frac{2\beta
_{1}}{\beta _{2}\Delta t}\boldsymbol{C+K} \\ 
\boldsymbol{\dot{a}}^{n+1}=-\frac{2\beta _{1}}{\beta _{2}\Delta t}%
\boldsymbol{a}^{n}+\left( 1-\frac{2\beta _{1}}{\beta _{2}}\right) 
\boldsymbol{\dot{a}}^{n}+\left( 1-\frac{\beta _{1}}{\beta _{2}}\right) 
\boldsymbol{\ddot{a}}^{n} \\ 
\boldsymbol{\ddot{a}}^{n+1}=-\frac{2}{\beta _{2}\Delta t^{2}}\boldsymbol{a}%
^{n}-\frac{2}{\beta _{2}\Delta t}\boldsymbol{\dot{a}}^{n}-\frac{1-\beta _{2}%
}{\beta _{2}}\boldsymbol{\ddot{a}}^{n}%
\end{array}
\label{3.4.4}
\end{equation}%
in which $\beta _{2}\geqslant \beta _{1}\geqslant 0.5$ are parameters for
the time integration.

\textcolor{blue}{However, when solving hydraulic fracturing problems by means of a \textquotedblleft monolithic\textquotedblright\ approach,} the system matrix
needs to be updated in each step, and in addition, the combined system
matrices will also be large, which will render the problem memory- and
time-intensive. 
\textcolor{blue}{Therefore, another alternative
option, the \textquotedblleft staggered approach\textquotedblright, is adopted here. In \citep{zhou2019phase}, a \textquotedblleft staggered coupling scheme\textquotedblright\ is implemented in COMSOL to solve the HF problems by using the phase-field model for saturated porous media. In each time step, the following solving scheme is adopted: the fields of pore pressure and displacement are calculated at first, then the historic strain field is updated according to the output pore pressure and displacement fields, whereafter the damage field is solved by using the fields already obtained. In addition, the Newton-Raphson algorithm is adopted in each staggered solution sequence, and a global relative error is estimated to determine whether to start a new iteration step or perform next time step. However, although the \textquotedblleft staggered coupling scheme\textquotedblright\ in \citep{zhou2019phase} can ensure the convergence of the coupled system, it seems to lose the ability to express some peculiar phenomena of the pore pressure, especially in dynamics. In this paper, at the expense of some robustness, a modified \textquotedblleft staggered approach\textquotedblright \ is adopted to capture the hydrodynamic phenomena such as pore pressure oscillations in the hydraulic fracturing process.}
In the adopted \textquotedblleft staggered approach\textquotedblright\ , the
coupled system is divided into two parts , hydraulic diffusion and
hydro-driven deformation of the porous solid %
\citep{yang2018hydraulic,wu2014extension,peng2018hydraulic}, which means
that the first row and second row of Eq.(\ref{3.4.1}) are solved
sequentially, and the previously solved results of $\boldsymbol{u}$ and $%
\boldsymbol{p}$ are used to evaluate the fluid volumetric source terms or
nodal forces for the next solving sequence. In each solution sequence, there
are two steps:

--step 1: solve the pressure field ($\boldsymbol{p}^{n+1}$) using the
following implicit time integration iteration:%
\begin{equation}
\boldsymbol{p}^{n+1}=\left[ \boldsymbol{S}+\vartheta \Delta t\boldsymbol{H}%
\right] ^{-1}\left\{ \left[ \boldsymbol{S}-\left( 1-\vartheta \right) \Delta
t\boldsymbol{H}\right] \boldsymbol{p}^{n}-\Delta t\boldsymbol{q}^{w}+%
\boldsymbol{Q}^{T}\left( \mathbf{u}^{n}-\mathbf{u}^{n-1}\right) \right\}
\label{3.4.5}
\end{equation}

--step 2: solve the displacement field ($\boldsymbol{u}^{n+1}$) of the solid
domain using the adaptive dynamic relaxation method presented in %
\citep{Underwood1983dynamic,Ni2019Coupling}.

\begin{figure}[h]
\begin{center}
\includegraphics[scale=0.55]{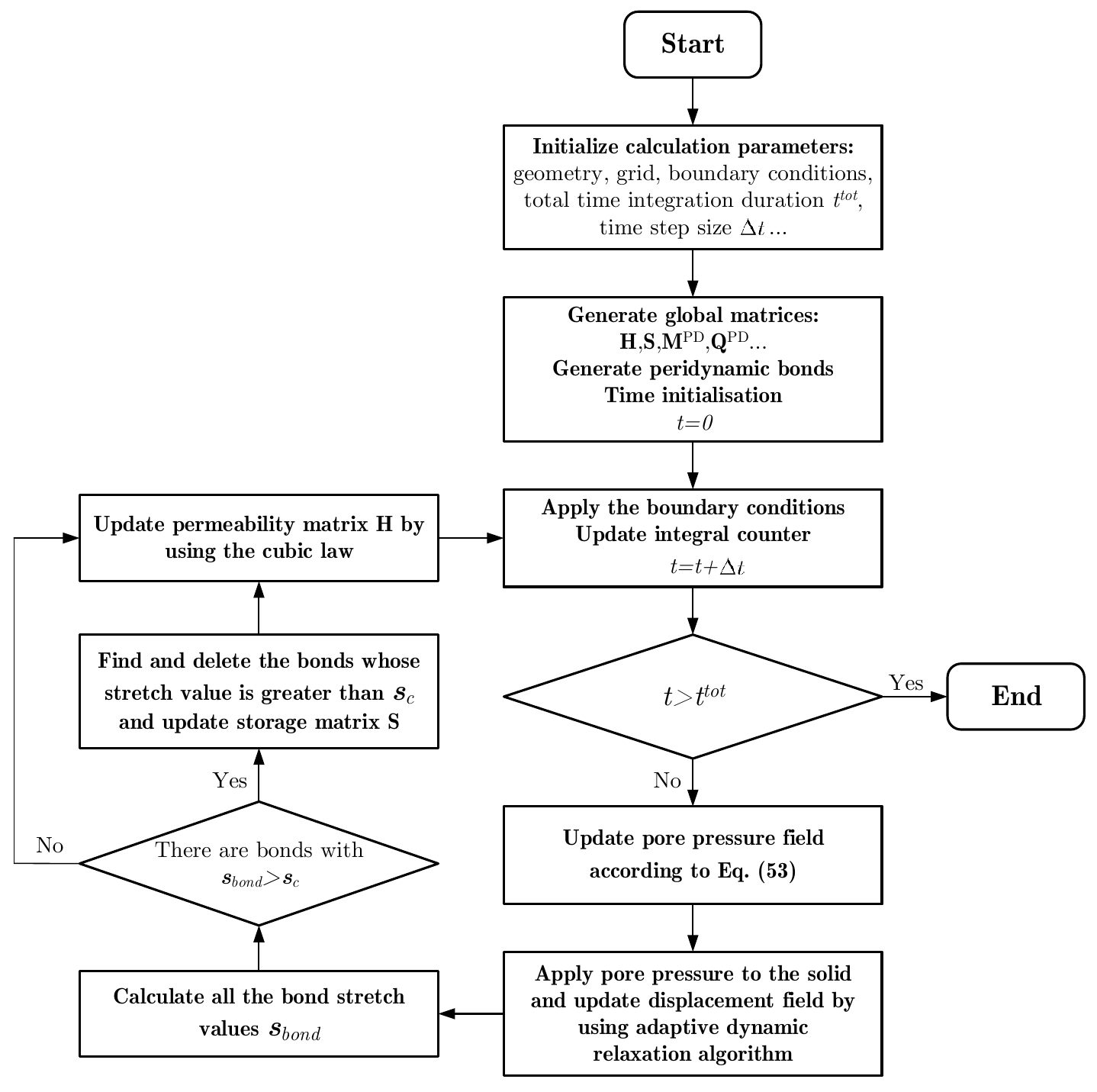} 
\end{center}
\caption{\textcolor{blue}{Flow chart of the solution algorithm for hydraulic
fracture problems.}}
\label{fig3}
\end{figure}

In addition, if there are initial or propagating cracks in the porous media,
the permeability and storage matrices ($\boldsymbol{H}$ and $\boldsymbol{S}$%
) need to be updated accordingly in each time step. \textcolor{blue}{The
flow chart of the solution algorithm by using the described \textquotedblleft staggered
approach\textquotedblright \  for HF problems is shown in Fig.\ref{fig3}.} 
\newpage

\section{Numerical examples}

\begin{figure}[thp!]
\begin{center}
\subfloat[Damage contour at the crack
tip]{\includegraphics[scale=0.3]{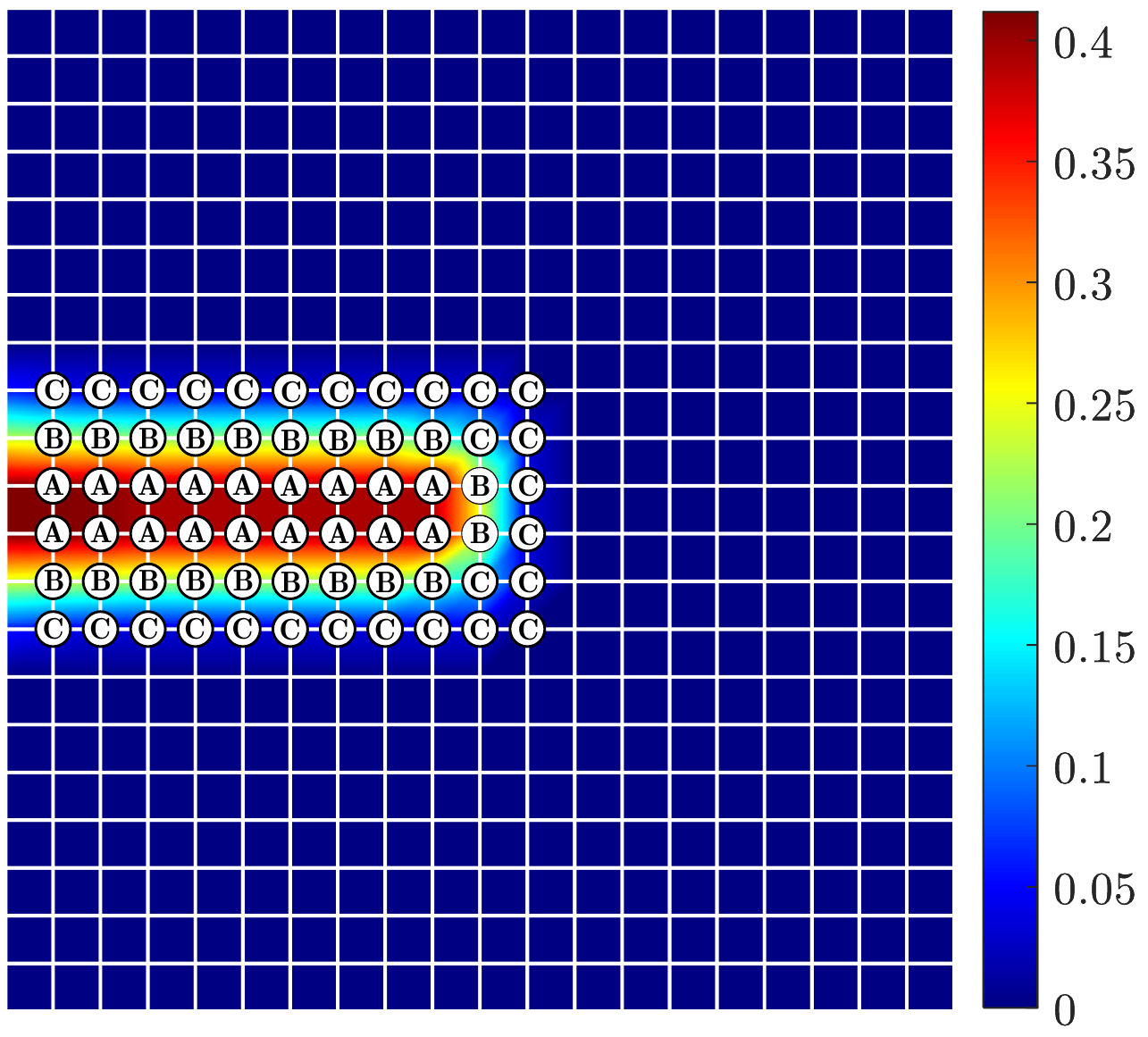}%
\label{fig4_1_1:sub1}} \subfloat[Isogram of the damage
contour]{\includegraphics[scale=0.675]{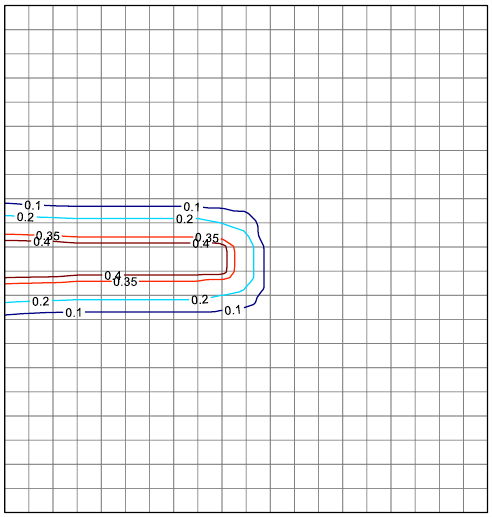}%
\label{fig4_1_1:sub2}} 
\end{center}
\caption{\textcolor{blue}{Damage value distribution at the crack tip in a
discrete peridynamic model.}}
\label{fig4_1_1}
\end{figure}

In this section, several numerical examples are presented to demonstrate the
capabilities of the proposed method. All examples are performed in plane
strain conditions. The $m$-ratio is always taken as $m_{ratio}=3 $. %
\textcolor{blue}{First we show in Fig. \ref{fig4_1_1} the distribution of
damage value at the crack tip in a discrete peridynamic model. As shown in
Fig. \ref{fig4_1_1:sub1}, the nodes where damage values are greater than 0
can be divided into three types A, B and C. The damage value at nodes of
type A is greater than 0.4, see Fig. \ref{fig4_1_1:sub2}, while the value at
nodes of type B is greater than 0.2. According to this characteristic, the
threshold values are set as $c_{1}=0.2$ and $c_{2}=0.35$ in the linear
indicator functions. Only one layer of elements would be identified in the
fracture domain, which can minimize the influence sphere of the crack on the
permeability parameters of elements. The smooth transition from fracture
domain to reservoir domain is guaranteed by identifying the second layer of
nodes in the transition domain. In addition, $c_{1}$  can also take other
values smaller than 0.2.}

\newpage

\subsection{One-dimensional consolidation problem}

\begin{figure}[h]
\begin{center}
\includegraphics[scale=0.6]{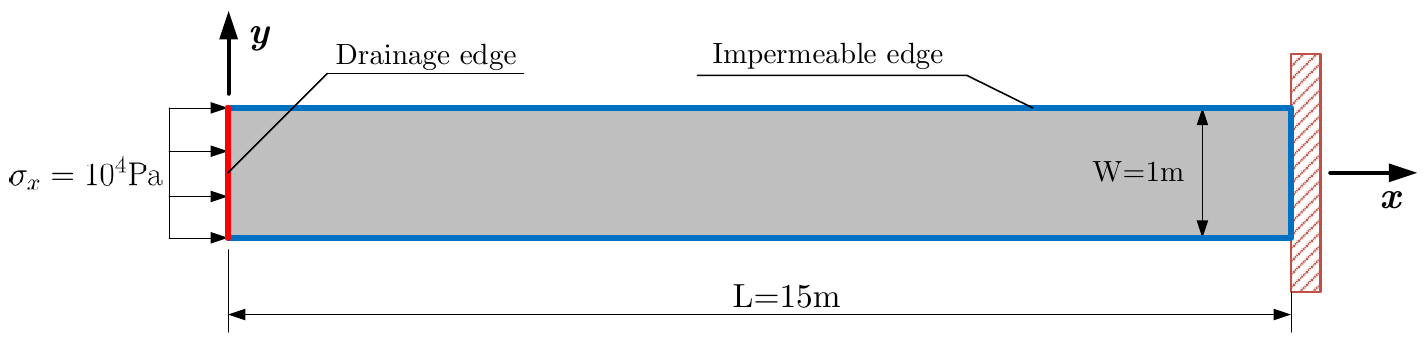}
\end{center}
\caption{Geometry and boundary conditions of the consolidation problem.}
\label{fig4}
\end{figure}

We solve an one-dimensional consolidation problem presented in %
\citep{zhang2019coupling}. The material parameters are, Young modulus: $%
E=10^{8}Pa$, Poisson's ratio: $\nu =0$, Biot constant: $\alpha =0.5$,
permeability coefficient: $k=10^{-12}m^{2}$, fluid viscosity: $\mu
^{w}=10^{-3}Pa\cdot s$, storage coefficient: $S=1/\left( 6.06\times
10^{9}Pa\right) $. Fig. \ref{fig4} shows the geometry and boundary
conditions. The analytical solutions of this problem are given as %
\citep{zhang2019coupling,wang2017theory}:%
\begin{equation}
p\left( x,t\right) =\frac{4vP_{0}}{\pi }\sum\limits_{m}^{N}\left\{ \frac{1}{%
2m+1}\exp \left( -\left( \frac{\left( 2m+1\right) \pi }{2L}\right)
^{2}ct\right) \times \sin \left( \frac{\left( 2m+1\right) \pi x}{2L}\right)
\right\}  \label{4.1}
\end{equation}%
\begin{equation}
\begin{array}{r}
u\left( x,t\right) =c_{m}vP_{0}\left\{ L-x-\frac{8L}{\pi ^{2}}%
\sum\limits_{m}^{N}\left\{ \frac{1}{\left( 2m+1\right) ^{2}}\exp \left(
-\left( \frac{\left( 2m+1\right) \pi }{2L}\right) ^{2}ct\right) \times \cos
\left( \frac{\left( 2m+1\right) \pi x}{2L}\right) \right\} \right\} \\ 
+bP_{0}\left( L-x\right)%
\end{array}
\label{4.2}
\end{equation}%
in which $a=10^{-8}Pa,$ $b=\frac{a}{1+a\alpha ^{2}/S},$ $v=\frac{a-b}{%
a\alpha },$ $c=\frac{\kappa }{\left( a\alpha ^{2}+S\right) \mu _{w}}$ and $%
c_{m}=\frac{a-b}{v}$.

\begin{figure}[h]
\begin{center}
\includegraphics[scale=0.8]{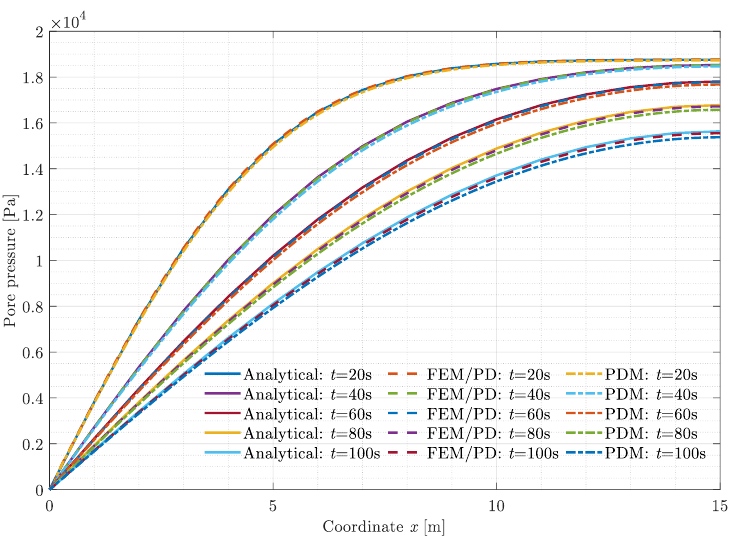}
\end{center}
\caption{\textcolor{blue}{Comparison of analytical and PD pore pressure
solutions along the
central axis.}}
\label{fig4_1}
\end{figure}

\begin{figure}[h!]
\begin{center}
\includegraphics[scale=0.8]{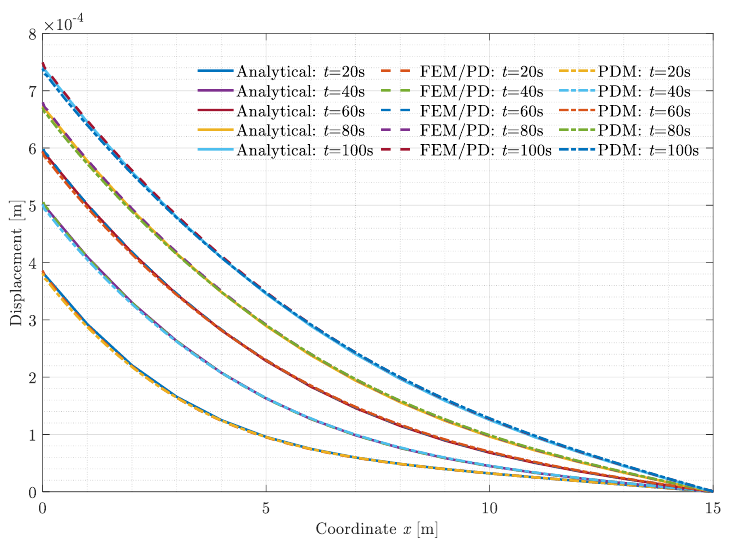}
\end{center}
\caption{\textcolor{blue}{Comparison of analytical and PD displacement
solutions along the
central axis.}}
\label{fig4_2}
\end{figure}

The discretization parameters are horizon $\delta =0.15m$, and the
corresponding grid spacing is $\Delta x=\delta /m_{ratio}=0.05m$. The
solution algorithm in $sect.$ $3.4.1$ is used to solve this coupled system,
and the time step is taken as $\Delta t=1s$. Figs. \ref{fig4_1} and \ref%
{fig4_2} show the distributions of the pore pressure and displacement along
the central axis at different times: $20s$, $40s$, $60s$, $80s$ and $100s$. 
\textcolor{blue}{The same problem has also solved by using a coupled peridynamic-only model presented in \citep{zhang2019coupling}. The numerical and analytical solutions are plotted in Figs. \ref{fig4_1} and \ref{fig4_2}. Tab. \ref{tab1} shows the computing costs of FEM/PD and PD-only solutions (PDM). 
The comparisons in Figs. \ref{fig4_1} and \ref{fig4_2} indicate that the
results obtained by the proposed method are in excellent agreement with the
analytical solutions, and the FEM/PD model has a higher efficiency than a PD-only model.}

\textcolor{blue}{
\begin{table}[h!]
	\caption{Computing costs of the consolidation problem by using different methods.}
	\label{tab1}
	\centering  
	{\scriptsize
		\begin{tabular}{p{2.5cm}<{\centering}p{3cm}<{\centering}p{3cm}<{\centering}}
			\toprule Methods & FEM/PD & PD-only \\ \hline
			CPU time [s] & 37.69 & 62.85 \\ 
			\bottomrule &  & 
		\end{tabular}}
\end{table}}

\subsection{Pressure distribution in a single crack}

\begin{figure}[h!]
\begin{center}
\includegraphics[scale=0.6]{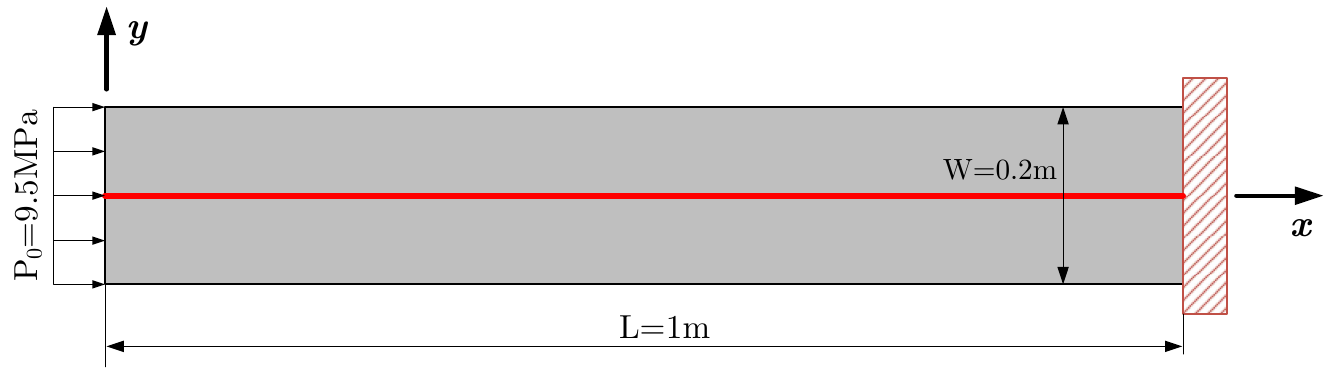}
\end{center}
\caption{Geometry and boundary conditions for a rock sample with a single
crack.}
\label{fig4_3_1}
\end{figure}

In this section, an example of variation of hydraulic pressure distribution
with time in a single crack \citep{zhou2019phase,yang2018hydraulic} is
studied by using the present method. Geometry and boundary conditions are
shown in Fig. \ref{fig4_3_1}.

A hydraulic pressure $P_{0}=9.5MPa$ is applied suddenly on the left edge,
while the other edges are impermeable. This is an "unsteady" or evolving
problem, and its analytical solution available in \citep{zhou2019phase} is
given as:%
\begin{equation}
\frac{P}{P_{0}}=1+\frac{4}{\pi }\sum\limits_{n=0}^{\infty }\left[ \exp
\left( -(2n+1)^{2}(T/4)\pi ^{2}\right) \cos \left( \frac{(2n+1)\pi }{2}\zeta
\right) \left( \frac{(-1)^{n+1}}{2n+1}\right) \right]  \label{4.3}
\end{equation}%
in which $\zeta =\frac{L-x}{L}$, $P$ is the pressure value at $x$, $L$ is
the length of the specimen, and $T_{d}$ is a dimensionless time, defined as:%
\begin{equation}
T_{d}=K_{w}\frac{\left( a^{2}/(12\mu _{w})\right) t}{L^{2}}  \label{4.4}
\end{equation}%
where $K_{w}$ and $\mu _{w}$ are the bulk modulus and viscosity coefficient
of the fluid, and $a$ is the aperture of the central crack. To compare with
the analytical solution, we only solve the fluid flow part by using the
solution algorithm in $sect.3.4.2$. The parameters associated to the fluid
flow are taken as: $K_{w}=2.2GPa$, $\mu_{w}=10^{-3}Pa\cdot s$, $\alpha =1$, $%
n=2\times 10^{-5}$ and $a=3\times 10^{-5}m$. Two sets of discretization
parameters have been used as a comparison: horizon $\delta _{1}=15mm$ and $%
\delta _{2}=30mm$, the corresponding grid spacings are $\Delta x=\delta
_{1}/m_{ratio}=5mm$ and $\Delta x=\delta _{2}/m_{ratio}=10mm$. The initial
crack is described by the peridynamic portion and the FEM model for the
calculation of the fluid flow is modified accordingly. The time integration
parameters are taken as: $\Delta t=2\times 10^{-8}s$ and $\vartheta =0.5$.

\begin{figure}[h]
\centering  
\subfloat[$T_{d}=0.1$.]{		%
\includegraphics[scale=0.6]{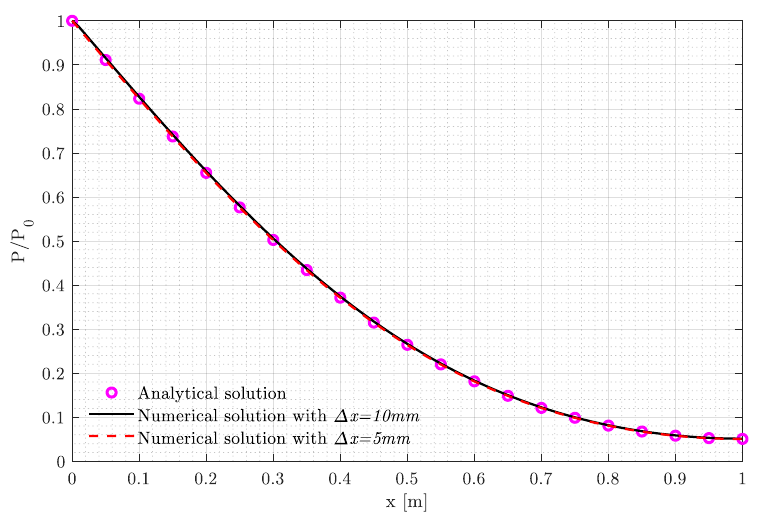}		%
\label{fig4_3:sub1}} \subfloat[$T_{d}=0.2$.]{		%
\includegraphics[scale=0.6]{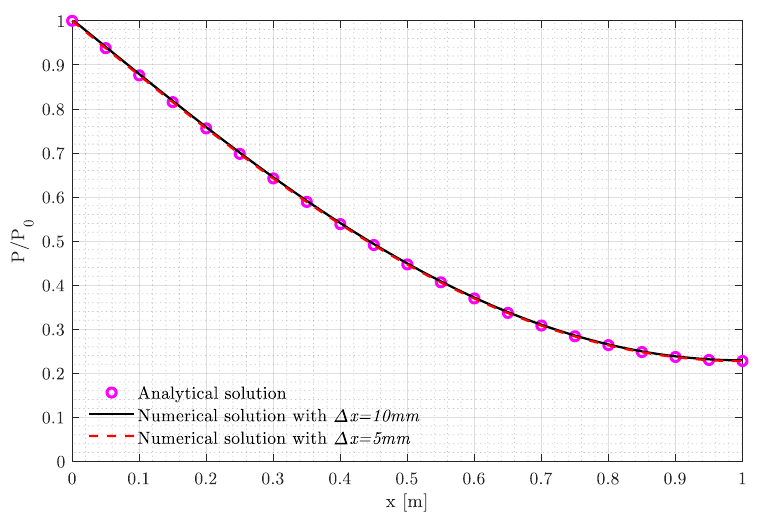}		%
\label{fig4_3:sub2}}\newline
\subfloat[$T_{d}=0.3$.]{		%
\includegraphics[scale=0.6]{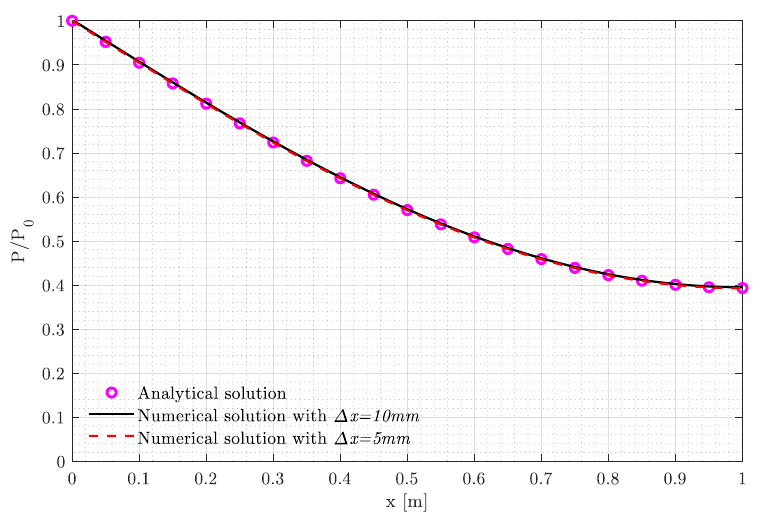}		%
\label{fig4_3:sub3}} \subfloat[$T_{d}=0.5$.]{		%
\includegraphics[scale=0.6]{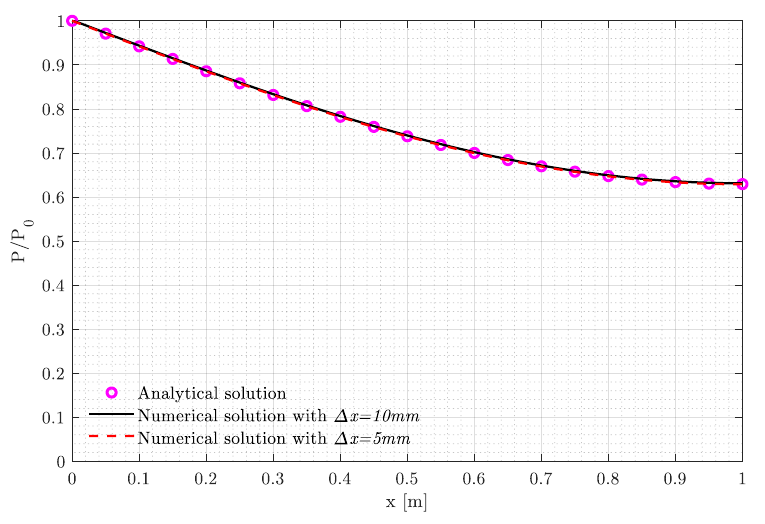}		%
\label{fig4_3:sub4}}
\caption{Comparison of the numerical and analytical pressure distribution
along the crack.}
\label{fig4_3}
\end{figure}

\begin{figure}[tph!]
\begin{center}
\includegraphics[scale=0.8]{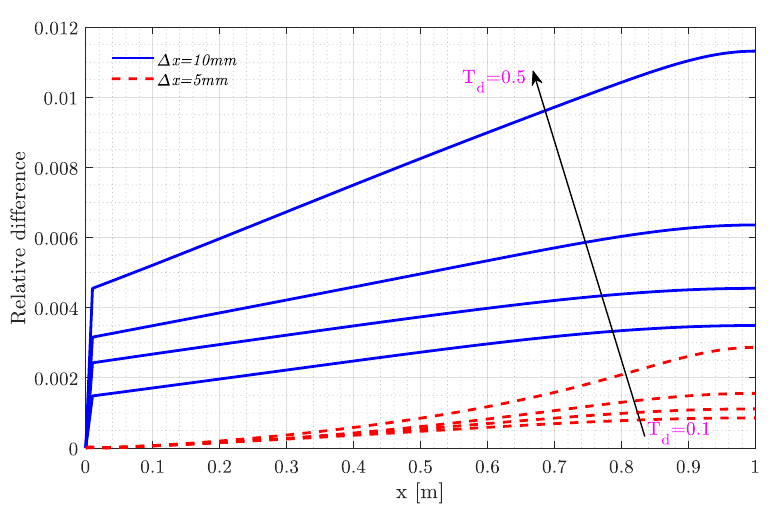}
\end{center}
\caption{The relative differences between numerical and analytical solutions
of the pressure distribution along the single crack.}
\label{fig4_42_1}
\end{figure}

Figs. \ref{fig4_3:sub1} to \ref{fig4_3:sub4} show the comparison between the
numerical and analytical solutions of the pressure distributions along the
initial crack at different times $T_{d}$. The relative differences between
the numerical and analytical solutions are plotted in Fig. \ref{fig4_42_1}.
In general, the numerical solutions are very close to the analytical ones
and a denser grid generate a more accurate solution. 


\newpage

\textcolor{blue}{\subsection{Pressure-driven fracture propagation: a centrally notched
specimen subjected to an increasing internal pore pressure}}

\label{pressure_drive}

\begin{figure}[h]
\begin{center}
\includegraphics[scale=0.5]{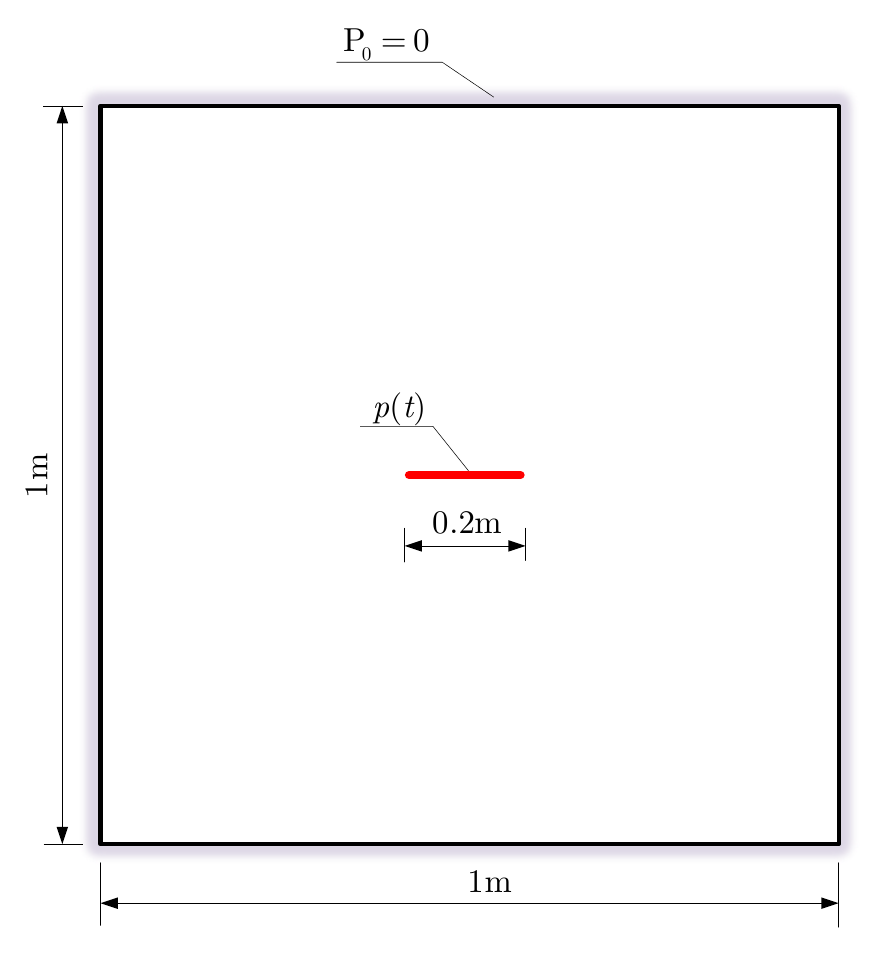}
\end{center}
\caption{Geometry and boundary conditions of the notched specimen subjected
to internal pressure.}
\label{fig4_43}
\end{figure}

In this section, a pressure-driven fracture propagation example presented in %
\citep{zhou2018phase} is solved by applying the present method. Geometry and
boundary conditions are shown in Fig. \ref{fig4_43}. An increasing pressure
is applied on the surface of the initial crack to deform the specimen.
Assume that the length of the initial crack is $2l_{c}$, the vertical
displacement of the crack surface driven by the pore pressure can be given
analytically as \citep{sneddon1969crack,zhou2018phase}:%
\begin{equation}
u\left( x,p\right) =\frac{2pl_{c}}{E_{p}}\left( 1-\frac{x^{2}}{l_{c}^{2}}%
\right) ^{1/2}  \label{4.5}
\end{equation}%
in which $E_{p}=E/(1-\upsilon ^{2})$ is the plane strain Young's modulus and 
$E$\ and $\upsilon $\ are the Young's modulus and Poisson's ratio,
respectively. 

The mechanical and fluid parameters used in the calculation are, Young's
modulus: $E=210GPa$, Poisson's ratio: $\upsilon =0.3$, critical energy
release rate: $G_{c}=2700J/m^{2}$, porosity: $n_{r}=0.002$, mass density: $%
\rho _{r}=\rho _{f}=1000kg/m^{3}$, biot constant: $\alpha =1$, bulk modulus
and viscosity coefficient of the fluid: $K_{w}=10^{8}Pa$ and $\mu
_{w}=10^{-3}Pa\cdot s$, permeability coefficient of the reservoir domain: $%
k_{r}=10^{-15}m^{2}$.

\begin{figure}[h]
\begin{center}
\includegraphics[scale=0.8]{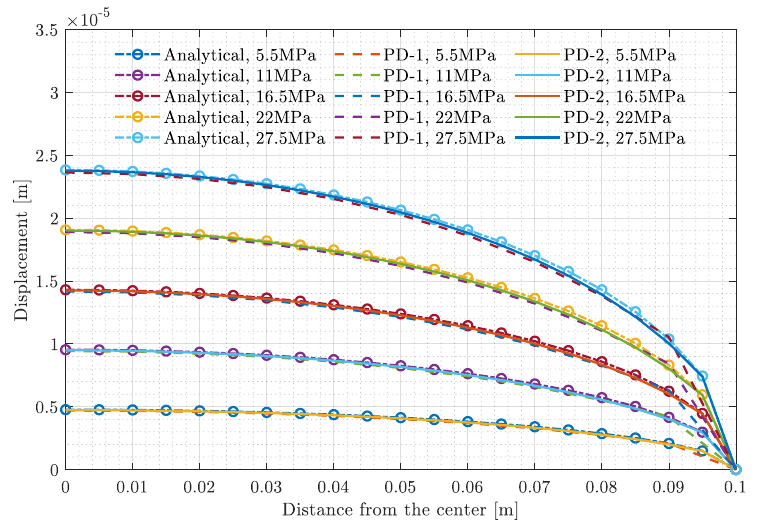}
\end{center}
\caption{Comparison of the displacement along the notch between the present
results and the analytical solution.}
\label{fig4_4_1}
\end{figure}

\begin{figure}[h!]
\centering  
\centering  
\subfloat[Damage level at
$2154.0s$.]{\includegraphics[scale=0.4]{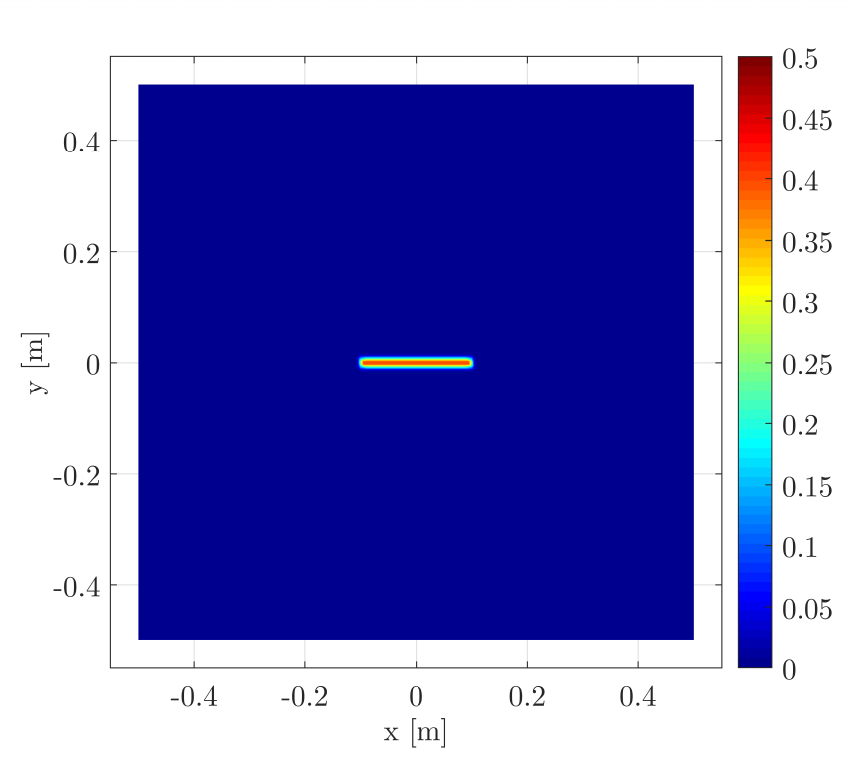}\label{fig4_4:sub1}}
\subfloat[Damage level at
$2154.4s$.]{\includegraphics[scale=0.4]{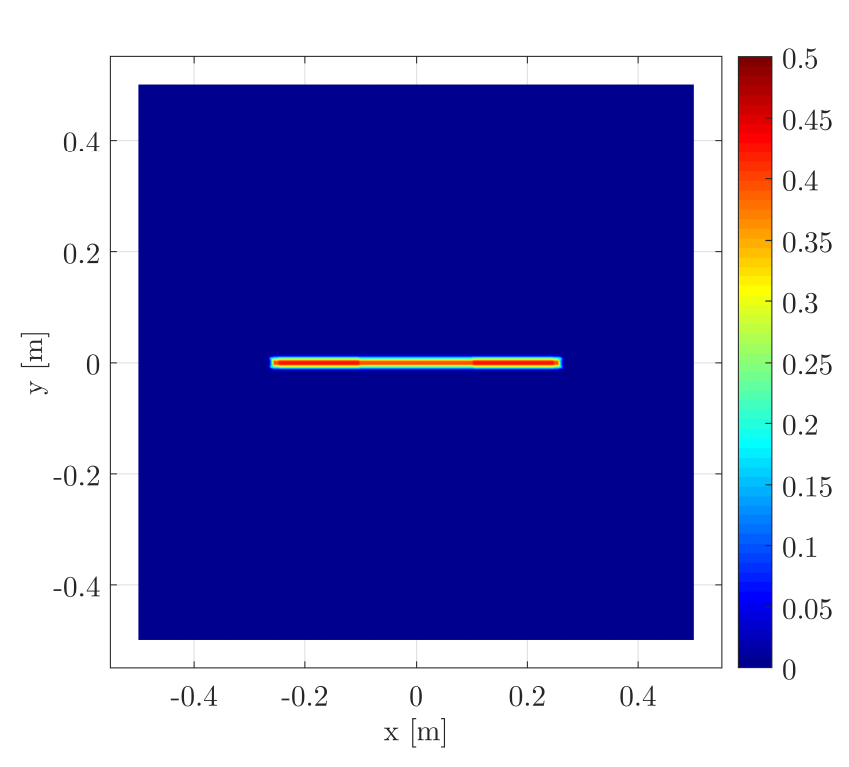}\label{fig4_4:sub3}}
\subfloat[Damage level at
$2154.8s$.]{\includegraphics[scale=0.4]{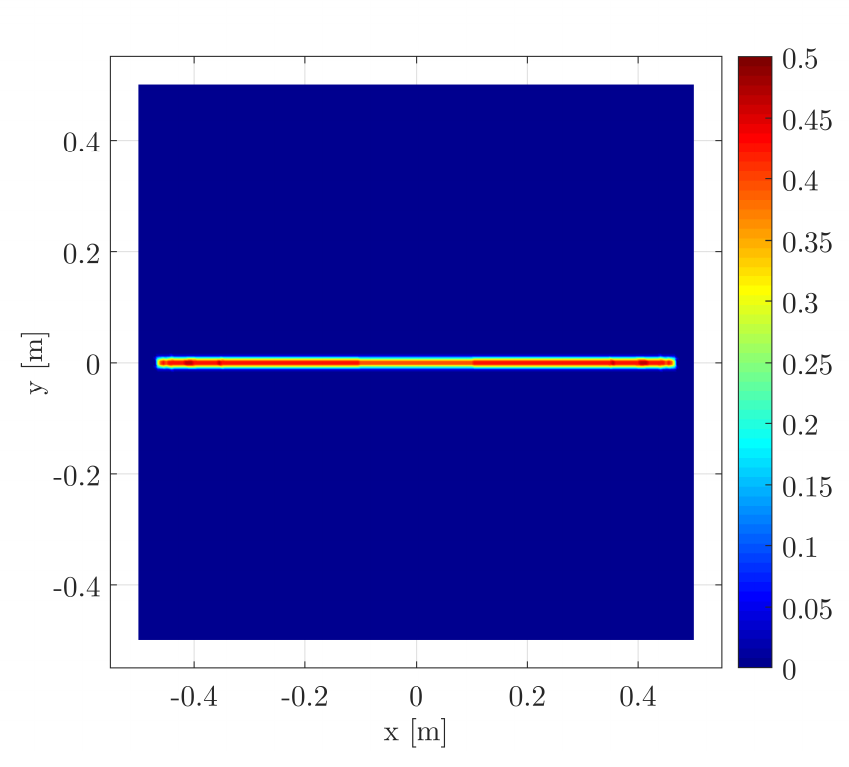}\label{fig4_4:sub5}}
\newline
\subfloat[Pressure distribution at
$2154.0s$.]{\includegraphics[scale=0.4]{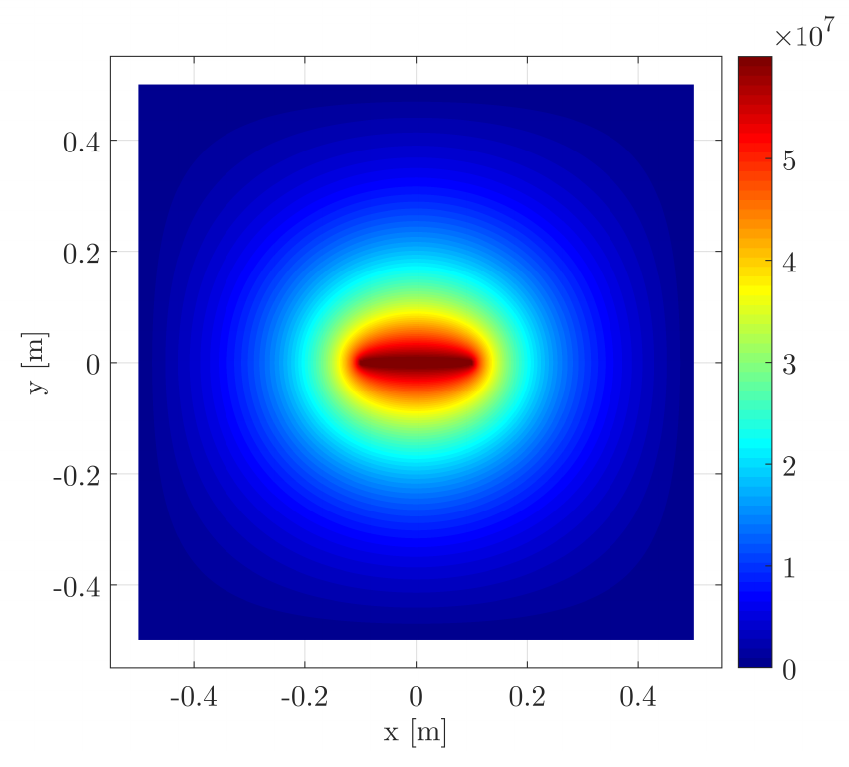}\label{fig4_4:sub2}}
\subfloat[Pressure distribution at
$2154.4s$.]{\includegraphics[scale=0.4]{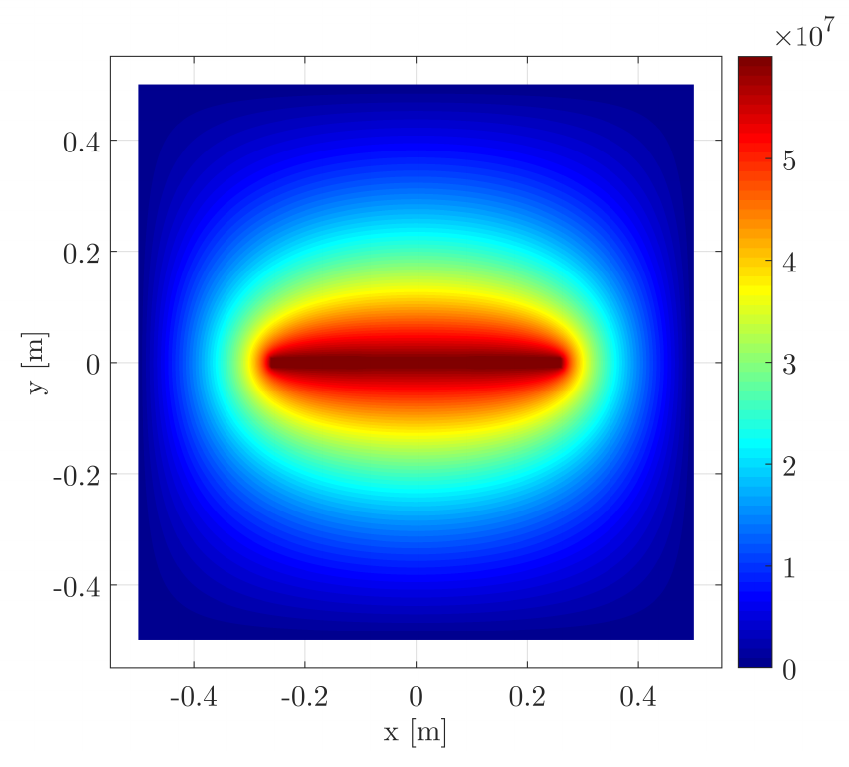}\label{fig4_4:sub4}}
\subfloat[Pressure distribution at
$2154.8s$.]{\includegraphics[scale=0.4]{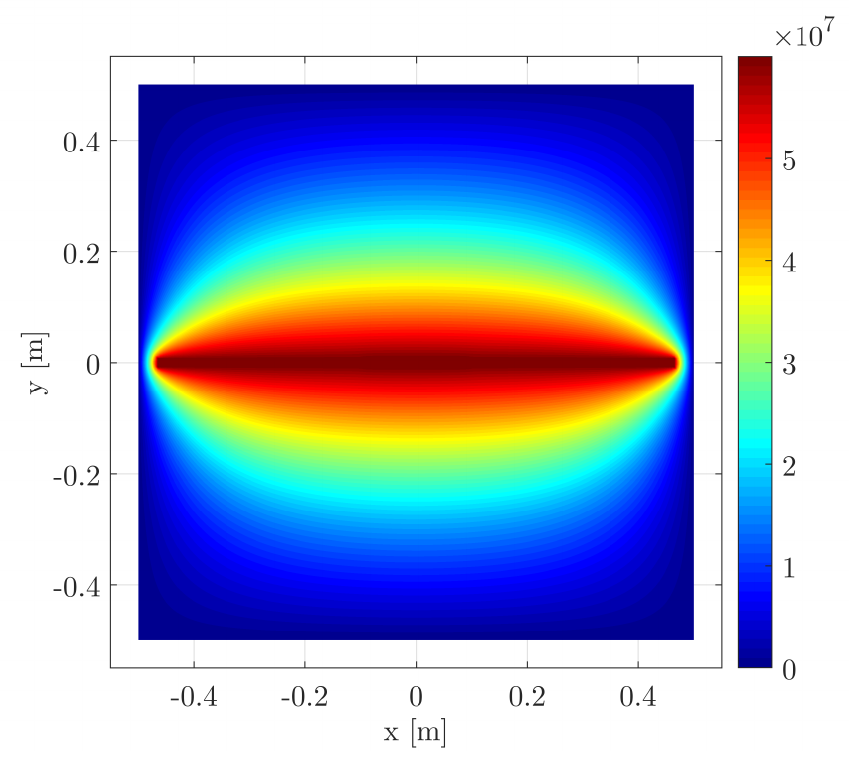}\label{fig4_4:sub6}}
\caption{Crack propagation and the variation of pressure distribution at
different time steps.}
\label{fig4_4}
\end{figure}

The whole domain is discretized with uniform quadrilateral grid. The FE
nodes of the fluid mesh and the PD nodes of the solid grid have the same
coordinates. Two sets of discretization parameters are taken for comparison:
Mesh $1$ with a size of $\Delta x_{1}=10^{-2}m$\ and Mesh $2$ with a size of 
$\Delta x_{2}=5\times 10^{-3}m$. The corresponding\ horizons are $\delta
_{1}=\Delta x_{1}\times m_{ratio}=3\times 10^{-2}m$ and $\delta _{2}=\Delta
x_{2}\times m_{ratio}=1.5\times 10^{-3}m$. As in \citep{zhou2018phase}, the
pressure leading to crack propagation with the adopted parameters is about $%
60MPa$. Thus, the final pressure $60.5MPa$ is applied on the crack surface
in $2200$ uniform time steps, and each time step is divided into $2500$
small incremental time steps with a size of $\Delta t=4\times 10^{-4}s$.
Fig. \ref{fig4_4_1} shows the comparison between the displacement along the
notch of the present results and the analytical solution. The numerical
results obtained with the denser mesh are closer to the analytical solutions
than those obtained with sparser mesh, and they are all in good agreement.
When the applied pressure reaches $59.235MPa$ at $2154s$, the crack starts
to propagate. The variations of the crack pattern and the distribution of
pressure are shown in Figs. \ref{fig4_4:sub1} to \ref{fig4_4:sub6}. 
\newpage

\textcolor{blue}{\subsection{Fluid-driven fracture propagation: reproducing the results of the KGD model}
In this example we deal with fluid driven fracture initiation and propagation. This is completely different from the case of $sect.$ \ref{pressure_drive} where a pressure driven fracture is investigated. In that case the driving force is applied directly to the momentum balance equation of the mixture, Eq.(\ref{1.7}), via the pressure term in the effective stress principle, Eq.(\ref{1.14}), while in case of fluid driven fracture the loading term is applied to the fluid mass balance equation Eq.(\ref{2.2.3}). Its effect is then transmitted to the momentum balance equation via the coupling terms. The loading cases for instance determine whether a power law behavior exists or is destroyed. This was clearly evidenced in \citep{milanese2016avalanches}. In general the  comparison of results for crack propagation between numerical models and exact solutions is not an easy task in fluid driven fracture. The reason is that the respective methods and underlying assumptions strongly impact the results. Exact solutions are based on many simplifying assumptions which are not always easy to implement in numerical methods based on more advanced approaches as explained below. Typical such assumptions are zero fluid lag, flow occurring predominantly in the fracture plane, finite and non-zero propagation speed, linear elastic fracture mechanics. They result in a steady propagation of a semi-infinite crack with constant or smoothly changing tip velocity \cite{linkov2014universal,bunger2005toughness}. As far as LEFM is concerned,  De Pater \citep{de2015hydraulic} states that \textquotedblleft carefully scaled laboratory studies have revealed that conventional models for fracture propagation which are based on LEFM cannot adequately describe propagation. Instead, propagation should at least be based on tensile failure over a cohesive zone\textquotedblright.
}

\textcolor{blue}{Boone and Ingraffea \citep{boone1990numerical} presented a
FEM/FD model based on LEFM with a moving crack along a path known a priori.
They assume negligible fracture toughness or fracture energy; further their
procedure does not allow to specify flow at the crack mouth hence some
iterations are needed in this case. These authors compare their solution
with the analytical solution of Spence and Sharp \citep{spence1985self}
based on an approximation to an asymptotic solution; in particular the fluid
is incompressible, the fracture impermeable and fluid lag in the tip region
is excluded. Secchi et al. \citep{secchi2007mesh} with their FEM model
based on cohesive tractions used the same analytical solution for a
comparison, but had to introduce LEFM to obtain favourable results. They
recovered however the stepwise advancement of the fracture tip, obviously
not present in the analytical solution. Some results are shown below for
comparison.}

\textcolor{blue}{Baykin and Golovin \citep{baykin2018application} reduced
for comparison reason their cohesive model to an elastic one with Biot
coefficient $\alpha=0$. Esfahani and Gracie \citep{parchei2019undrained} use
XFEM with a cohesive model and drained/undrained split. They validate their
model with respect to the solution for the problem of a planar hydraulic
fracture (HF) that propagates under plane strain condition, presented by
Geertsma and De Klerk \citep{geertsma1969rapid}, known as the KGD HF model.
As Secchi et al. \citep{secchi2007mesh} they find that the KGD solution
cannot be directly compared with the results of a cohesive HF model. The
reason is that the KGD model is independent of the fracture energy (or,
equivalently, fracture toughness), whereas the solution of a cohesive HF
model strongly depends on the value of the fracture energy and the length of
the cohesive zone. The KGD solution represents the propagation of a plane
strain HF in a viscosity dominated regime, whereas the solution of a
cohesive HF gets close to a viscosity dominated regime only when its
cohesive zone vanishes (or when the fracture toughness tends to zero).
Therefore, the KGD and cohesive models are only comparable when the fracture
propagates to such an extent that the cohesive zone is a negligible portion
of the total fracture length, or when the fracture energy (fracture
toughness) tends to zero. Finally, Chen et al. \citep{chen2020explicit} use
a Boundary Element Model with LEFM for modeling planar 3D hydraulic fracture
growth. Their results compare well with analytical solutions but evidence
stepwise fracture advancement as opposed to smooth behaviour.}

\begin{figure}[h!]
\centering  
\includegraphics[scale=0.65]{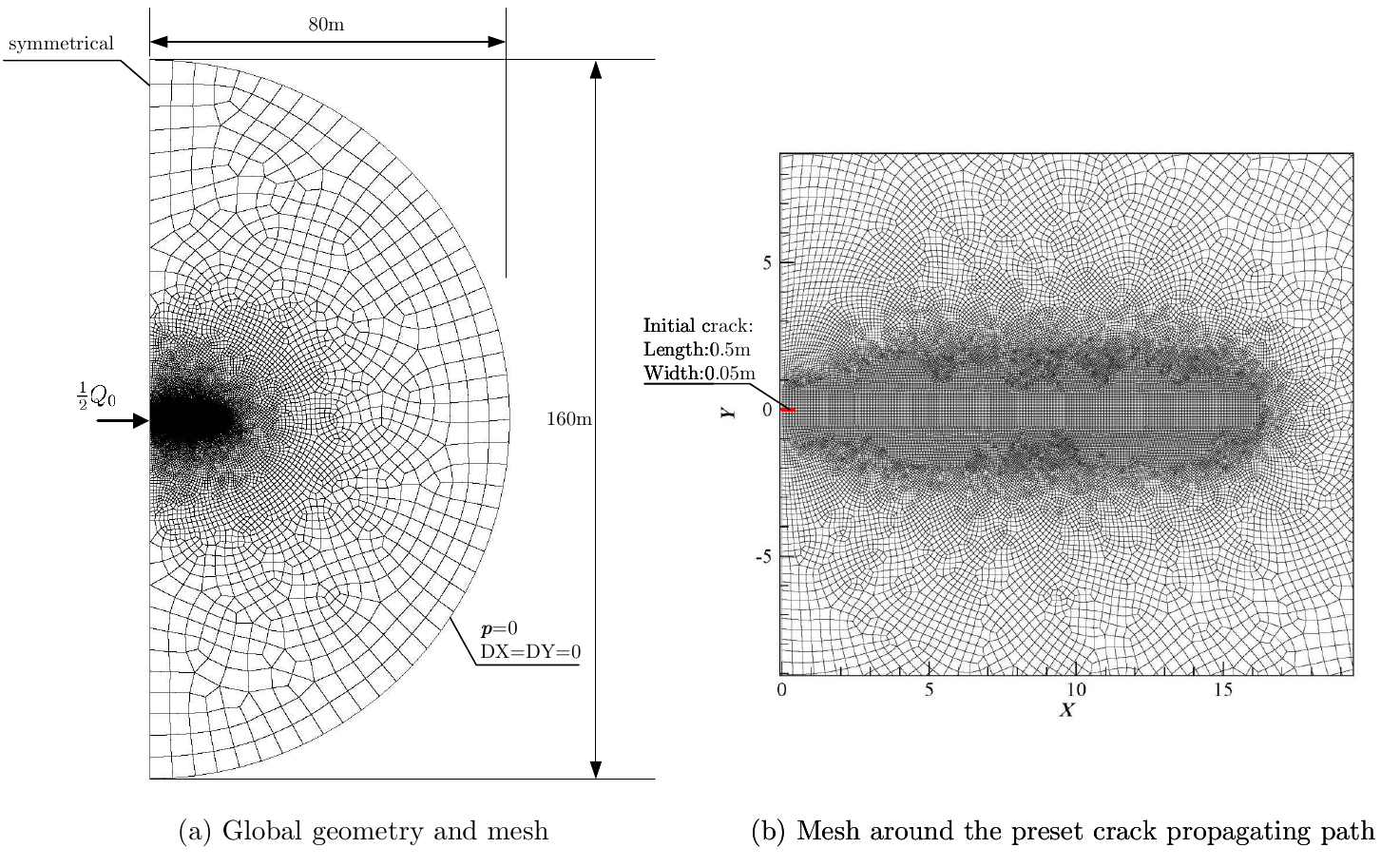}  
\caption{\textcolor{blue}{Diagrammatic sketch of the KGD model.}}
\label{fig4_22}
\end{figure}

\textcolor{blue}{To assess the accuracy of the proposed method in case of
fluid driven crack propagation we have simulated the case of Spence and
Sharp \citep{spence1985self}, keeping in mind that different methods give
different answer as just explained above i.e. with different constitutive
models it is difficult to represent the so-called exact solutions. The
problem is considered in plane strain conditions and the fluid is injected
at the center of an initial crack with a constant injection rate of
$Q_{0}=10^{-4}m^{3}/s$. The infinite medium is approximated as in Fig.
\ref{fig4_22}a where the geometry and boundary conditions are shown, the grid size along the crack tip is $\Delta x=0.05m$. The
parameters used in the simulation are listed in Tab. \ref{Tab.3}. The fluid flow is governed by FEM only, the solid part is modeled by the coupling approach described in \citep{ni2019static,Ni2019Coupling}, the global coupling matrix is generated by using a similar way, correspondingly.} 
\begin{table}[h!]
\caption{\textcolor{blue}{Mechanical and fluid parameters used in the
reproduction of the KGD model.}}
\label{Tab.3}
\centering
\setlength{\tabcolsep}{2mm}{\  
\begin{tabular}{cccccccccc}
\toprule $E$ & $\upsilon$ & $G_c$ & $\rho$ & $\alpha$ & $n$ & $K_{w}$ & $\mu$
& $k$ &  \\ \hline
$14.4 GPa$ & $0.2$ & $10 J/m^2$ & $1000kg/m^{3}$ & $0.7883$ & $0.19$ & $3 GPa
$ & $10^{-3}Pa\cdot s$ & $2\times10^{-14}m^{2}$ &  \\ 
\bottomrule &  &  &  &  &  &  &  &  & 
\end{tabular}
}
\end{table}

\begin{figure}[h!]
\centering  
\subfloat[The result of hybrid FEM/PD
model.]{\includegraphics[scale=0.6]{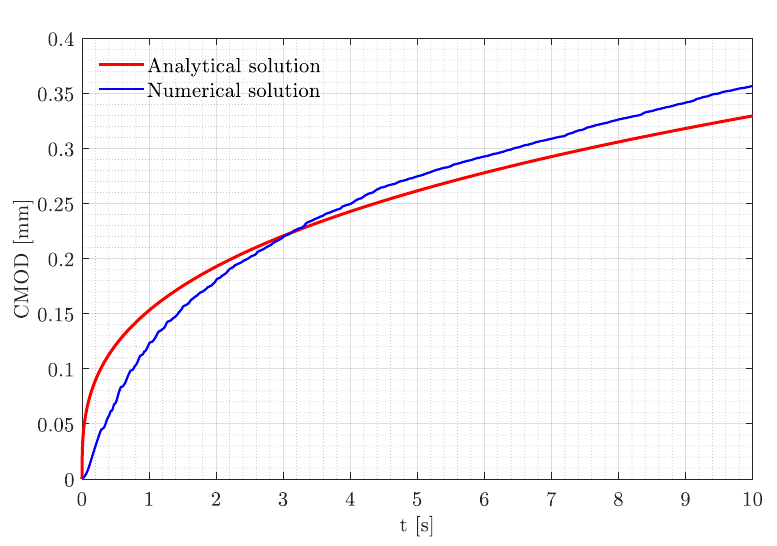}%
\label{fig4_5_1:sub1}}  \subfloat[The result reported in Secchi et al.
2007.]{%
\includegraphics[width=3.35in,height=2.1in]{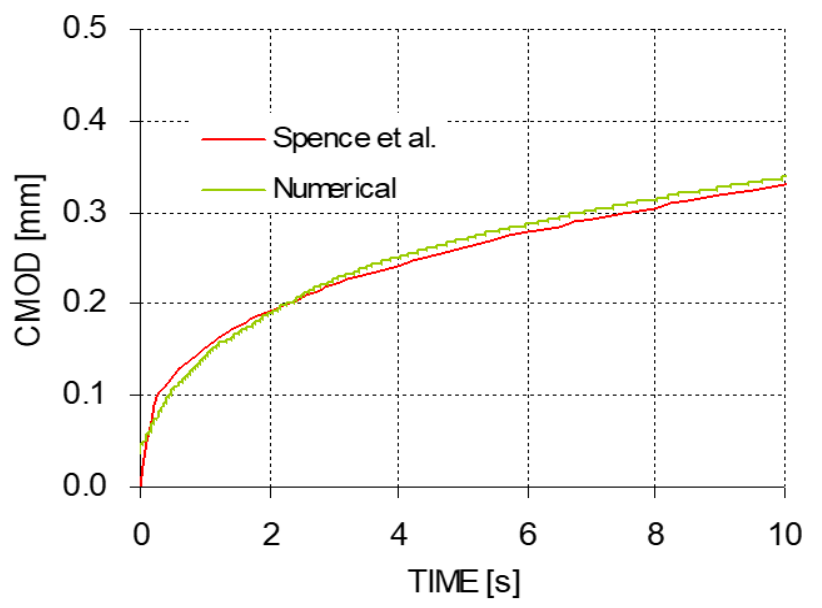}%
\label{fig4_5_1:sub2}} 
\caption{CMOD vs. time: (a) results obtained by FEM/PD model ($\Delta
t=10^{-3}s$); (b) results from \citep{secchi2007mesh} obtained by FEM
model with LEFM ($\Delta t=10^{-2}s$), redrawn with permission from Secchi
et al. 2007.}
\label{fig4_5_1}
\end{figure}

\begin{figure}[h!]
\centering  
\subfloat[The result of hybrid FEM/PD
model.]{\includegraphics[scale=0.6]{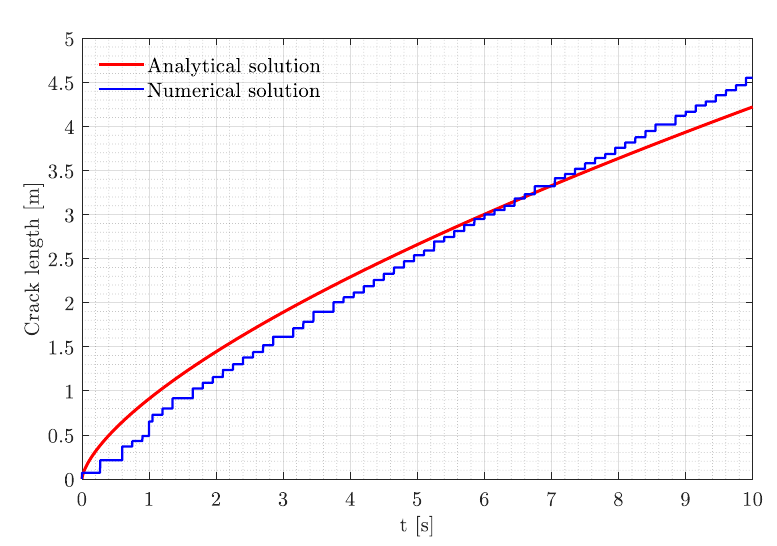}%
\label{fig4_5_2:sub1}}  \subfloat[The result reported in Secchi et al.
2007.]{%
\includegraphics[width=3.35in,height=2.1in]{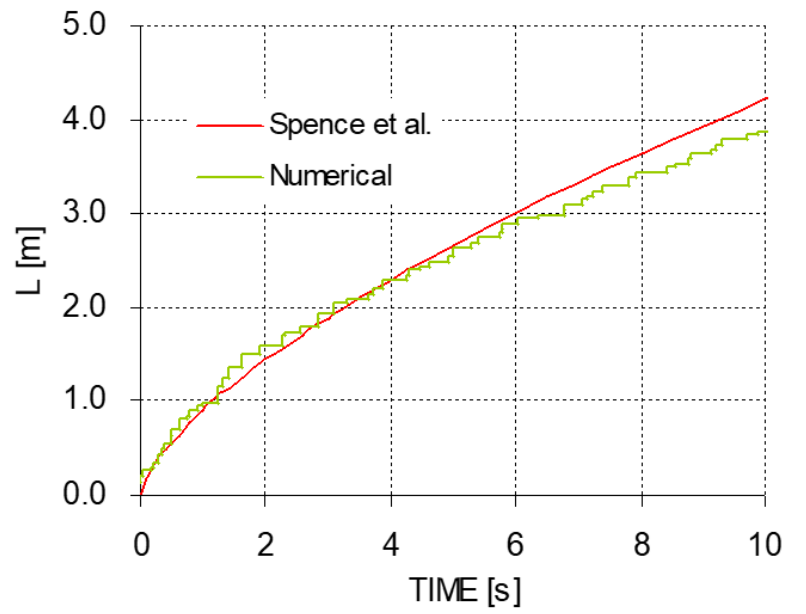}%
\label{fig4_5_2:sub2}} 
\caption{Crack length vs. time: (a) results obtained by FEM/PD model ($%
\Delta t=10^{-3}s$); (b) results from \citep{secchi2007mesh} obtained
by FEM model with LEFM ($\Delta t=10^{-2}s$), redrawn with permission from
Secchi et al. 2007. Note that the advancement steps in one sec are much less than
the number of time steps. This is a reliable indication that the behaviour
is truly stepwise.}
\label{fig4_5_2}
\end{figure}

\begin{figure}[h!]
\centering  
\subfloat[The result of hybrid FEM/PD
model.]{\includegraphics[scale=0.6]{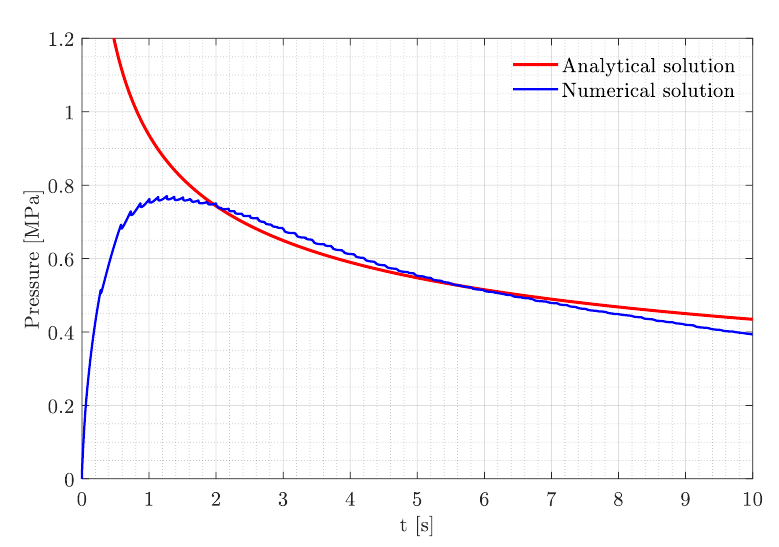}%
\label{fig4_5_3:sub1}}  \subfloat[The result reported in Peruzzo et al. 2019.]{%
\includegraphics[width=3.35in,height=2.1in]{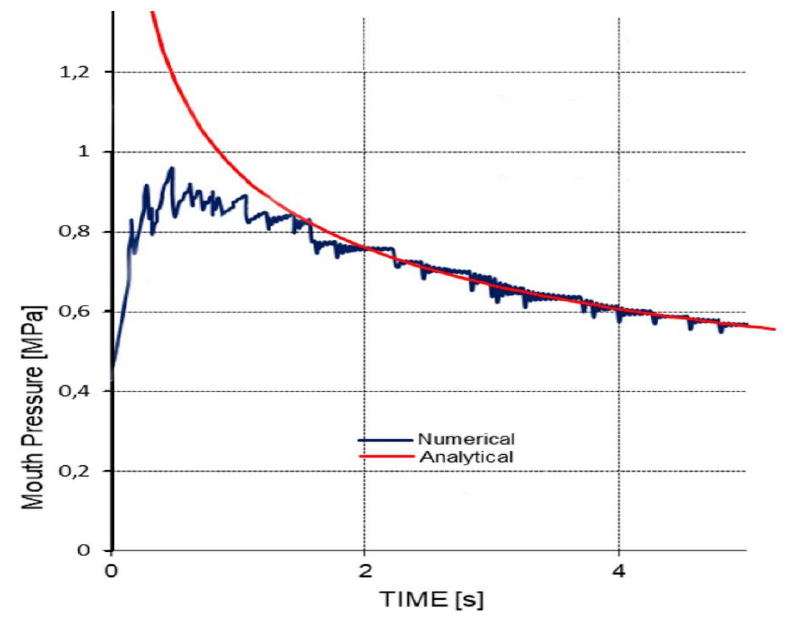}%
\label{fig4_5_3:sub2}} 
\caption{Pressure at the injection point vs. time: (a) results obtained by
FEM/PD model; (b) results from \citep{peruzzo2019dynamics} obtained by FEM
model with LEFM , redrawn with permission from Peruzzo et al. 2019.}
\label{fig4_5_3}
\end{figure}

\textcolor{blue}{The obtained results with the hybrid FEM/PD model are
compared with the exact solution and with the results of the FEM solution of
Secchi et al. \citep{secchi2007mesh} and Peruzzo et al.
\citep{peruzzo2019dynamics} in Figs. \ref{fig4_5_1} to \ref{fig4_5_3} for
respectively crack mouth opening displacement, crack length and pressure at
the injection point.}


\textcolor{blue}{Given the difference in the adopted constitutive models the
results are acceptable. In particular the stepwise crack tip advancement
with irregular steps is recovered. This shows that the hybrid FEM/PD method
is able to represent experimentally observed features, which are rather
difficult to reproduce. Fewer crack advancement steps in a time interval
than the number of time steps as in Fig. \ref{fig4_5_2} is a reliable
indicator that the behavior is truly stepwise. The stepwise advancement has
be explained by Biot's theory in \citep{milanese2016avalanches}. A further
explanation can be borrowed from Huang and Detournay
\citep{huang2013discrete} who obtained in a fracturing DEM model \textquotedblleft plateaus\textquotedblright \
as in Fig. \ref{fig4_5_2} and observe that  \textquotedblleft when the crack is temporarily
arrested, another sudden extension in the crack length occurs only after
enough energy has been accumulated in the system\textquotedblright. The example also
evidences the impact of the assumed solid behaviour in HF modelling.}
\newpage
\subsection{Fluid-driven fracture propagation: interaction between HF crack
and natural cracks}

\begin{figure}[tbh]
\centering  
\subfloat[Case
1.]{\includegraphics[scale=0.4]{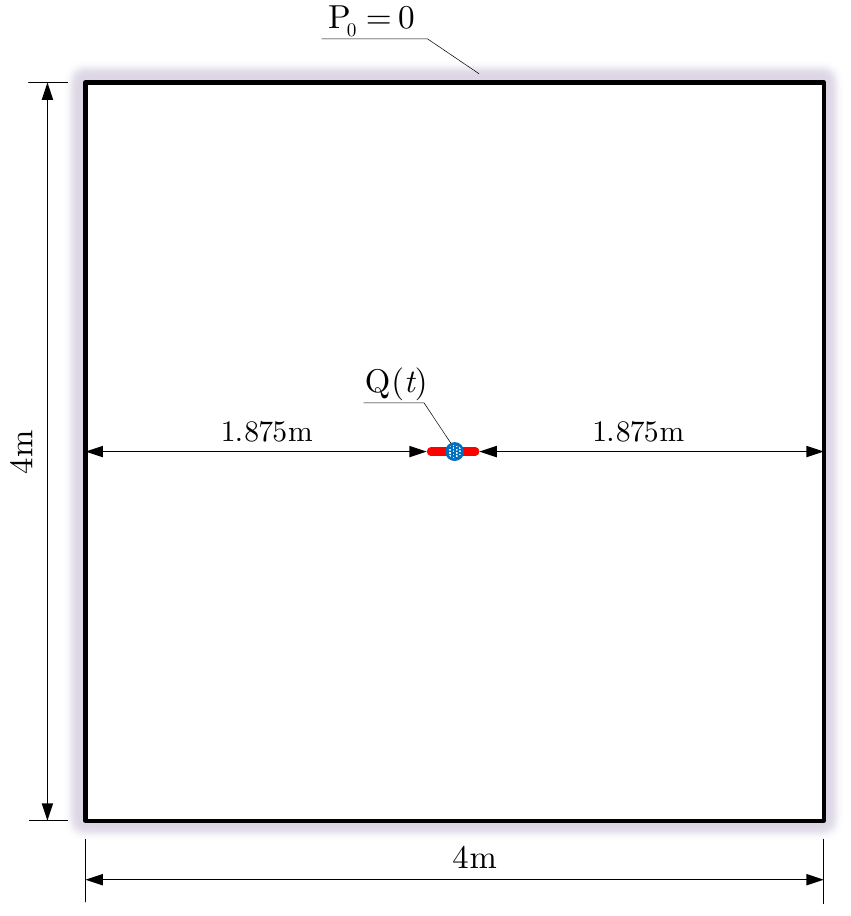}\label{fig4_4_2:sub1}}
\subfloat[Case
2.]{\includegraphics[scale=0.4]{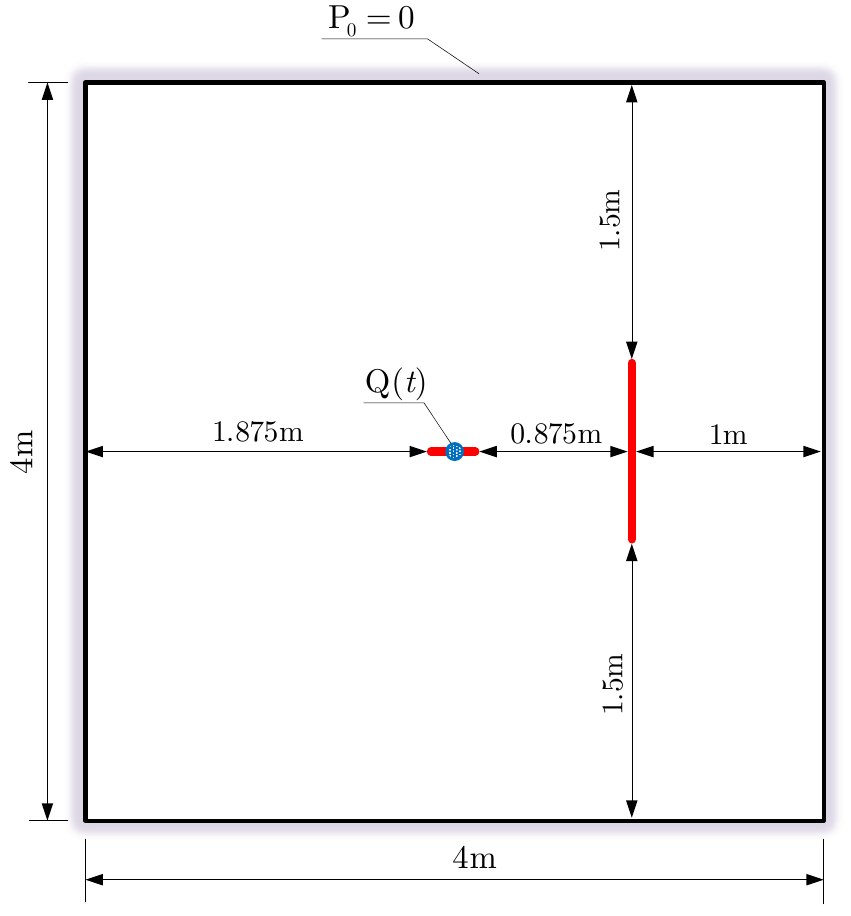}\label{fig4_4_2:sub2}}%
\newline
\subfloat[Case
3.]{\includegraphics[scale=0.4]{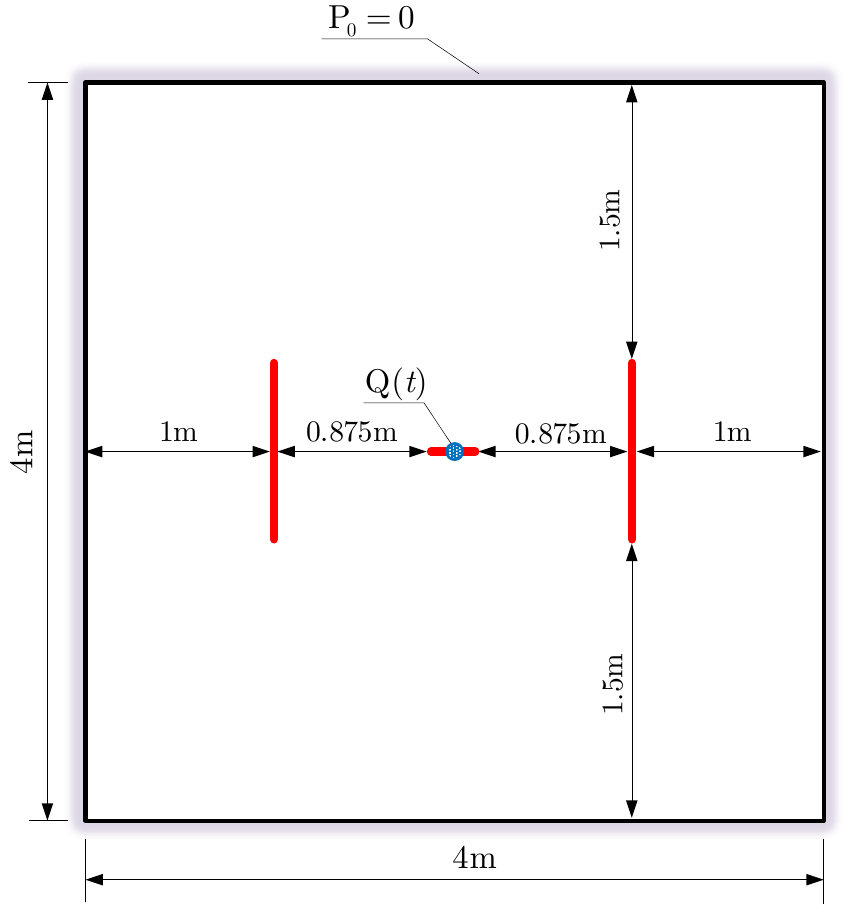}\label{fig4_4_2:sub3}}
\caption{Geometry and boundary conditions of the notched specimens subjected
to internal fluid injection.}
\label{fig4_4_2}
\end{figure}

In this section, three different fluid-driven fracture cases are studied to
evidence the characteristics and capabilities of the present method.
Geometric parameters and boundary conditions are shown in Fig. \ref{fig4_4_2}%
. Referring to \citep{lee2016pressure}, the mechanical and fluid parameters
used in the calculations are shown in Tab. \ref{Tab.2}. In case 1, a crack
is located at the center of the specimen with the initial length $l=0.25m$,
while in cases $2$ and $3$, additional vertical cracks with a length of $1m$
are set in the specimen to study the interaction between hydraulic cracks
and pre-existing natural cracks. The fluid is injected at the center of the
initial crack with a constant volume rate of $Q=1\times 10^{-3}m^{3}/s$.
This and the following flow rates, together with the assumed dynamic
viscosity, are in the range of high flow rate and low viscosity,
investigated in \citep{lhomme2002experimental}. The whole model is
discretized using uniform quadrilateral elements for fluid flow with a grid
size of $\Delta x=1\times 10^{-2}m$, while the PD grid shares the same node
positions. The corresponding horizon is $\delta =m_{ratio}\times \Delta
x=3\times 10^{-2}m$. The time integral parameters are taken as: $\Delta
t=1\times 10^{-3}s$ and $\vartheta =1$. 

\begin{table}[tbh!]
\caption{Mechanical and fluid parameters used in fluid-driven fracture
examples.}
\label{Tab.2}\centering
\setlength{\tabcolsep}{2mm}{\ 
\begin{tabular}{cccccccccc}
\toprule $E$ & $\upsilon$ & $G_c$ & $\rho$ & $\alpha$ & $n$ & $K_{w}$ & $\mu$
& $k$ &  \\ \hline
$10^8 Pa$ & $0.2$ & $100 J/m^2$ & $1000kg/m^{3}$ & $1$ & $0.4$ & $10^{8}Pa$
& $10^{-3}Pa\cdot s$ & $10^{-12}m^{2}$ &  \\ 
\bottomrule &  &  &  &  &  &  &  &  & 
\end{tabular}
}
\end{table}

\begin{figure}[tbh]
\centering  
\includegraphics[scale=0.8]{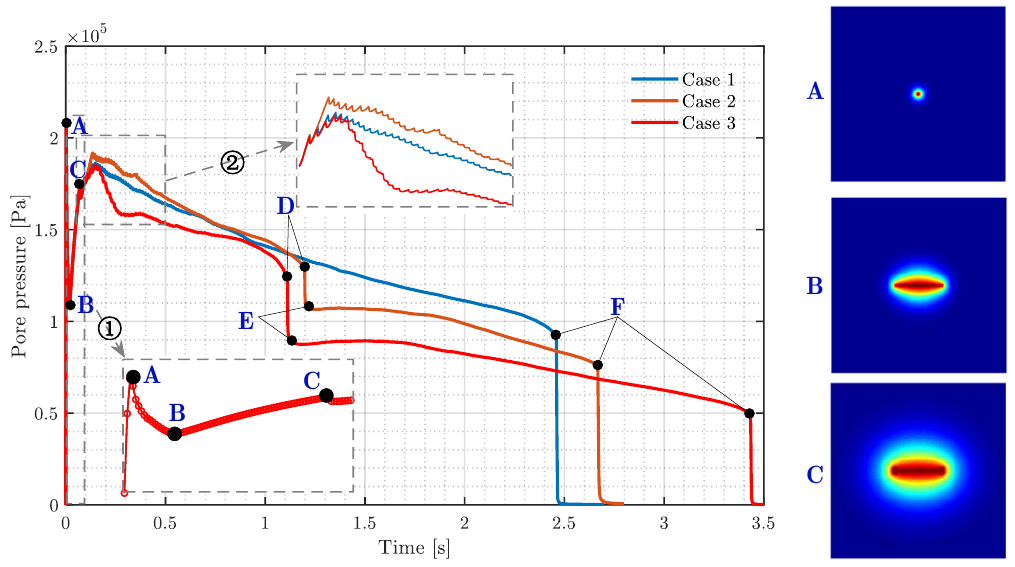}
\caption{Variation of pressure value at injection point with time.}
\label{fig4_12}
\end{figure}

The values of the pressure at the injection point are recorded in time and
plotted in Fig. \ref{fig4_12}. As shown in the magnifying frame \ding{172}
of Fig. \ref{fig4_12}, at the moment of injection, the pressure value at the
injection point suddenly increases to a large value, under its action, the
initial central crack opens gradually, and the permeability in the initial
central crack domain increases correspondingly, which leads to a steep drop
of the pressure. Shown at the right half part of Fig. \ref{fig4_12}, the
fluid diffusion is isotropic because of the small crack opening width (see
the pressure distribution at instant A). With the increase of the crack
opening width, the fluid diffusion takes precedence over crack direction
(instant B), until the pressure values at the crack tips are approximately
equal to that at the injection point (instant C). After instant C, hydraulic
fracture occurs in all the cases, and the pressure values can rise a bit at
the incipient stage because the crack propagates slowly, soon afterwards,
hydraulic cracks propagate rapidly, the pressure values drop gradually. In
cases 2 and 3, after the hydraulic cracks touch the natural vertical crack
sets and force them to open, the pressure values suddenly drop at extreme
speed (between instants D and E), afterwards, the pressure values in the
vertical cracks increase gradually until the new cracks start to propagate
from their tips and finally extend to the boundary of the specimen. After
cracks touch the 0-pressure boundary of the specimen, the pressure value at
the injection point drops again at an extreme rate (see the parts after the
instant F in Fig. \ref{fig4_12}).

\begin{figure}[tbh!]
\centering  
\subfloat[Damage level at
$0.1s$.]{\includegraphics[scale=0.4]{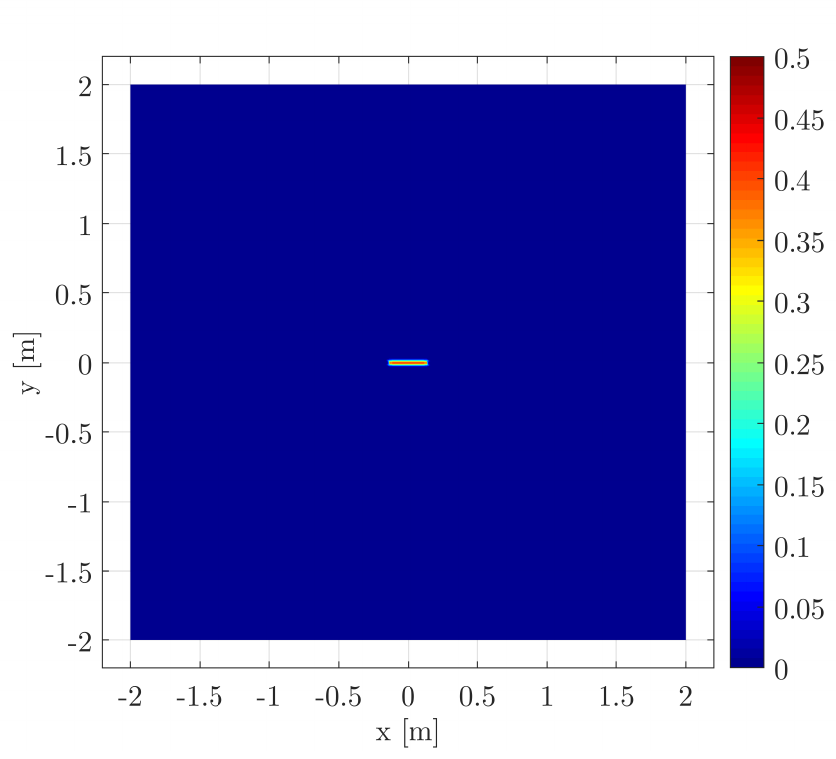}\label{fig4_6:sub1}}
\subfloat[Damage level at
$1s$.]{\includegraphics[scale=0.4]{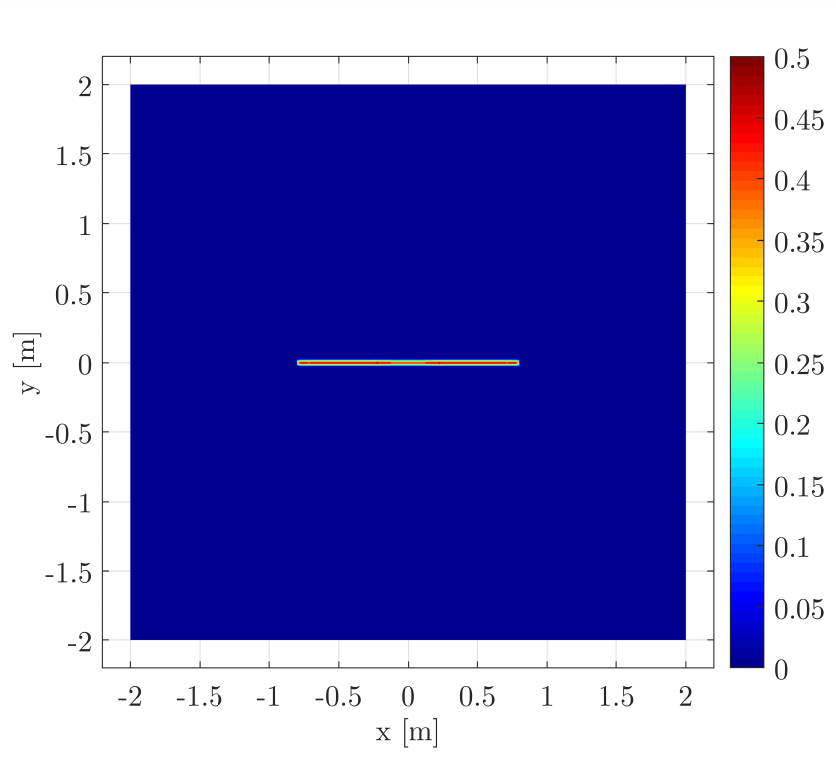}\label{fig4_6:sub3}}
\subfloat[Damage level at
$2.4s$.]{\includegraphics[scale=0.4]{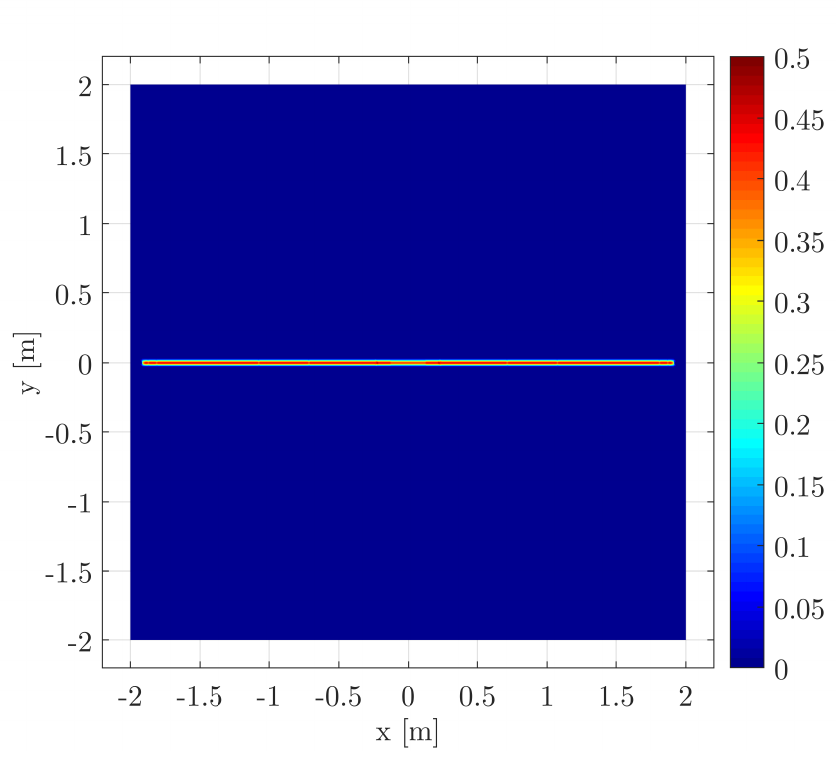}\label{fig4_6:sub5}}
\newline
\subfloat[Pressure distribution at
$0.1s$.]{\includegraphics[scale=0.4]{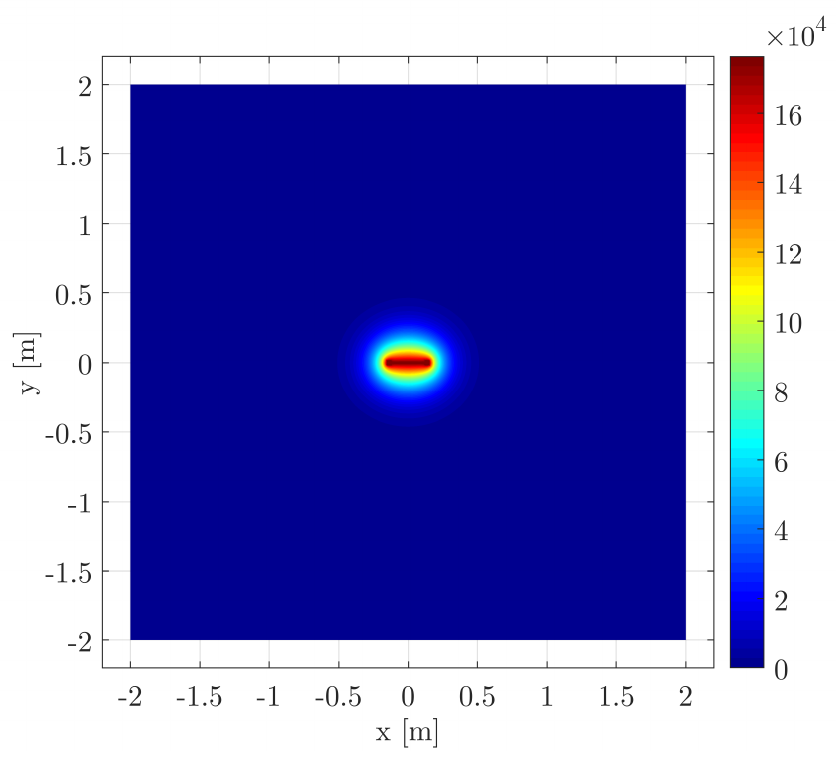}\label{fig4_6:sub2}}
\subfloat[Pressure distribution at
$1s$.]{\includegraphics[scale=0.4]{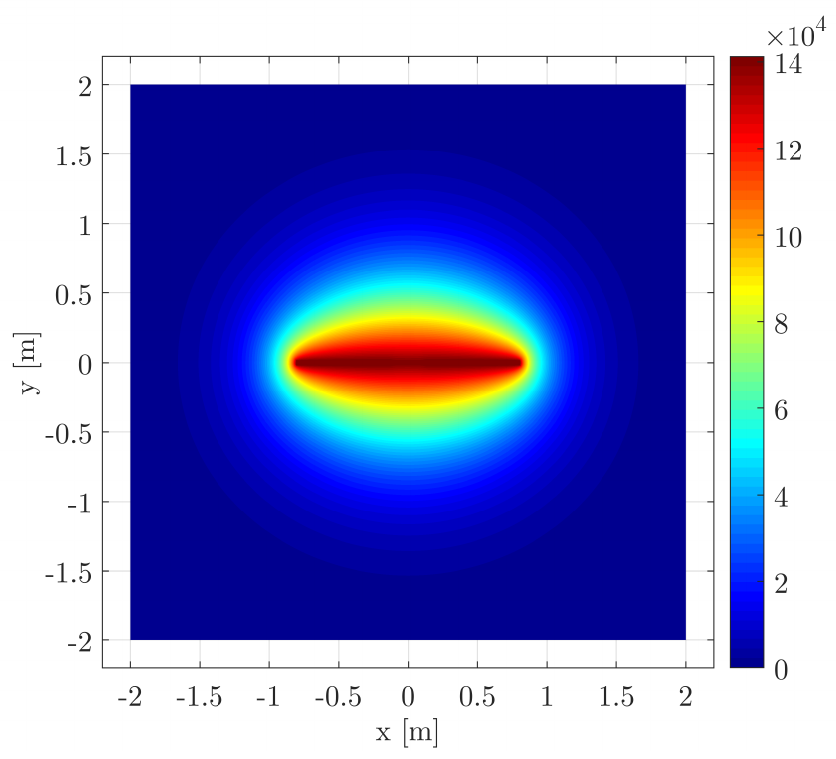}\label{fig4_6:sub4}}
\subfloat[Pressure distribution at
$2.4s$.]{\includegraphics[scale=0.4]{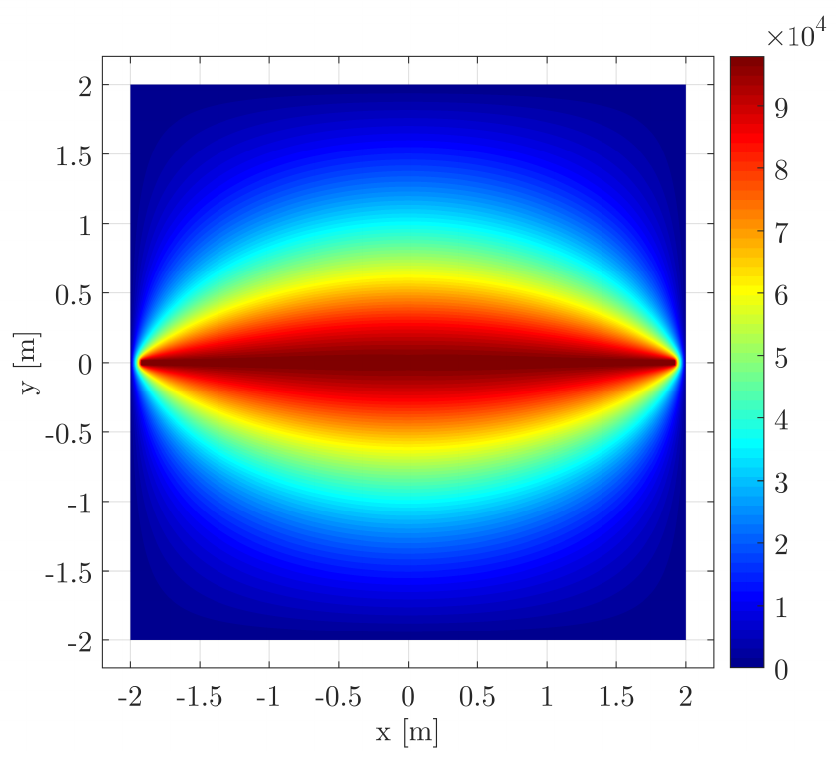}\label{fig4_6:sub6}}
\caption{Crack propagation and variation of pressure distribution in case 1
at three time instants.}
\label{fig4_6}
\end{figure}

\begin{figure}[tbh!]
\centering  
\subfloat[Damage level at
$0.1s$.]{\includegraphics[scale=0.4]{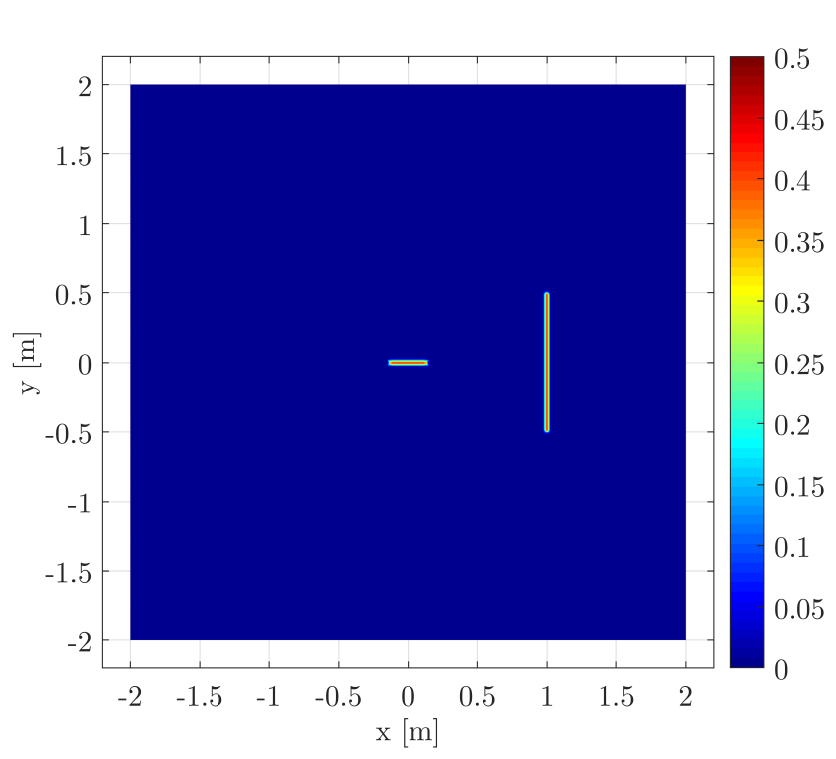}\label{fig4_7:sub1}}
\subfloat[Damage level at
$1.2s$.]{\includegraphics[scale=0.4]{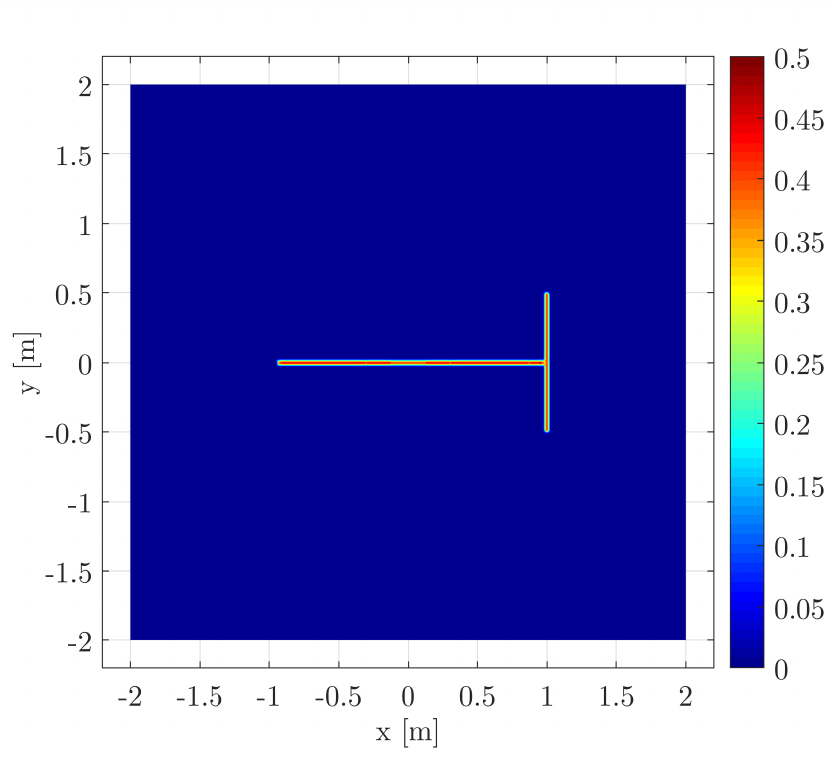}\label{fig4_7:sub3}}
\subfloat[Damage level at
$2.4s$.]{\includegraphics[scale=0.4]{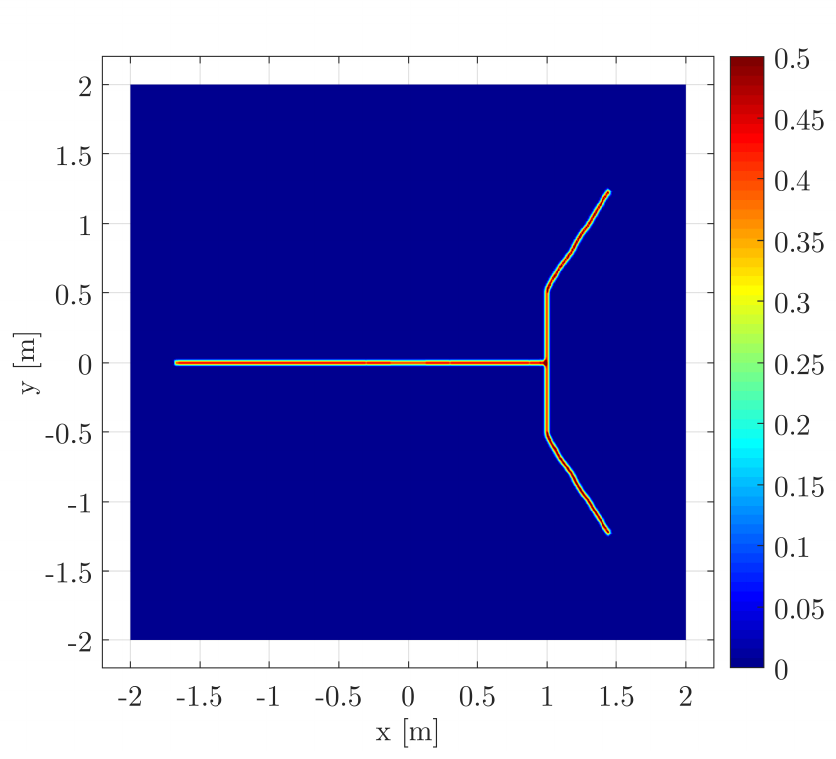}\label{fig4_7:sub5}}
\newline
\subfloat[Pressure distribution at
$0.1s$.]{\includegraphics[scale=0.4]{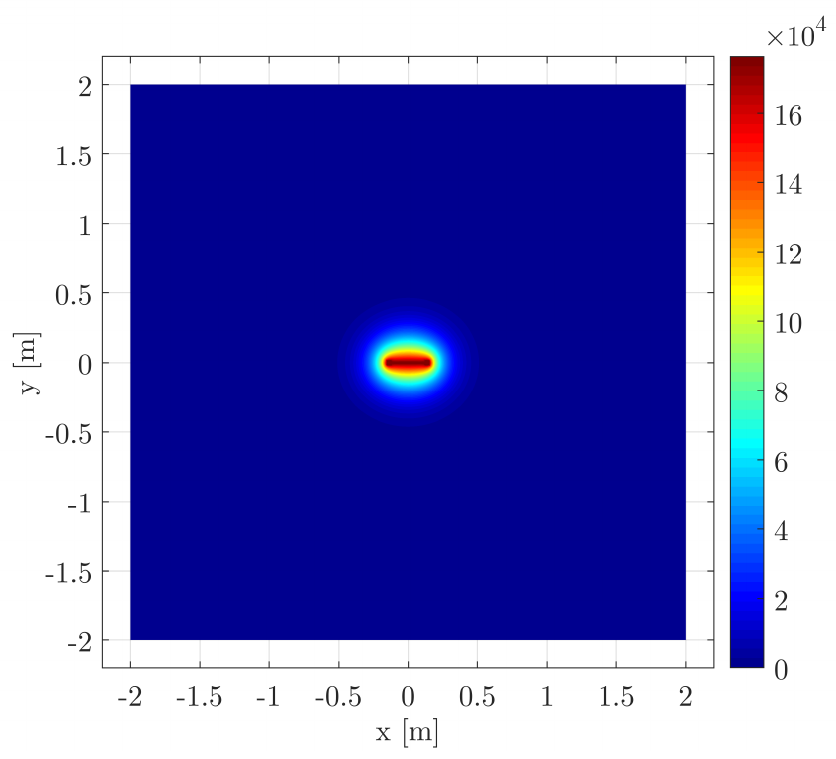}\label{fig4_7:sub2}}
\subfloat[Pressure distribution at
$1.2s$.]{\includegraphics[scale=0.4]{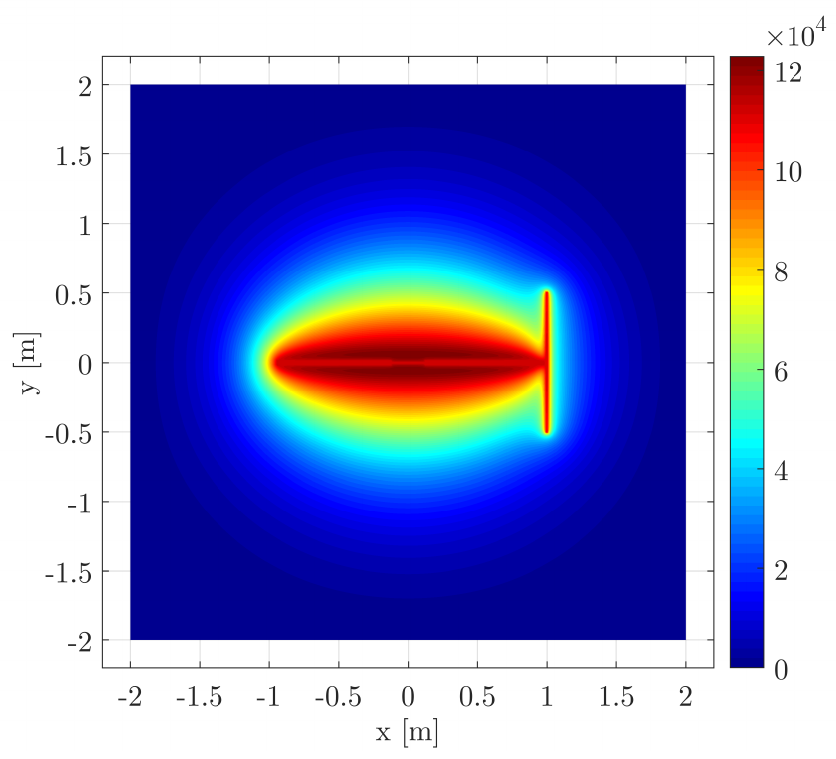}\label{fig4_7:sub4}}
\subfloat[Pressure distribution at
$2.4s$.]{\includegraphics[scale=0.4]{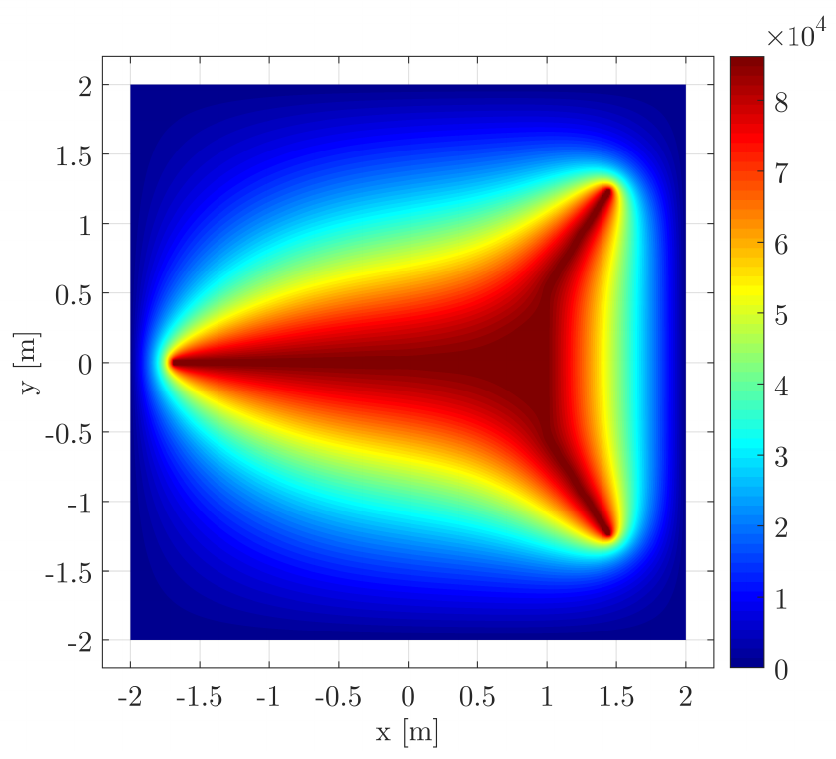}\label{fig4_7:sub6}}
\caption{Crack propagation and variation of pressure distribution in case 2
at three time instants.}
\label{fig4_7}
\end{figure}

\begin{figure}[tbh]
\centering  
\subfloat[Damage level at
$0.1s$.]{\includegraphics[scale=0.4]{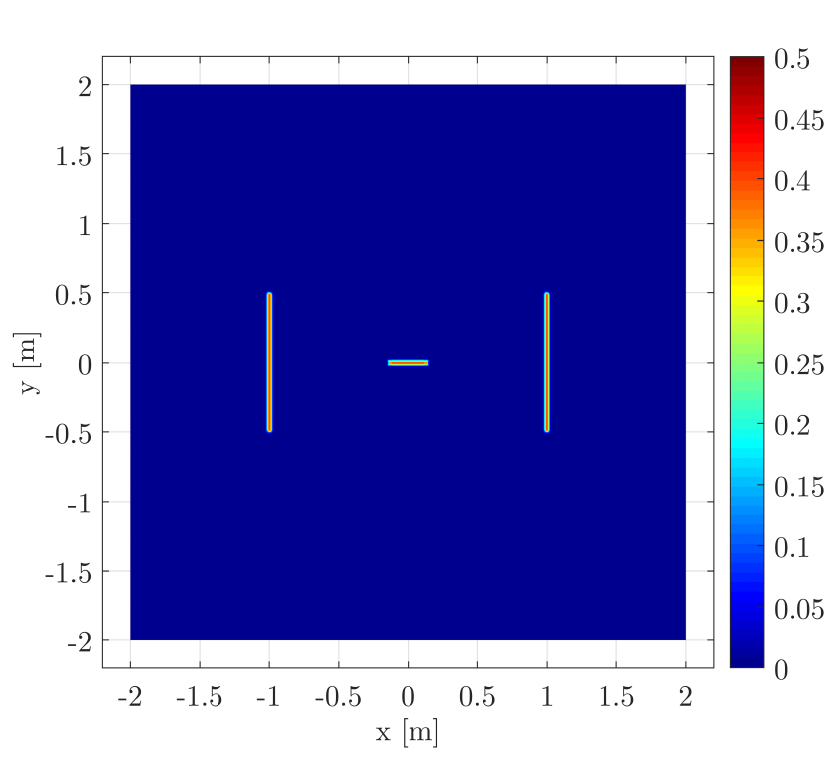}\label{fig4_8:sub1}}
\subfloat[Damage level at
$1.2s$.]{\includegraphics[scale=0.4]{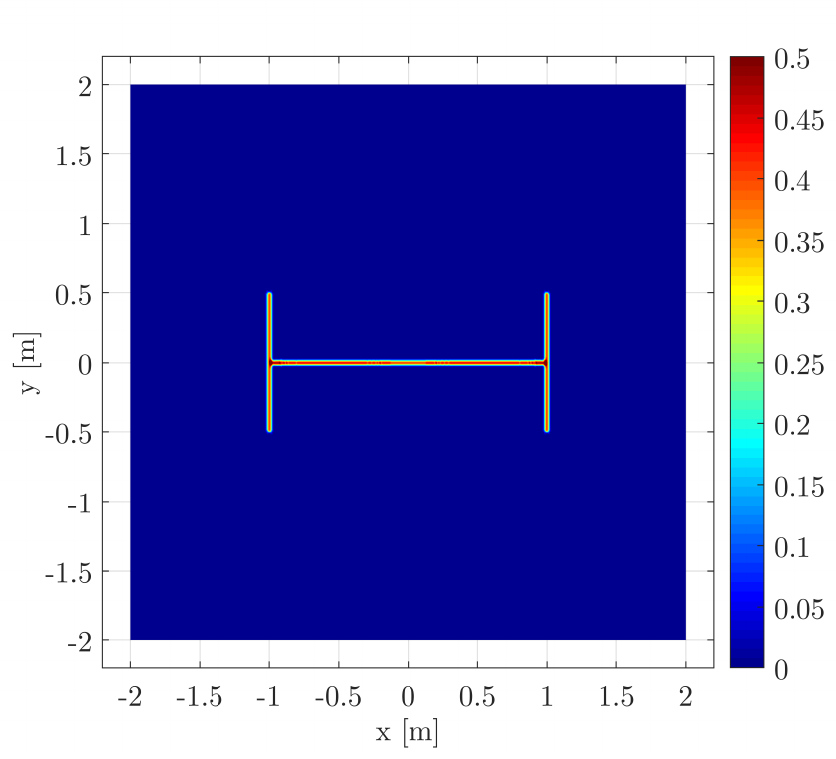}\label{fig4_8:sub2}}
\subfloat[Damage level at
$2.4s$.]{\includegraphics[scale=0.4]{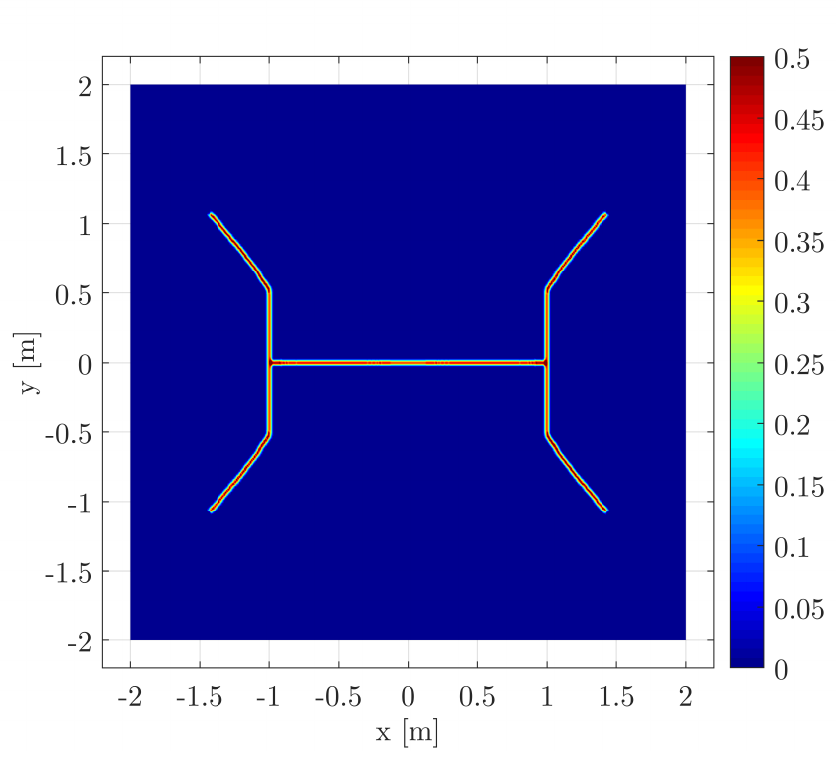}\label{fig4_8:sub3}}
\newline
\subfloat[Pressure distribution at
$0.1s$.]{\includegraphics[scale=0.4]{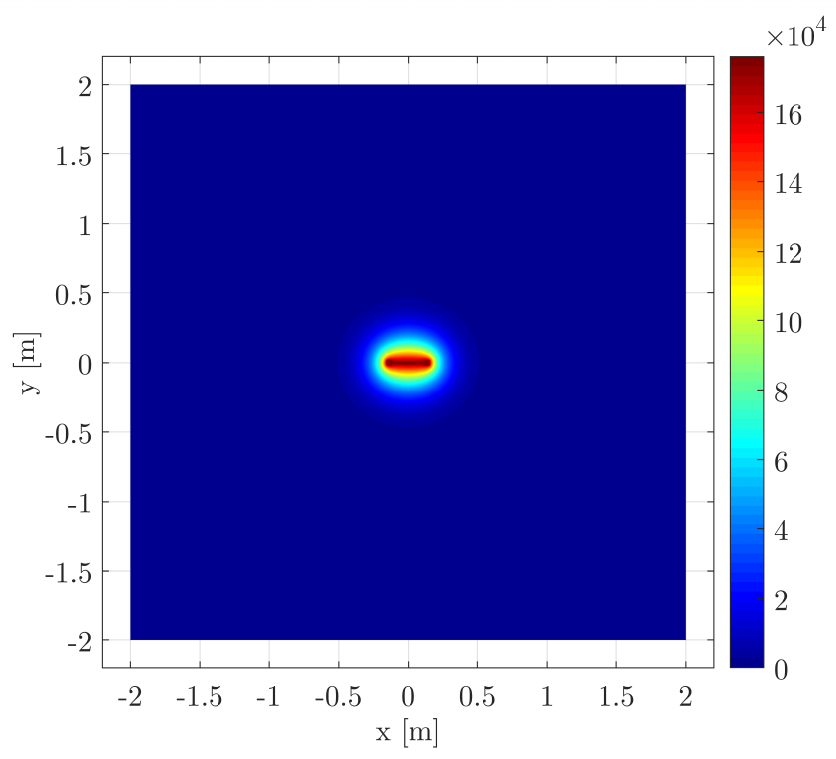}\label{fig4_8:sub4}}
\subfloat[Pressure distribution at
$1.2s$.]{\includegraphics[scale=0.4]{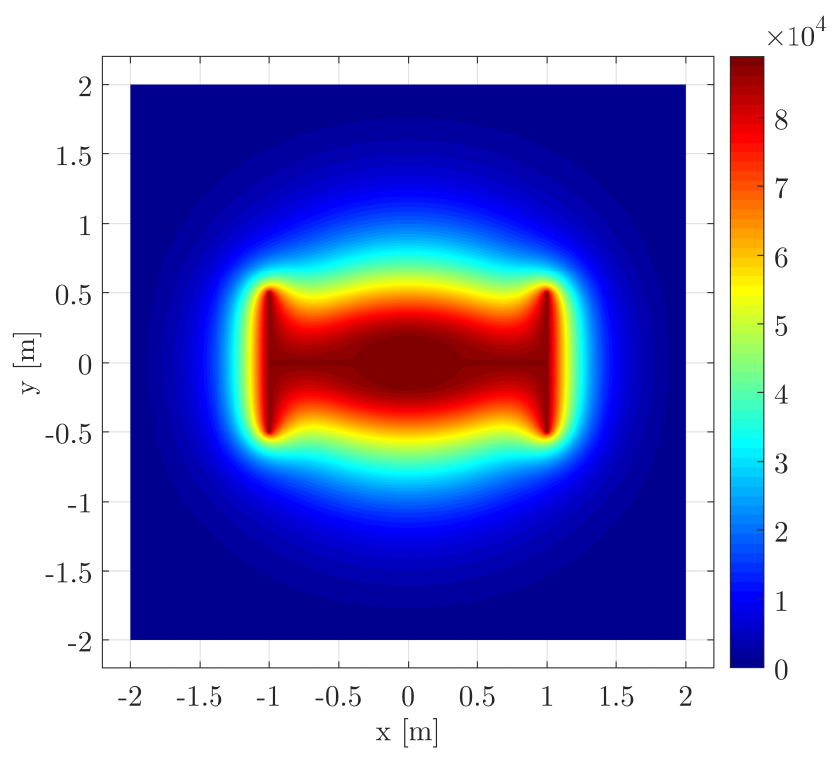}\label{fig4_8:sub5}}
\subfloat[Pressure distribution at
$2.4s$.]{\includegraphics[scale=0.4]{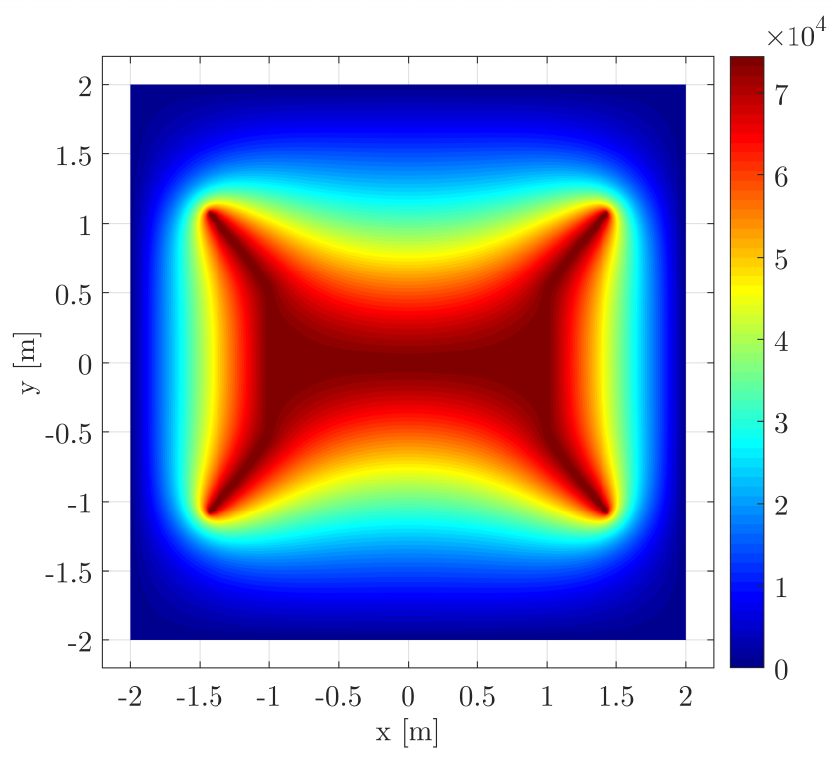}\label{fig4_8:sub6}}
\caption{Crack propagation and variation of pressure distribution in case 3
at three time instants.}
\label{fig4_8}
\end{figure}

Figs. \ref{fig4_6}, \ref{fig4_7} and \ref{fig4_8} show the crack patterns
and the corresponding variation of pressure distributions in cases 1, 2 and
3 at three different time instants during the fracturing process. With the
injection of fluid, the injected fluid flows along the cracks and forces the
fluid pressure near the crack area to increase, correspondingly, the solid
skeleton deforms and crack propagates. The present method can describe the
phenomena of spontaneous crack propagation and the fluid flow in fractured
saturated porous media in a simple but reasonable way.

\begin{figure}[tbh]
\centering  
\subfloat[Crack pattern in case
1-1.]{\includegraphics[scale=0.4]{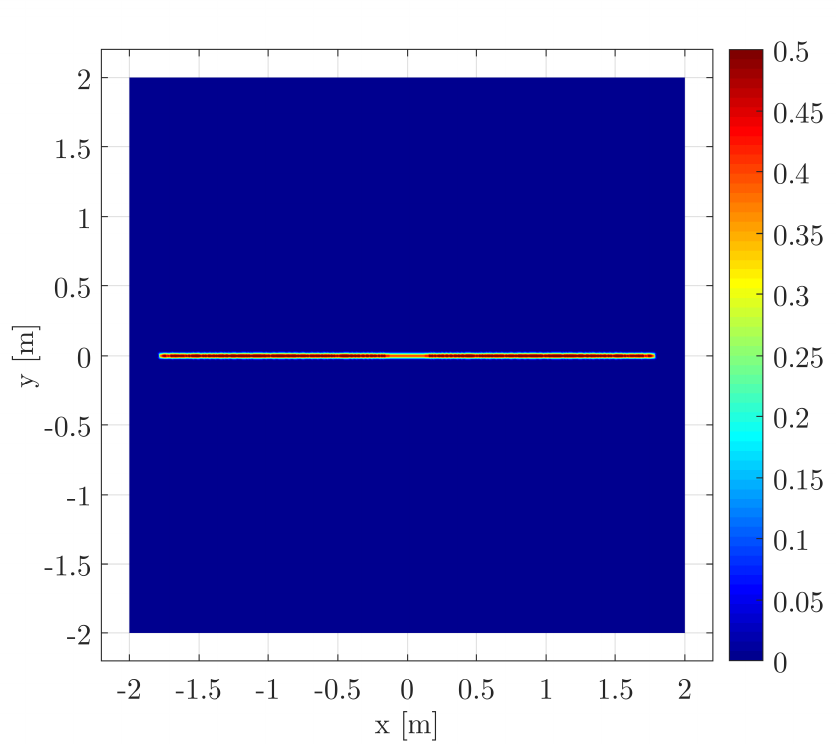}\label{fig4_9:sub1}}
\subfloat[Crack pattern in case
1-2.]{\includegraphics[scale=0.4]{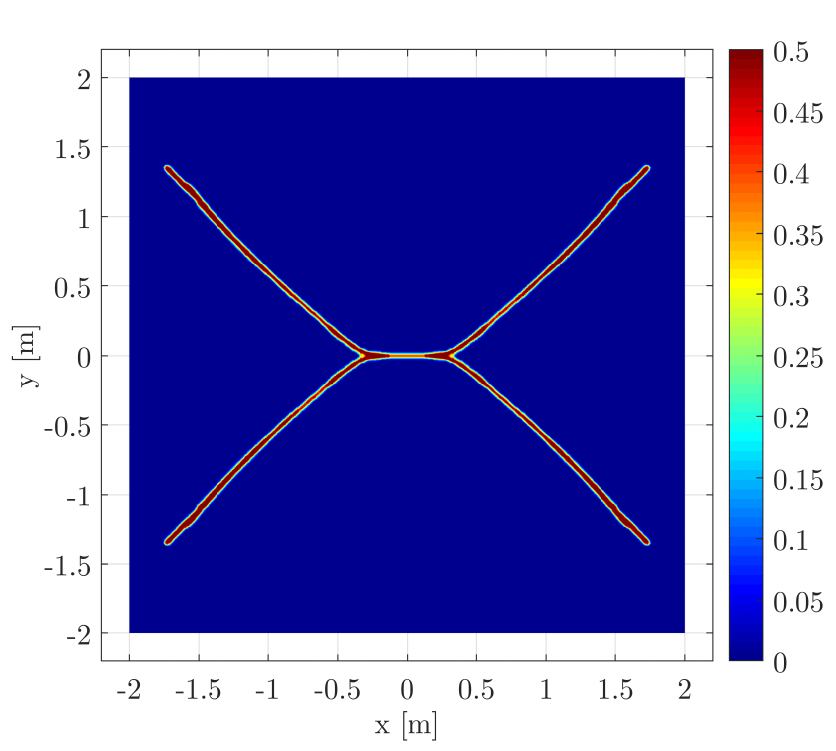}\label{fig4_9:sub2}}
\subfloat[Crack pattern in case
1-3.]{\includegraphics[scale=0.4]{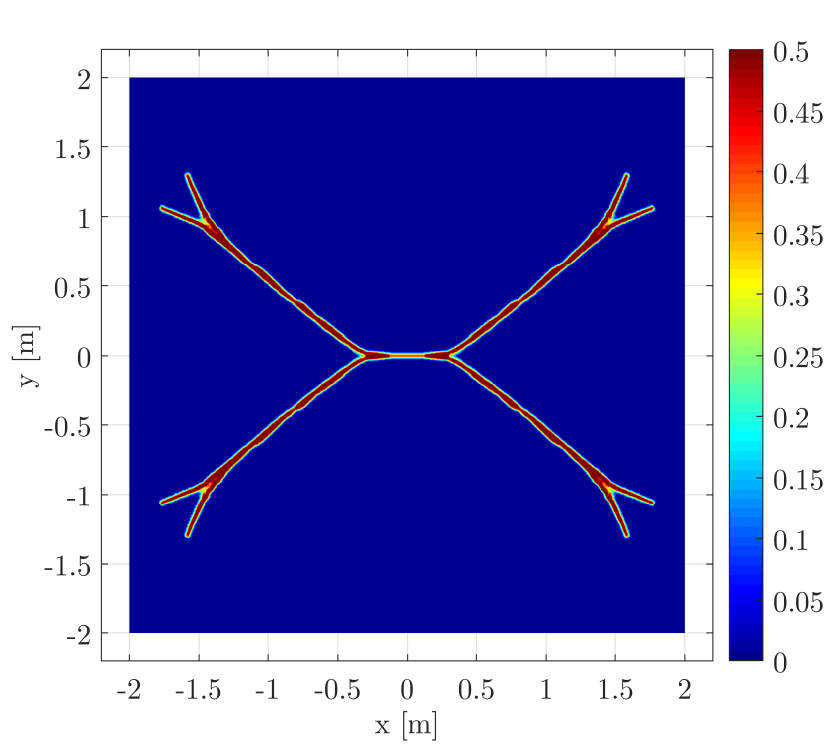}\label{fig4_9:sub3}}
\newline
\subfloat[Pressure distribution in case
1-1.]{\includegraphics[scale=0.4]{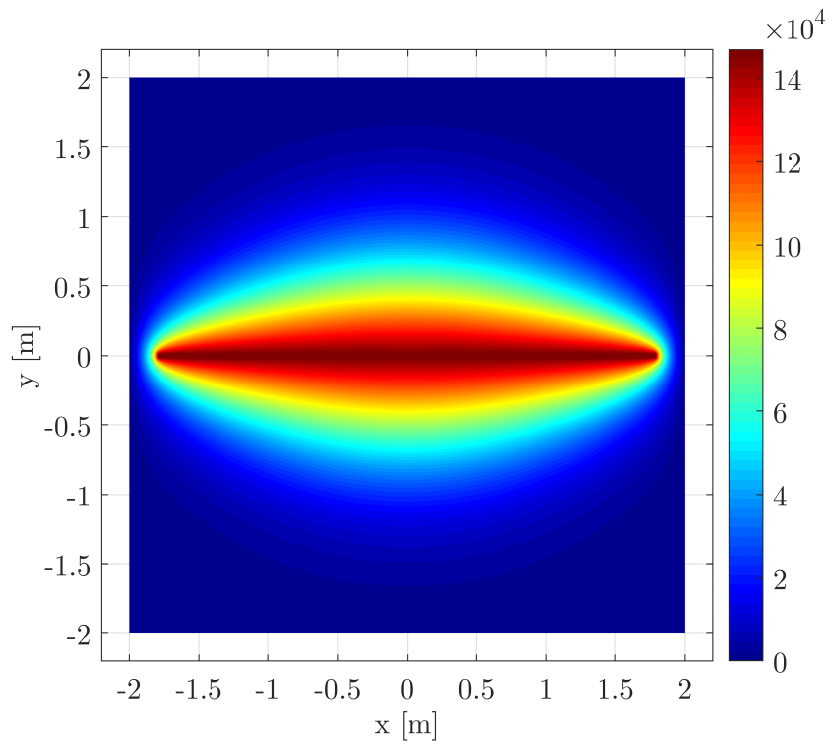}\label{fig4_9:sub4}}
\subfloat[Pressure distribution in case
1-2.]{\includegraphics[scale=0.4]{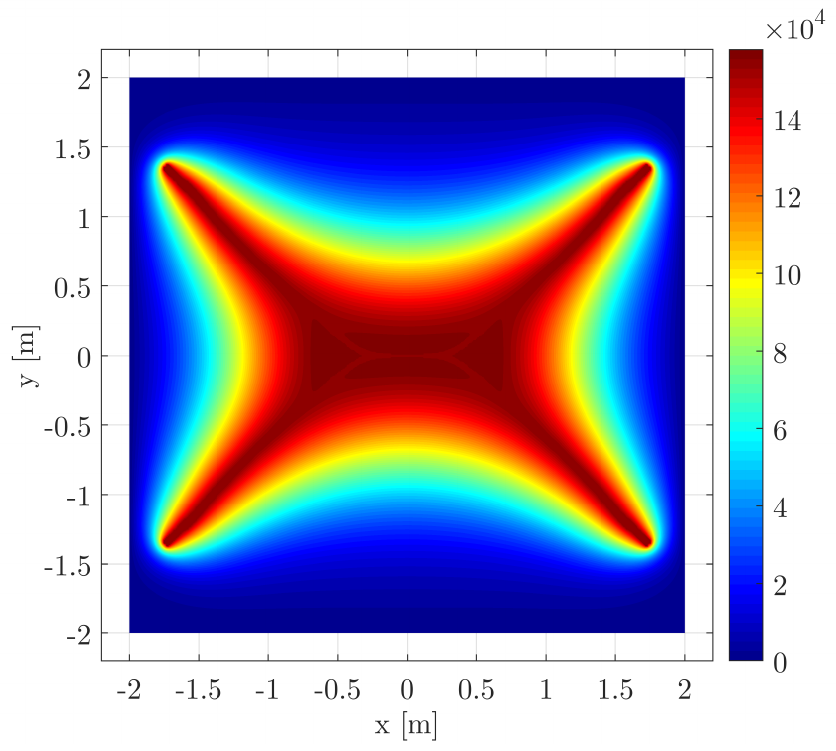}\label{fig4_9:sub5}}
\subfloat[Pressure distribution in case
1-3.]{\includegraphics[scale=0.4]{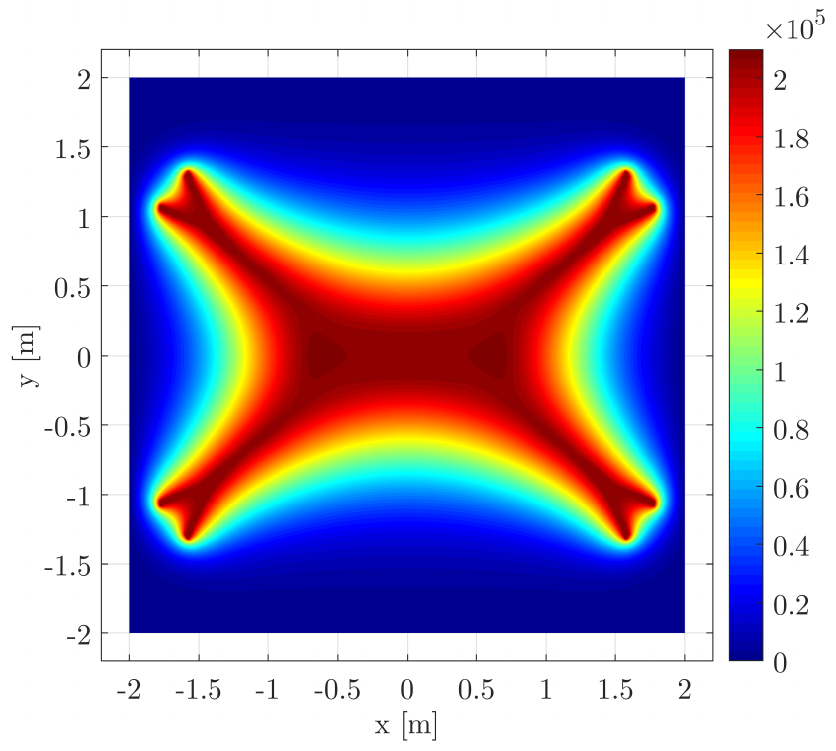}\label{fig4_9:sub6}}
\caption{Crack patterns and pressure distributions at the final step of case
1 with different injection rates.}
\label{fig4_9}
\end{figure}

\begin{figure}[tbh!]
\centering  
\includegraphics[scale=0.8]{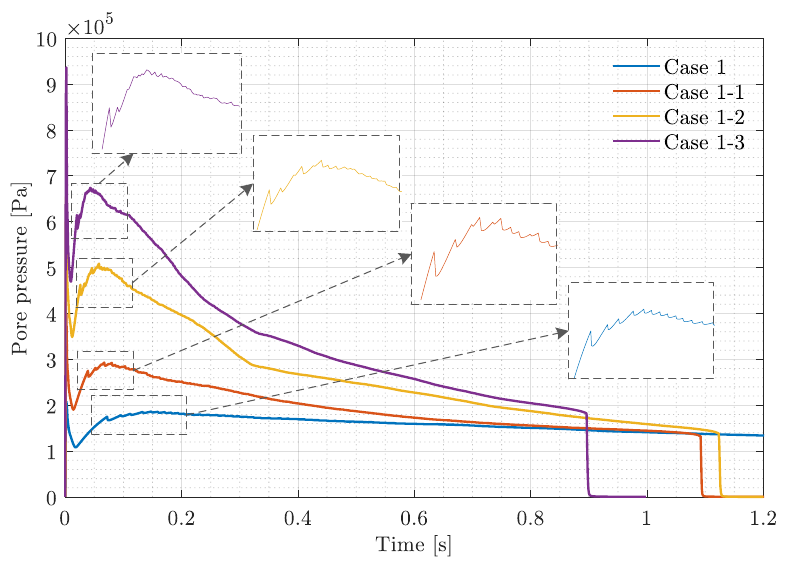}
\caption{The variation of pressure value at injection point with time of
cases 1, 1-1, 1-2 and 1-3.}
\label{fig4_10}
\end{figure}

\begin{figure}[tbh!]
\centering  
\includegraphics[scale=0.8]{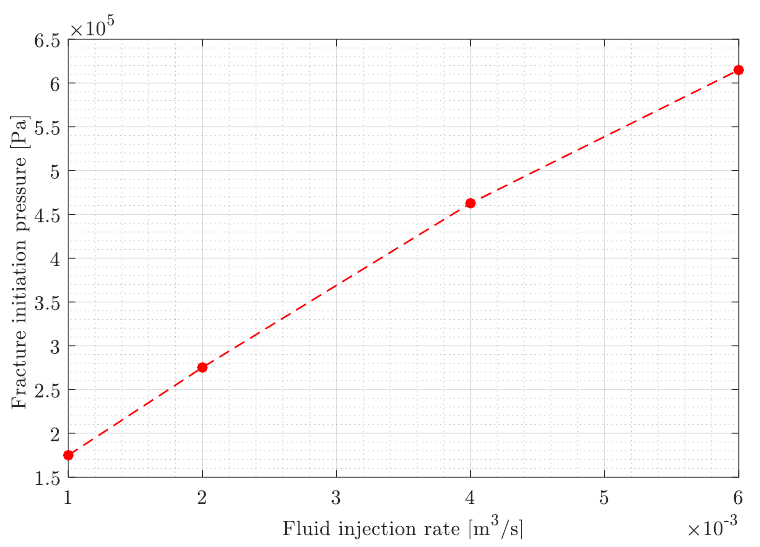}
\caption{The relationship between the injection rates and the fracture
initiation pressures.}
\label{fig4_11}
\end{figure}

\clearpage To investigate the effects of the injection rate on hydraulic
fracturing, the geometry in Fig. \ref{fig4_4_2:sub1} is adopted with three
other different injection rates, which are marked as case 1-1 ( $Q=2\times
10^{-3}m^{3}/s$), case 1-2 ( $Q=4\times 10^{-3}m^{3}/s$) and case 1-3 ( $%
Q=6\times 10^{-3}m^{3}/s$). Figs. \ref{fig4_9:sub1} to \ref{fig4_9:sub3}
show the final crack patterns while Figs. \ref{fig4_9:sub4} to \ref%
{fig4_9:sub6} show the corresponding distributions of the pressure. With
certain high injection rates, hydraulic crack bifurcation can be observed,
which is a typical phenomenon in the dynamic hydraulic fracture. The
variation of pressure value at the injection point with time in cases 1,
1-1, 1-2 and 1-3 are plotted in Fig. \ref{fig4_10}, from which it can be
observed that the fracture initiation pressure at the injection point
increases with the increase of injection rate and their numerical
relationship is plotted in Fig. \ref{fig4_11}. In addition, as shown in the
magnifying frame \ding{173} of Fig. \ref{fig4_12} and the magnifying frames
in Fig. \ref{fig4_10}, the fluid pressure at the injection point presents a
characteristic oscillation, which is consistent with the typical pattern
observed experimentally in \citep{lhomme2002experimental} and the numerical
evidence in \citep{cao2018porous}. In the fracture event the pressure falls
rapidly, indicating that the volume of the induced cracks increases faster
than the injection rate \citep{okland2002importance}. 

\newpage

\section{Conclusions}

This paper presents a hybrid FEM and PD modelling approach for simulating
hydraulic fracture propagation in saturated porous media. In the proposed
approach, FE equations are used to govern the fluid flow, while PD is used
to describe the deformation of the solid phase and capture the crack
propagation.

Several numerical examples are carried out. First, the proposed approach is
validated by two benchmark examples, and the numerical results agree well
with the analytical solutions. Then, several numerical examples are
presented to further demonstrate the capabilities and the main features of
the proposed approach in simulating crack propagation and bifurcation in
saturated porous media under pressure- and fluid-driven conditions.

The phenomenon of fluid pressure oscillations is observed in the
fluid-driven hydraulic fracture examples, in agreement with experimental
observations.

\textcolor{blue}{The presented method is one of the first to solve multi-physics problems involving
discontinuities by using local models coupled with PD.} Similar approaches
could be used in the future to effectively simulate crack propagation due to
other multi-physics fields: electricity, thermal, chemical, ... the FEM
method can be used to solve the multi-physics equations and the coupled
Peridynamics model can solve the structural problem in the solid. 
\textcolor{blue}{
\section*{List of Symbols}
\begin{tabular}[t]{|lp{14cm}|}
	\hline
	$\boldsymbol{x},\,\boldsymbol{x}' $ & location vector of material points \\ 
	$\boldsymbol{u} $                 & {displacement vector} \\ 	
	$\boldsymbol{\xi} $          & {relative position vector of two material points} \\	
	$\boldsymbol{\eta} $          & {relative displacement vector of two material points} \\	
	$\underline{\boldsymbol{X}}$ & {reference position vector state} \\
	$\underline{\boldsymbol{Y}}$ & {deformation vector state} \\
	$\underline{x}$          & {reference position scalar state} \\
	$\underline{y}$          & {deformation scalar state} \\
	$\underline{\boldsymbol{M}}$  & {unit state in the direction of deformed bond} \\ 
	$\underline{e} $          & {extension scalar state} \\  
	$\textbf{\b{T}}$ & {force density vector states}\\
	$\boldsymbol{\ddot{u}}$ & {acceleration vector of material point}\\
	$\boldsymbol{b} $          & {body force density} \\	
	$\mathcal{H}_{x}$ & {neighborhood associated with the material point $\boldsymbol{x}$}\\
	$\theta $ & {volume dilatation value at material point}\\
	$\underline{e}^{d}$  & {deviatoric extension state of peridynamic bond}\\
	$m$ & {weighted volume of material point}\\
	$\underline{\mathit{w}}$ & {influence function}\\
	$\underline{t}$ & {force density scalar state}\\
	$\kappa,\,\mu$ & {bulk modulus and shear modulus of solid material, respectively }\\
	$s$ & {stretch value of peridynamic bond}\\
	$s_{c}$  & {critical stretch of bonds beyond which bond is broken} \\
    $G_{c} $          & {critical energy release rate for mode I fracture} \\ 
    $\varrho $          & {characteristic function describing the connection status of bonds} \\
    $\varphi _{x}$          & {damage value at point $\boldsymbol{x}$} \\
    $s $          & {local normal deformation (stretch) of bond} \\
    $s_c $          & {critical stretch of bonds beyond which bond is broken} \\  
	$\boldsymbol{\sigma}^{tot}, \boldsymbol{\sigma}^{eff}  $ & {total and effective stress tensors} \\  
	$\boldsymbol{I}$ & {unit tensor} \\
	$p$ & {pore pressure}\\  
	$\alpha_{r}, \alpha_{T}, \alpha_{f}$ & {Biot coefficient in reservoir, transition and fracture domains} \\ 
	$\rho_{r}, \rho_{T}, \rho_{f}  $  & {mass density of media in reservoir, transition and fracture domains} \\		
	\hline
\end{tabular}
\newpage
\begin{tabular}[t]{|lp{12.5cm}|}
	\hline       
	$s_{r}, s_{T}, s_{f}  $  & {storage coefficient of media in reservoir, transition and fracture domains} \\	  
	$n_{r}, n_{T}, n_{f}  $  & {porosity of media in reservoir, transition and fracture domains} \\	
	$k_{r}, k_{T}, k_{f}  $  & {permeability of media in reservoir, transition and fracture domains} \\	
	$q_{r}, q_{T}, q_{f}  $  & {source term in reservoir, transition and fracture domains} \\	
	$K_{r}, K_{w} $  & {bulk moduli of solid skeleton and fluid, respectively} \\	
	$\varepsilon _{v}$ &  {volumetric strain}	\\	
	$\mu_{w} $          & {viscosity coefficient of the fluid} \\
	$g$          & {gravity acceleration: $9.8m^{2}/s$} \\	
	$c_{1},c_{2}$ & {threshold values for identifying the three flow domains}\\
	$\chi _{r},\chi _{f}$ & {linear indicator functions defined to connect the three flow domains}\\
	$a$          & {aperture of the crack} \\      
	$\boldsymbol{S}, \boldsymbol{Q}, \boldsymbol{H}$ & {compressibility, coupling and permeability matrices of FE equations}\\        
	$\boldsymbol{N}_{u}, \boldsymbol{N}_{p}$ & {shape functions for displacement and pressure, respectively} \\ 
	$\boldsymbol{m}$ & {unit vector replacing $\boldsymbol{I}$ used in finite element equations}\\       
	$\boldsymbol{M}^{PD}, \boldsymbol{K}^{PD}, \boldsymbol{Q}^{PD}$ & {mass, stiffness and coupling matrices of PD equations}\\   	    
	$E$          & {Young's modulus of solids} \\	
	$\nu$          & {Poisson's ratio of solids} \\
	\hline
\end{tabular}} \clearpage

\section*{Acknowledgements}

This work has been jointly supported by the National Key Research and
Development Program of China (2017YFC1501102), the National Natural Science
Foundation of China (Grant Nos. 51679068 and 11872172), the Fundamental
Research Funds for Central Universities (2017B704X14), the Postgraduate
Research \& Practice Innovation Program of Jiangsu Province (Grant No.
KYCX17\_0479), and the China Scholarship Council (No. 201706710018).

\textcolor{blue}{U. Galvanetto and M. Zaccariotto acknowledge the support they received from MIUR under the research project PRIN2017-DEVISU and from University of Padua under the research projects BIRD2018 NR.183703/18 and BIRD2017 NR.175705/17. }

F. Pesavento would like to acknowledge the project 734370-BESTOFRAC
\textquotedblleft Environmentally best practices and optimisation in
hydraulic fracturing for shale gas/oil development\textquotedblright
-H2020-MSCA-RISE-2016 and the support he received from University of Padua
under the research project BIRD197110/19 \textquotedblleft Innovative models
for the simulation of fracturing phenomena in structural engineering and
geomechanics\textquotedblright .

B.A. Schrefler gratefully acknowledges the support of the Technische
Universität München - Institute for Advanced Study, funded by the German
Excellence Initiative and the TUV SÜD Foundation.

\textcolor{blue}{The authors thank also C. Peruzzo for useful discussions.}

\clearpage

\bibliographystyle{elsarticle-num}
\bibliography{mybib}

\newpage


%
%

\end{document}